\documentclass{article}

\usepackage{amsfonts}
\usepackage{amsmath}
\usepackage{mathrsfs}
\usepackage{amssymb}
\usepackage{authblk}
\usepackage{graphics}
\usepackage{epstopdf}
\usepackage{hyperref}
\usepackage{subcaption}
\usepackage{tikz}
\usepackage{comment}
\usepackage{soul}
\usepackage{color}
\usepackage{upgreek}

\usepackage[square,numbers]{natbib}
\usepackage{geometry}
\newcommand{\bm}[1]{\boldsymbol{#1}}
\newcommand{\wtilde}{\widetilde}
\newcommand{\what}{\widehat}
\newcommand{\dd}{\mathrm{d}}
\newcommand{\pd}{\partial}
\newcommand{\Tr}{\intercal}

\newcommand{\fpd}[2]{\frac{\pd #1}{\pd #2}}

\newcommand{\sGamma}{\mathit{\Gamma}}

\DeclareMathAlphabet{\mathsfit}{T1}{\sfdefault}{\mddefault}{\sldefault}
\SetMathAlphabet{\mathsfit}{bold}{T1}{\sfdefault}{\bfdefault}{\sldefault}

\providecommand{\keywords}[1]{\textbf{Keywords: } #1}

\title{Modelling wave propagation in elastic solids via high-order accurate implicit-mesh discontinuous Galerkin methods}
\author[1,*]{Vincenzo Gulizzi}
\author[2]{Robert Saye}
\affil[1]{Center for Computational Sciences and Engineering (CCSE), Lawrence Berkeley National Laboratory MS 50A-3111, 1 Cyclotron Rd, Berkeley, CA 94720, USA}
\affil[2]{Mathematics Group, Lawrence Berkeley National Laboratory MS 50A-1148, 1 Cyclotron Rd, Berkeley, CA 94720, USA}
\affil[*]{\texttt{vgulizzi@lbl.gov}}

\begin{document}

\maketitle

\begin{abstract}
A high-order accurate implicit-mesh discontinuous Galerkin framework for wave propagation in single-phase and bi-phase solids is presented.
The framework belongs to the embedded-boundary techniques and its novelty regards the spatial discretization, which enables boundary and interface conditions to be enforced with high-order accuracy on curved embedded geometries.
High-order accuracy is achieved via high-order quadrature rules for implicitly-defined domains and boundaries, whilst a cell-merging strategy addresses the presence of small cut cells.
The framework is used to discretize the governing equations of elastodynamics, written using a first-order hyperbolic momentum-strain formulation, and an exact Riemann solver is employed to compute the numerical flux at the interface between dissimilar materials with general anisotropic properties.
The space-discretized equations are then advanced in time using explicit high-order Runge-Kutta algorithms.
Several two- and three-dimensional numerical tests including dynamic adaptive mesh refinement are presented to demonstrate the high-order accuracy and the capability of the method in the elastodynamic analysis of single- and bi-phases solids containing complex geometries.

\end{abstract}

\section{Introduction}\label{sec:INTRO}
The propagation of elastic waves is the subject of many research fields in science and engineering, such as geophysics \cite{chapman2004fundamentals}, structural health monitoring \cite{mitra2016guided} and metamaterials design \cite{srivastava2015elastic,wu2020brief}.
In these areas, computational methods are well-established tools, especially when complicated geometries and/or heterogeneous materials are involved.
However, a common burden in computational modelling is the generation of a high-quality mesh of the domains of analysis, which often results in the most laborious part of the development of a numerical scheme \cite{cottrell2009isogeometric}.

Generally, there are two approaches to meshing irregular domains.
A common approach is to use body-fitted meshes, whereby the mesh elements are generated to conform to the boundaries of the domain;
this approach can be flexible in resolving complex geometrical features but often requires a non-trivial effort to provide high-quality elements and may become more demanding for moving meshes or dynamic adaptive mesh refinements.
On the other hand, one may use embedded-boundary (EB) methods, sometimes also referred to as cut-cell, immersed-boundary or fictitious-domain methods, where a curved geometry is represented on a regular background grid and the boundaries of the grid's cells do not need to conform to the boundaries of the geometry;
this facilitates a more algorithmic or automated approach to mesh generation, data storage and adaptive mesh refinement (AMR), see e.g.~Ref.\cite{aftosmis2000parallel,burman2015cutfem}.
Additionally, in EB methods, the majority of the mesh elements are regular elements and their properties (such as the mass matrix, its inverse and/or the interpolation operator for AMR applications) can be routinely, efficiently and accurately handled.
However, it is clear that EB methods require additional considerations to handle the presence of curved geometries, especially if high-order accuracy is desired.

One of the most widely employed numerical techniques in solid and structural dynamics is the Finite Element Method (FEM) \cite{komatitsch2000simulation,bathe2006finite,ham2012finite}.
FEM models are typically based on body-fitted meshing strategies but have also been developed in combination with EB approaches.
A notable example is the Finite Cell Method (FCM) \cite{parvizian2007finite,duster2008finite,joulaian2014finite,elhaddad2015finite} for two- and three-dimensional problems.
In the FCM, the integrand functions used to evaluate all elemental quantities of a standard FEM (such as the mass matrix, the stiffness matrix or the body forces) are multiplied by an indicator function $\psi(\bm{x})$, which is introduced to define the embedded domain $\mathscr{D}$ as $\psi(\bm{x}) = 1$ for $\bm{x} \in \mathscr{D}$ and $\psi(\bm{x}) = \psi_0$, with $0 < \psi_0 \ll 1$, for $\bm{x} \notin \mathscr{D}$.
This approach simplifies the mesh generation procedure but requires suitable integration schemes, or quadrature rules, for those elements where the indicator function jumps from 1 to $\psi_0$ and vice versa.
In Ref.\cite{abedian2013performance}, the authors compared different integration approaches for the FCM and concluded that the most effective strategy to compute the integrals was to partition the elements into quadtree or octree subgrids (in 2D or 3D, respectively) suitably refined in proximity of the embedded boundaries.
Indeed, the accurate evaluation of the elemental matrices in EB methods is not a simple task, particularly in 3D.
The use of high-order quadrature rules for cut cells is one of the key aspects of the method presented in this work, as discussed shortly.

Various modifications of and/or alternatives to the FEM have been proposed in the literature to reduce the meshing effort.
A few examples are the extended FEM \cite{motamedi2010dynamic,chin2021spectral} and the extended Ritz method \cite{benedetti2019x}, whereby the space of basis functions is enriched to automatically account for the presence of embedded boundaries or interfaces, the Virtual Element Method (VEM) \cite{park2019nonconvex}, which allows the use of general (also non-convex) polygonal elements, or the Boundary Element Method (BEM) \cite{manolis1983comparative,benedetti2010fast}, which is based on an integral formulation and allows solving the equations of elastodynamics by discretizing the domain's boundaries only.

The discontinuous Galerkin (DG) method has also proven to be a very powerful and flexible numerical technique for solving different classes of PDEs \cite{arnold2002unified,cockburn1998runge}.
In a DG method, the numerical solution is represented in a space of polynomial basis functions that are discontinuous among the mesh elements.
Then, the inter-element continuity, the interface conditions between dissimilar materials and the boundary conditions are enforced in weak sense by introducing suitable boundary integrals.
This naturally enables high-order accuracy, the treatment of generally-shaped elements and $hp$ AMR with conventional and non-conventional (e.g.~polytopic) mesh elements \cite{hartmann2002adaptive,zanotti2015solving,cangiani2017hp,antonietti2019v}.
Moreover, unlike other numerical schemes based on continuous approximations, DG methods feature block-structured mass matrices, which are highly desirable in explicit time-stepping schemes and parallel computations as they can be easily inverted on an element-by-element basis.
Thus, thanks to its discontinuous nature, the DG method has been employed in combination with the EB approach for high-order accurate solution of
elliptic PDEs \cite{johansson2013high,saye2017implicitI},
incompressible and compressible fluid flow \cite{saye2017implicitII,gulizzi2021coupled}, and
statics of thin structures with cut-outs \cite{gulizzi2020implicit}, among many other applications.
However, in the context of solving elastodynamics problems, while the literature offers several DG schemes using body-fitted meshes, which can be classified into those derived from the second-order hyperbolic formulations \cite{antonietti2016high,appelo2018energy,antonietti2018high} and those derived from the first-order hyperbolic formulations \cite{dumbser2006arbitrary,de2007arbitrary,wilcox2010high,zhan2018exact}, DG methods using embedded boundaries appear less investigated.
Recently, Tavelli et al.\cite{tavelli2019simple} proposed a diffused-interface DG scheme for elastic waves where, similar to the FCM, the curved geometry is represented by an indicator function that takes value $1$ within the solid and $0$ outside the solid; as discussed by the authors, their scheme is high-order accurate far from the embedded boundary but only first-order accurate in the transition region of the indicator function, i.e.~in proximity of the embedded boundary.
Although high-order accurate embedded-boundary methods for the 2D scalar wave equation \cite{sticko2019higher} and the 2D acoustic equation \cite{adjerid2019immersed} have been developed, to the best of the authors' knowledge, a high-order accurate embedded-boundary method for elastic wave propagation in two- and three-dimensional geometries has not before been investigated in the literature.

In this work, we present high-order accurate DG methods for elastodynamics in embedded geometries.
In prior work, these kinds of DG methods have been successfully employed to model with high-order accuracy free-surface flow and rigid body-fluid interaction in the incompressible regime \cite{saye2016interfacial,saye2017implicitI,saye2017implicitII,saye2020fast}, the static response of thin multilayered structures with cut-outs \cite{gulizzi2020implicit,gulizzi2020high,guarino2021high}, gas dynamics problems \cite{gulizzi2021coupled}; here, the methods are extended to model elastic waves propagating in single- and bi-phase solids characterized by general anisotropic behavior.
Geometries are represented via a level set function, whose zero-contour denotes either the curved boundary of a single-phase solid or the interface between the phases of a two-phase solid.
A key feature of the proposed approach is the use of high-order quadrature rules for implicitly-defined domains and boundaries stemming from the intersection between the grids and the level set function.
These quadrature rules are generated using the algorithm developed in Ref.\cite{saye2015high} and enable the resolution of the embedded geometry as well as the enforcement of boundary and interface conditions with high-order accuracy.
To avoid the presence of arbitrarily small cut cells, which would lead to overly restrictive times steps and ill-conditioned discrete operators, the framework developed in this work uses a cell-merging technique, whereby those cells with a volume fraction smaller than a user-defined threshold are merged with their neighbors.
The DG method is used to discretize the governing equations of elastodynamics in space and an exact Riemann solver \cite{zhan2018exact} is employed to compute the numerical flux at the interface between materials with dissimilar elastic properties.
Meanwhile, explicit high-order Runge-Kutta algorithms \cite{cockburn1998runge} serve as time integrators.

The paper is organized as follows:
Sec.(\ref{sec:IBVP}) introduces the momentum-strain formulation of elastodynamics for single- and bi-phase solids considered in this work;
Sec.(\ref{sec:DG}) presents the implicit-mesh discontinuous Galerkin framework including the generation of the implicitly-defined meshes, the weak form of the governing equations and the $hp$ adaptive mesh refinement;
Sec.(\ref{sec:RESULTS}) demonstrates high-order accuracy and the capability of the method by discussing numerical results obtained for several two- and three-dimensional wave propagation problems in single- and bi-phase solids.
Conclusions and discussions for further developments are given in Sec.(\ref{sec:CONCLUSIONS}).

\section{Elastodynamic formulation}\label{sec:IBVP}

\subsection{Geometry description}

\begin{figure}
\centering
\includegraphics[width = \textwidth]{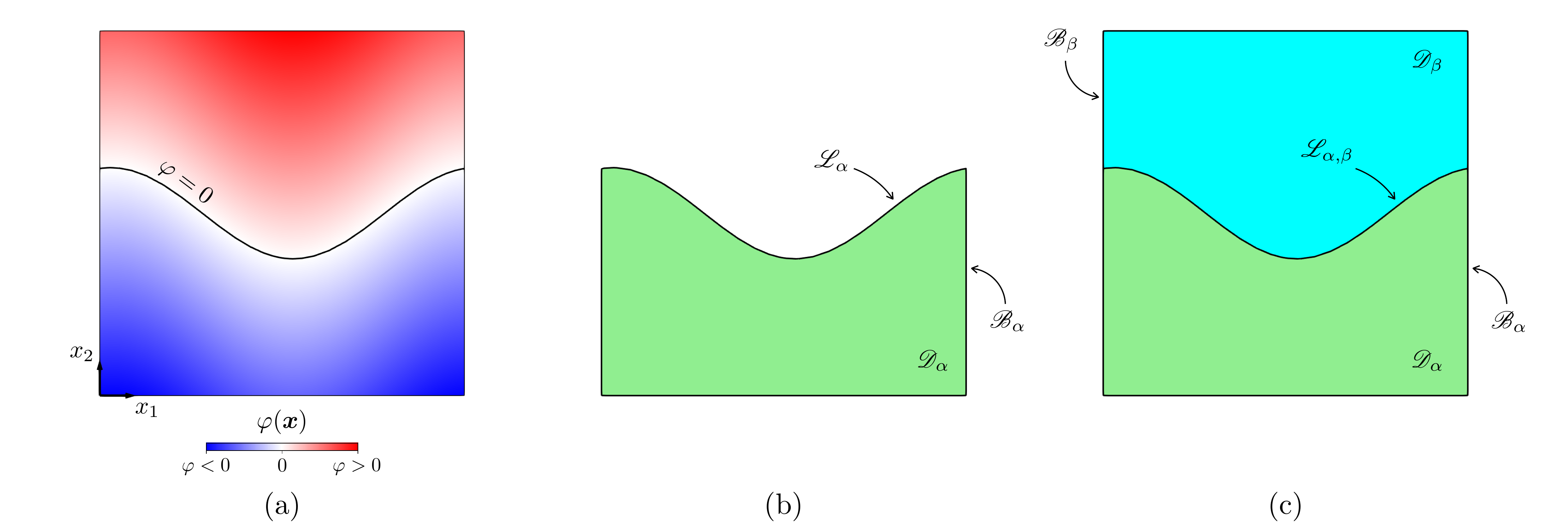}
\caption
{
    (a) Sample level set function $\varphi$ partitioning a square into a region where $\varphi$ is negative and a region where $\varphi$ is positive.
    (b) Single-phase solid and (c) bi-phase solid implicitly defined using the level set function of figure (a).
}
\label{fig:IBVP - LEVEL SET GEOMETRY}
\end{figure}

The geometry of a single-phase or a bi-phase solid is described implicitly by a level set function.
Consider a $d$-dimensional rectangle $\mathscr{R} \subset \mathbb{R}^d$ and its outer boundary $\pd\mathscr{R}$.
Consider also a level set function $\varphi: \mathscr{R} \rightarrow \mathbb{R}$ and let $\mathscr{R}^- \equiv \{\bm{x} \in \mathscr{R} : \varphi(\bm{x}) < 0\}$ be the portion of $\mathscr{R}$ where $\varphi$ is negative, $\mathscr{R}^+ \equiv \{\bm{x} \in \mathscr{R} : \varphi(\bm{x}) > 0\}$ be the portion of $\mathscr{R}$ where $\varphi$ is positive and $\mathscr{L} \equiv \{\bm{x}\in\mathscr{R} : \varphi(\bm{x}) = 0\}$ be the zero-contour of $\varphi$.
Moreover, let $\pd\mathscr{R}^- \equiv \{\bm{x} \in \pd\mathscr{R} : \varphi(\bm{x}) < 0\}$ be the portion of $\pd\mathscr{R}$ where $\varphi$ is negative and $\pd\mathscr{R}^+ \equiv \{\bm{x} \in \pd\mathscr{R} : \varphi(\bm{x}) > 0\}$ be the portion of $\pd\mathscr{R}$ where $\varphi$ is positive.

Then, we define the domain $\mathscr{D}_\alpha$ and the boundary $\pd\mathscr{D}_\alpha$ of a solid consisting of one phase $\alpha$ as $\mathscr{D}_\alpha \equiv \mathscr{R}^-$ and $\pd\mathscr{D}_\alpha \equiv \mathscr{B}_\alpha\cup\mathscr{L}_\alpha$, respectively, where $\mathscr{B}_\alpha \equiv \pd\mathscr{R}^-$ and $\mathscr{L}_\alpha \equiv \mathscr{L}$.
Alternatively, for a bi-phase solid consisting of two distinct phases $\alpha$ and $\beta$, we define the domain and the outer boundary of the phase $\alpha$ as $\mathscr{D}_\alpha \equiv \mathscr{R}^-$ and $\mathscr{B}_\alpha \equiv \pd\mathscr{R}^-$, respectively, the domain and the outer boundary of the phase $\beta$ as $\mathscr{D}_\beta \equiv \mathscr{R}^+$ and $\mathscr{B}_\beta \equiv \pd\mathscr{R}^+$, respectively, and the interface between the two phases as $\mathscr{L}_{\alpha,\beta} \equiv \mathscr{L}$.

As an example, Fig.(\ref{fig:IBVP - LEVEL SET GEOMETRY}a) shows a level set function $\varphi$ defined over a square, while Fig.(\ref{fig:IBVP - LEVEL SET GEOMETRY}b) shows the domain and the boundaries of the corresponding two-dimensional single-phase solid, and Fig.(\ref{fig:IBVP - LEVEL SET GEOMETRY}c) shows the domains, boundaries and interface of the corresponding two-dimensional bi-phase solid.

\subsection{Governing equations}\label{ssec:Governing equations}
The governing equations of linear elastodynamics can be stated using different formulations, which, on the basis of the chosen set of primary variables, include the displacement formulation, see e.g.~\cite{bathe2006finite,ham2012finite}, the velocity-stress formulation, see e.g.~\cite{de2007arbitrary}, or the momentum-strain formulation, see e.g.~\cite{zhan2018exact}.
Here, we use the momentum-strain formulation since it allows for the consideration of materials with space-dependent constitutive properties.
In the following, the governing equations, the initial conditions and the boundary conditions, are written for the domain $\mathscr{D}_\alpha$ and its boundary $\pd\mathscr{D}_\alpha$, but are also valid for the domain $\mathscr{D}_\beta$ and its boundary $\pd\mathscr{D}_\beta$, if one replaces the subscript $\alpha$ with $\beta$.
However, note that both subscripts $\alpha$ and $\beta$ appear explicitly in the equations governing the interface conditions for bi-phase solids discussed at the end of this section.

Let $\bm{v}_\alpha$ and $\bm{m}_\alpha$ denote the velocity field and the momentum field, respectively, and let $\bm{\gamma}_\alpha$ and $\bm{\sigma}_\alpha$ denote the strain field and the stress field, respectively.
In $\mathbb{R}^d$, upon defining $N_v \equiv d$ and $N_\sigma \equiv d(d+1)/2$, $\bm{v}_\alpha$ and $\bm{m}_\alpha$ are $N_v$-dimensional vectors, whereas $\bm{\gamma}_\alpha$ and $\bm{\sigma}_\alpha$ are $N_\sigma$-dimensional vectors containing the strain and the stress components, respectively, in Voigt notation \cite{carcione2007wave}; for example, in 2D we have $\bm{\gamma}_\alpha = (\gamma_{\alpha 11},\gamma_{\alpha 22},\gamma_{\alpha 12})^\Tr$, while in 3D we have $\bm{\gamma} = (\gamma_{\alpha 11},\gamma_{\alpha 22},\gamma_{\alpha 33},\gamma_{\alpha 23},\gamma_{\alpha 13},\gamma_{\alpha 12})^\Tr$, where $\bullet^\Tr$ denotes the transpose of $\bullet$.

The domain $\mathscr{D}_\alpha$ is characterized by the density $\rho_\alpha$, such that $\bm{m}_\alpha = \rho_\alpha\bm{v}_\alpha$, and by the (positive-definite) $N_\sigma\times N_\sigma$ matrix $\bm{c}_\alpha$ of elastic stiffness constants that links the stress and strain via the general Hooke's law $\bm{\sigma}_\alpha = \bm{c}_\alpha \bm{\gamma}_\alpha$.
Then, in the momentum-strain formulation, the governing equations of elastodynamics are written as the following first-order hyperbolic system of PDEs
\begin{equation}\label{eq:governing equations - vector form}
\fpd{\bm{U}_\alpha}{t}+\fpd{\bm{F}_{\alpha i}}{x_i} = \bm{S}_{\alpha}.
\end{equation}
In Eq.(\ref{eq:governing equations - vector form}) and in the remainder of the paper, the subscript $i$ takes values in $\{1,\dots,d\}$ and implies summation when repeated, $t$ is the time, and $x_i$ is the $i$-th coordinate of the $d$-dimensional space location vector $\bm{x}$.
Additionally, $\bm{U}_{\alpha}$, $\bm{F}_{\alpha i}$ and $\bm{S}_{\alpha}$ are $N_U$-dimensional vectors, with $N_U \equiv N_v+N_\sigma$, representing the conserved variables, the fluxes in the $i$ direction and the source terms, respectively; they are defined as
\begin{equation}\label{eq:governing equations - vectors}
\bm{U}_\alpha\equiv
\left(\begin{array}{c}
\bm{m}_\alpha\\
\bm{\gamma}_\alpha
\end{array}
\right),
\quad
\bm{F}_{\alpha i}\equiv
\left(\begin{array}{c}
-\bm{I}_i^\Tr\bm{c}_\alpha\bm{\gamma}_\alpha\\
-\rho_\alpha^{-1}\bm{I}_i\bm{m}_\alpha\\
\end{array}
\right)
\quad\mathrm{and}\quad
\bm{S}_\alpha\equiv
\left(\begin{array}{c}
\rho_\alpha\bm{b}_\alpha\\
\bm{0}\\
\end{array}
\right),
\end{equation}
where $\bm{b}$ is the $d$-dimensional vector of body forces and the matrices $\bm{I}_i$ are given by
\begin{equation}
\bm{I}_1\equiv
\left[\begin{array}{ccc}
1&0\\
0&0\\
0&1
\end{array}
\right]
\quad\mathrm{and}\quad
\bm{I}_2\equiv
\left[\begin{array}{ccc}
0&0\\
0&1\\
1&0
\end{array}
\right],
\end{equation}
in 2D, and by
\begin{equation}
\bm{I}_1\equiv
\left[\begin{array}{ccc}
1&0&0\\
0&0&0\\
0&0&0\\
0&0&0\\
0&0&1\\
0&1&0
\end{array}
\right],
\quad
\bm{I}_2\equiv
\left[\begin{array}{ccc}
0&0&0\\
0&1&0\\
0&0&0\\
0&0&1\\
0&0&0\\
1&0&0
\end{array}
\right]
\quad\mathrm{and}\quad
\bm{I}_3\equiv
\left[\begin{array}{ccc}
0&0&0\\
0&0&0\\
0&0&1\\
0&1&0\\
1&0&0\\
0&0&0
\end{array}
\right],
\end{equation}
in 3D.

Eq.(\ref{eq:governing equations - vector form}) is assumed to be valid for $(t,\bm{x}) \in \mathscr{T}\times\mathscr{D}_\alpha$, where $\mathscr{T} \equiv [0,T]$ is the time interval and $T$ is the final time, and is supplemented by initial, boundary and interface conditions.
Initial conditions are given as
\begin{equation}\label{eq:initial conditions}
\bm{U}_\alpha = \bm{U}_{\alpha 0}(\bm{x})\quad\mathrm{for}~t = 0~\mathrm{and}~\bm{x}\in\mathscr{D}_\alpha,
\end{equation}
where $\bm{U}_{\alpha 0}(\bm{x})$ contains the known values of $\bm{U}_\alpha$ at $t = 0$.

Boundary conditions are prescribed at the outer boundary of the domain $\mathscr{D}_\alpha$.
Here, we consider three types of boundary conditions, namely \emph{i}) prescribed values of the velocity field $\bm{v}$, \emph{ii}) prescribed values of the traction field $\bm{t}$ and \emph{iii}) absorbing boundary conditions, which are typical for elastodynamics problems.
Prescribing the value of the velocity field $\bm{v}$ is equivalent to prescribing the momentum field $\bm{m}$ as 
\begin{equation}\label{eq:boundary conditions - v}
\bm{m}_\alpha = \rho_\alpha\overline{\bm{v}}\quad\mathrm{for}~(t,\bm{x})\in\mathscr{T}\times\pd\mathscr{D}_{\alpha v},
\end{equation}
where $\overline{\bm{v}}$ is the prescribed value of the velocity field, which includes $\overline{\bm{v}} = \bm{0}$ in case of a fixed boundary, and $\pd\mathscr{D}_{\alpha v}$ is the portion of the outer boundary of $\mathscr{D}_\alpha$ where $\overline{\bm{v}}$ is prescribed.
Prescribing the value of the traction field $\bm{t}$ is expressed in terms of the strain $\bm{\gamma}_\alpha$ as 
\begin{equation}\label{eq:boundary conditions - t}
\bm{I}_{n_\alpha}^\Tr\bm{c}_\alpha\bm{\gamma}_\alpha = \overline{\bm{t}}\quad\mathrm{for}~(t,\bm{x})\in\mathscr{T}\times\pd\mathscr{D}_{\alpha t},
\end{equation}
where $\overline{\bm{t}}$ is the prescribed value of the traction field, which includes $\overline{\bm{t}} = \bm{0}$ in case of a free boundary, $\pd\mathscr{D}_{\alpha t}$ is the portion of the outer boundary of $\mathscr{D}_\alpha$ where $\overline{\bm{t}}$ is prescribed and $\bm{I}_{n_\alpha} \equiv n_{\alpha i}\bm{I}_i$, being $n_{\alpha i}$ the $i$-th component of the outer unit normal of $\pd\mathscr{D}_{\alpha t}$.
Absorbing boundary conditions refer to a boundary that does not reflect the incoming waves.

Finally, for bi-phase solids only, interface conditions are prescribed at the interface $\mathscr{L}_{\alpha,\beta}$.
Here, we consider perfect interface conditions, i.e.~continuity of the velocity field and equilibrium of the traction field, which are given in terms of momentum and strain as follows
\begin{equation}\label{eq:interface conditions}
\left\{\begin{array}{l}
\rho_\alpha^{-1}\bm{m}_\alpha = \rho_\beta^{-1}\bm{m}_\beta\\
\bm{I}_{n_\alpha}^\Tr\bm{c}_\alpha\bm{\gamma}_\alpha+\bm{I}_{n_\beta}^\Tr\bm{c}_\beta\bm{\gamma}_\beta = \bm{0}
\end{array}\right.
\quad\mathrm{for}~(t,\bm{x})\in\mathscr{T}\times\mathscr{L}_{\alpha,\beta}.
\end{equation}

Within the present DG framework, the three types of boundary conditions (including the absorbing boundary conditions) and the interface conditions are enforced via suitable definitions of the numerical flux as proposed by Zhan et al.\cite{zhan2018exact}.
The numerical flux is a key ingredient of DG formulations and will be introduced in Sec.(\ref{sec:DG}).
\section{Implicit-mesh discontinuous Galerkin methods}\label{sec:DG}

\subsection{Implicitly-defined meshes}\label{ssec:IM}

The discontinuous Galerkin method typically requires a suitable mesh of the domain under analysis.
Here, we use the implicitly-defined mesh technique \cite{saye2017implicitI, saye2017implicitII, gulizzi2021coupled}, whereby the domain discretization is obtained by intersecting the implicitly-defined phases and a structured background grid that is easily generated for the rectangle $\mathscr{R}$ containing the solid.

Consider a structured grid $\mathscr{G} \equiv \bigcup_{\bm{j}}\mathscr{C}^{\bm{j}}$, where $\mathscr{C}^{\bm{j}} \equiv [x_1^{\bm{j}},x_1^{\bm{j}}+h_1]\times\dots\times [x_d^{\bm{j}},x_d^{\bm{j}}+h_d] \subset \mathscr{R}$ is a $d$-dimensional rectangular cell, $\bm{j}$ is the $d$-tuple identifying the location of the cell within the grid and $x_i^{\bm{j}}$ and $h_i$ are the cell's lower end and the cell's size in the $i$-th direction, respectively.
Each cell $\mathscr{C}$ is intersected with the implicitly-defined phases of the solid and is classified on the basis of its volume fraction $f_\alpha$ given by
\begin{equation}\label{eq:volume fraction}
f_\alpha \equiv \frac{1}{V_{\mathscr{C}}}\int\limits_{\mathscr{C}\cap\mathscr{D}_\alpha}1~\dd V,
\end{equation}
where $V_\mathscr{C}$ is the volume of the cell, that is $V_\mathscr{C} = h_1h_2$ in 2D and $V_\mathscr{C} = h_1h_2h_3$ in 3D.
Referring to the phase $\alpha$, 
\emph{entire} cells are those cells falling entirely inside $\mathscr{D}_\alpha$ and having volume fraction $f_\alpha = 1$;
\emph{empty} cells are those cells falling entirely outside $\mathscr{D}_\alpha$ and having volume fraction $f_\alpha = 0$;
\emph{large} cells are those cells cut by $\mathscr{L}$ and having volume fraction $\overline{f} < f_\alpha < 1$;
\emph{small} cells are those cells cut by $\mathscr{L}$ and having volume fraction $0 < f_\alpha \le \overline{f}$.
The same classification is performed by intersecting the cells with the phase $\beta$ in bi-phase solids.
Henceforth, entire and large cells are collectively referred to as \emph{primary} cells.

In the classification above, the parameter $\overline{f}$ denotes a user-defined volume fraction threshold introduced to identify the small cells, i.e.~those cells whose presence would lead to overly-small time-step restrictions and ill-conditioned discrete operators.
Here, to address the small-cell problem, we employ a cell-merging strategy whereby small cells are merged with their neighbors.
In particular, each small cell is merged with one primary cell among the neighboring cells of its $3\times 3$ neighborhood, in 2D, or $3\times 3\times 3$ neighborhood, in 3D.
To select the neighbor for merging, the neighboring cells are grouped in the following order: in 3D, we consider first the cells sharing a face with the small cell, second the cells sharing an edge with the small cell and third the cells sharing a corner with the small cell; within each group the neighboring cells are ordered according to their volume fraction.
Then, the neighbor targeted for merging is the first cell in the first non-empty group.
In 2D, the search for the target neighbor starts from the neighboring cells that share an edge with the small cell.
It is worth noting that that, although each small cell is merged with one primary cell only, multiple small cells are allowed to target the same primary cell, possibly leading to a cluster of multiple cells merged together.

Once the cell-merging procedure is completed, the cells of the grid consist of \emph{non-merged} cells, i.e.~the primary cells that have not been targeted during the merging process, \emph{merged} cells, i.e.~the union of small cells and their merging neighbors, and empty cells.
Then, the implicitly-defined mesh $\mathcal{M}_\alpha$ of the phase $\alpha$ is written $\mathcal{M}_\alpha \equiv \bigcup_{e = 1}^{N_\alpha^e} \mathscr{D}_\alpha^{e}$, where $\mathscr{D}_\alpha^{e}$ is the $e$-th implicitly-defined mesh element of the phase $\alpha$ and $N_\alpha^e$ is the number of mesh elements.
Finally, the implicitly-defined mesh $\mathcal{M}$ of a single-phase solid coincides with the implicitly-defined mesh of the phase $\alpha$, i.e.~$\mathcal{M} \equiv \mathcal{M}_\alpha$, whereas, the implicitly-defined mesh $\mathcal{M}$ of a bi-phase solid is the collection of mesh elements of the phase $\alpha$ and the phase $\beta$, i.e.~$\mathcal{M} \equiv \mathcal{M}_\alpha\bigcup\mathcal{M}_\beta$.
Ultimately, these implicitly-defined meshes consist of a collection of standard $d$-dimensional rectangular elements and a relatively smaller number of curved elements that conform with the curvature of the zero-contour of the level set function.

\begin{figure}
\centering
\includegraphics[width = \textwidth]{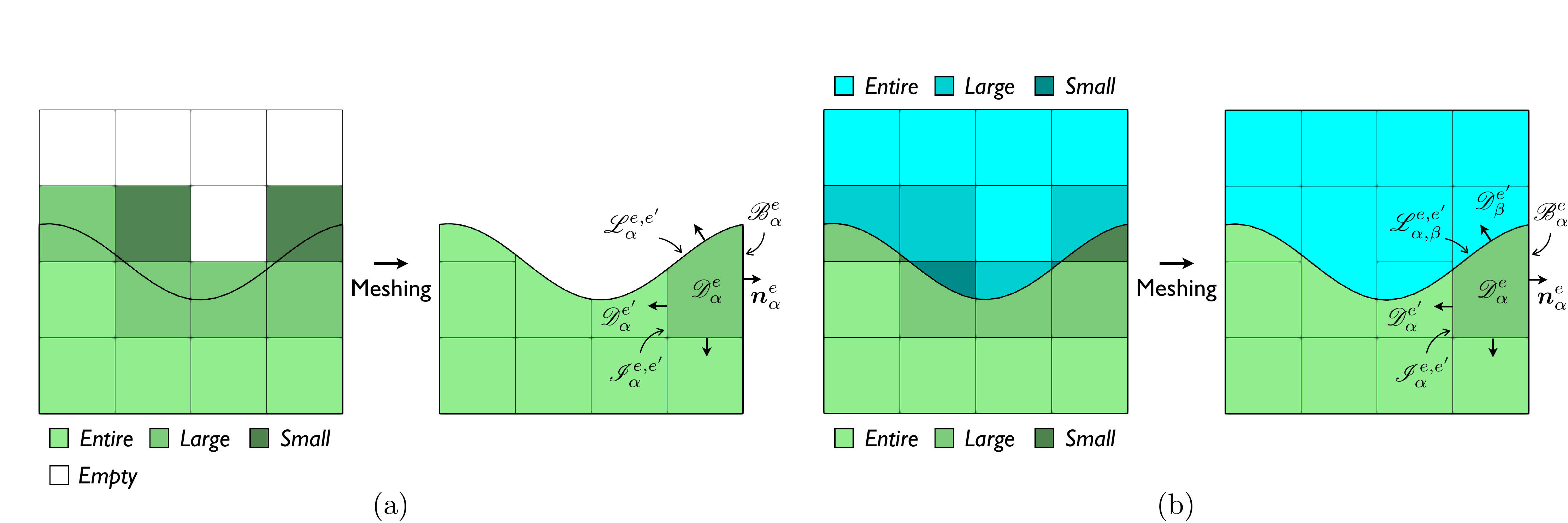}
\caption
{
    Cell classification and implicitly-defined meshes for (a) the single-phase solid and (b) the bi-phase solid shown in Fig.(\ref{fig:IBVP - LEVEL SET GEOMETRY}) intersected with a $4\times 4$ background grid.
}
\label{fig:IBVP - GRID MESH}
\end{figure}

\begin{figure}
\centering
\includegraphics[width = \textwidth]{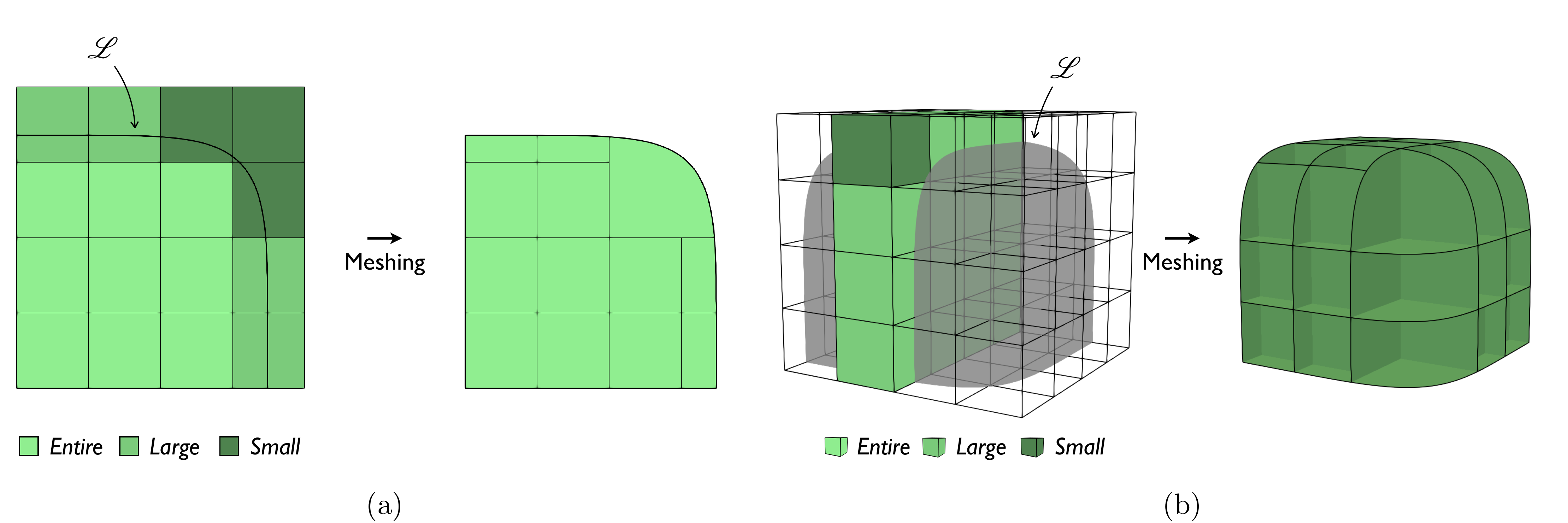}
\caption
{
    Examples of (a) a two-dimensional implicitly-defined mesh and (b) a three-dimensional implicitly-defined mesh where multiple small cells are merged with the same cell.
}
\label{fig:IBVP - MERGING DETAIL}
\end{figure}

Figures (\ref{fig:IBVP - GRID MESH}) and (\ref{fig:IBVP - MERGING DETAIL}) show a few implicitly-defined meshes obtained with the procedure described above.
Figure (\ref{fig:IBVP - GRID MESH}a) shows the classification of the cells of a $4\times 4$ grid when intersected with the single-phase solid of Fig.(\ref{fig:IBVP - LEVEL SET GEOMETRY}b) and the corresponding implicitly-defined mesh; Fig.(\ref{fig:IBVP - GRID MESH}a) also highlights an implicitly-defined element $\mathscr{D}_\alpha^e$ and its outer boundary $\mathscr{B}_\alpha^e\cup\mathscr{L}_\alpha^e$, intraphase boundary $\mathscr{I}_\alpha^{e,e'}$ shared with the neighboring element $\mathscr{D}_\alpha^{e'}$, and outer unit normal $\bm{n}^e_\alpha$.
Similarly, Fig.(\ref{fig:IBVP - GRID MESH}b) shows the classification of the cells of a $4\times 4$ grid when intersected with the bi-phase solid of Fig.(\ref{fig:IBVP - LEVEL SET GEOMETRY}c) and the corresponding implicitly-defined mesh; Fig.(\ref{fig:IBVP - GRID MESH}b) also highlights an implicitly-defined element $\mathscr{D}_\alpha^e$ and its outer boundary $\mathscr{B}_\alpha^e$, intraphase boundary $\mathscr{I}_\alpha^{e,e'}$ shared with the neighboring element $\mathscr{D}_\alpha^{e'}$ of the same phase, interface boundary $\mathscr{L}_{\alpha,\beta}^{e,e'}$ shared with the neighboring element $\mathscr{D}_\beta^{e'}$ of the phase $\beta$, and outer unit normal $\bm{n}^e_\alpha$.
Finally, Fig.(\ref{fig:IBVP - MERGING DETAIL}) shows a 2D example and a 3D example of implicitly-defined meshes where multiple small cells have targeted the same nearby cell for merging.

\subsection{Discontinuous Galerkin formulation}
Once the mesh of the domain is generated, discontinuous Galerkin formulations are derived by introducing of a suitable space of discontinuous basis functions and stating the governing equations in weak form.
Here, given that the discretization is constructed using Cartesian grids and the majority of the mesh elements are standard $d$-dimensional rectangles, it is natural to define the local basis functions as tensor-product polynomials.

Let $\mathscr{D}^{e}$ be a mesh element of $\mathcal{M}$ associated with the primary grid cell $\mathscr{C}^{\bm{j}}$ (and all the small cells that are merged with $\mathscr{C}^{\bm{j}}$) and let $\mathcal{P}_{hp}^e$ be the space of tensor-product polynomials of degree $p$ in the (hyper)rectangular volume occupied by $\mathscr{C}^{\bm{j}}$.
Then, the space $\mathcal{V}_{hp}$ of discontinuous basis functions for the mesh $\mathcal{M}$ based on the grid $\mathscr{G}$ is
\begin{equation}\label{eq:space}
\mathcal{V}_{hp} \equiv \left\{v:\mathscr{G}\rightarrow\mathbb{R}~|~v|_{\mathscr{D}^e}\in\mathcal{P}_{hp}^e,~\forall \mathscr{D}^e\in\mathcal{M}\right\},
\end{equation}
while the related space $\mathcal{V}_{hp}^{N}$ of discontinuous polynomials vector fields is $\mathcal{V}_{hp}^N \equiv (\mathcal{V}_{hp})^N$.

For a single-phase solid, the weak form of the governing equations is obtained by multiplying Eq.(\ref{eq:governing equations - vector form}) by the test functions $\bm{V}\in\mathcal{V}_{hp}^{N_U}$, integrating over a generic mesh element $\mathscr{D}_\alpha^e$ and performing integration by parts in space, which yield
\begin{equation}\label{eq:governing equations - DG weak form - single-phase solid}
\int_{\mathscr{D}_\alpha^e}\bm{V}^\Tr\frac{\pd\bm{U}_\alpha}{\pd t}\dd{V} = 
\int_{\mathscr{D}_\alpha^e}\bm{V}^\Tr\bm{S}_\alpha\dd{V}+
\int_{\mathscr{D}_\alpha^e}\frac{\pd\bm{V}^\Tr}{\pd x_i}\bm{F}_{\alpha i}\dd{V}-
\int_{\mathscr{B}_\alpha^{e}\cup\mathscr{L}_\alpha^e}\bm{V}^\Tr\what{\bm{F}}_n\dd{S}-
\sum_{e'\in\mathcal{N}_\alpha^{e}}\int_{\mathscr{I}_{\alpha}^{e,e'}}\bm{V}^\Tr\what{\bm{F}}_n\dd{S}.
\end{equation}
Similarly, for a bi-phase solid, one obtains
\begin{multline}\label{eq:governing equations - DG weak form - bi-phase solid}
\int_{\mathscr{D}_\alpha^e}\bm{V}^\Tr\frac{\pd\bm{U}_\alpha}{\pd t}\dd{V} =
\int_{\mathscr{D}_\alpha^e}\bm{V}^\Tr\bm{S}_\alpha\dd{V}+
\int_{\mathscr{D}_\alpha^e}\frac{\pd\bm{V}^\Tr}{\pd x_i}\bm{F}_{\alpha i}\dd{V}+\\-
\int_{\mathscr{B}_\alpha^{e}}\bm{V}^\Tr\what{\bm{F}}_n\dd{S}-
\sum_{e'\in\mathcal{N}_\alpha^{e}}\int_{\mathscr{I}_{\alpha}^{e,e'}}\bm{V}^\Tr\what{\bm{F}}_n\dd{S}-
\sum_{e'\in\mathcal{N}_{\alpha,\beta}^{e}}\int_{\mathscr{L}_{\alpha,\beta}^{e,e'}}\bm{V}^\Tr\what{\bm{F}}_n\dd{S}.
\end{multline}
Note that Eq.(\ref{eq:governing equations - DG weak form - bi-phase solid}) is valid for a mesh element $\mathscr{D}_\beta^e$ of $\mathcal{M}_\beta$ if one switches the subscripts $\alpha$ and $\beta$ and considers that $\mathscr{L}_{\beta,\alpha}^{e,e'}$ coincides with $\mathscr{L}_{\alpha,\beta}^{e,e'}$ but has opposite unit normal.
In Eqs.(\ref{eq:governing equations - DG weak form - single-phase solid}) and (\ref{eq:governing equations - DG weak form - bi-phase solid}),
$\mathcal{N}_\alpha^{e}$ denotes the set of mesh elements of $\mathcal{M}_\alpha$ that are neighbors of $\mathscr{D}_\alpha^e$,
$\mathcal{N}_{\alpha,\beta}^{e}$ denotes the set of mesh elements of $\mathcal{M}_\beta$ that are neighbors of $\mathscr{D}_\alpha^e$, and
$\what{\bm{F}}_n$ is the so-called \emph{numerical flux}.
The expression of $\what{\bm{F}}_n$ depends on:
the solution state $\bm{U}_\alpha^e$ and the boundary conditions at $\mathscr{B}_\alpha^{e}\cup\mathscr{L}_\alpha^e$,
the adjacent solution states $\bm{U}_\alpha^{e}$ and $\bm{U}_\alpha^{e'}$ of neighboring mesh elements of the same phase at $\mathscr{I}_{\alpha}^{e,e'}$,
and the adjacent solution states $\bm{U}_\alpha^{e}$ and $\bm{U}_\beta^{e'}$ of neighboring mesh elements of different phases at $\mathscr{L}_{\alpha,\beta}^{e,e'}$.
In all cases mentioned above, the expression of the numerical flux used in this work is based on the exact Riemann solver developed by Zhan et al.\cite{zhan2018exact} for anisotropic elastodynamics wherein adjacent elements are allowed to have different values of density and/or elastic stiffness constants.

The final semidiscrete evolution equation is obtained by expressing $\bm{U}$ over each mesh element as a linear combination of the spatial basis functions with time-dependent coefficients.
Using a compact notation, this is written as
\begin{equation}\label{eq:ansatz}
\bm{U}(t,\bm{x}) = \mathbb{B}_\alpha^e(\bm{x})\mathbb{X}_\alpha^e(t),\quad\mathrm{for}~\bm{x}\in\mathscr{D}_\alpha^e,
\end{equation}
where $\mathbb{B}_\alpha^e$ is a $N_U\times N_UN_p$ matrix containing the basis functions, $\mathbb{X}_\alpha^e$ is the vector of coefficients of length $N_UN_p$ and $N_p \equiv (1+p)^d$ is the number of basis functions.
For a bi-phase solid, substituting Eq.(\ref{eq:ansatz}) into Eq.(\ref{eq:governing equations - DG weak form - bi-phase solid}) and letting $\bm{V}$ range over the basis functions, one obtains
\begin{equation}\label{eq:governing equations - semidiscrete form - bi-phase solid}
\mathbb{M}_\alpha^e\dot{\mathbb{X}}_\alpha^e = 
\int_{\mathscr{D}_\alpha^e}\mathbb{B}_\alpha^{e\Tr}\bm{S}_\alpha\dd{V}+
\int_{\mathscr{D}_\alpha^e}\frac{\pd\mathbb{B}_\alpha^{e\Tr}}{\pd x_i}\bm{F}_{\alpha i}\dd{V}-
\int_{\mathscr{B}_\alpha^{e}}\mathbb{B}_\alpha^{e\Tr}\what{\bm{F}}_n\dd{S}-
\sum_{e'\in\mathcal{N}_\alpha^{e}}\int_{\mathscr{I}_{\alpha}^{e,e'}}\mathbb{B}_\alpha^{e\Tr}\what{\bm{F}}_n\dd{S}-
\sum_{e'\in\mathcal{N}_{\alpha,\beta}^{e}}\int_{\mathscr{L}_{\alpha,\beta}^{e,e'}}\mathbb{B}_\alpha^{e\Tr}\what{\bm{F}}_n\dd{S},
\end{equation}
where the superimposed dot denotes the time derivative and $\mathbb{M}_\alpha^e$ is the mass matrix of the element $\mathscr{D}_\alpha^e$ given by
\begin{equation}\label{eq:governing equations - mass matrix}
\mathbb{M}_\alpha^e \equiv 
\int_{\mathscr{D}_\alpha^e}\mathbb{B}_\alpha^{e\Tr}\mathbb{B}_\alpha^{e}\dd{V}.
\end{equation}
An expression similar to the one given in Eq.(\ref{eq:governing equations - semidiscrete form - bi-phase solid}) is obtained for single-phase solids if one substitutes Eq.(\ref{eq:ansatz}) into Eq.(\ref{eq:governing equations - DG weak form - single-phase solid}).

\subsection{High-order quadrature rules for implicitly-defined elements}\label{ssec:Q}

\begin{figure}
\centering
\includegraphics[width = \textwidth]{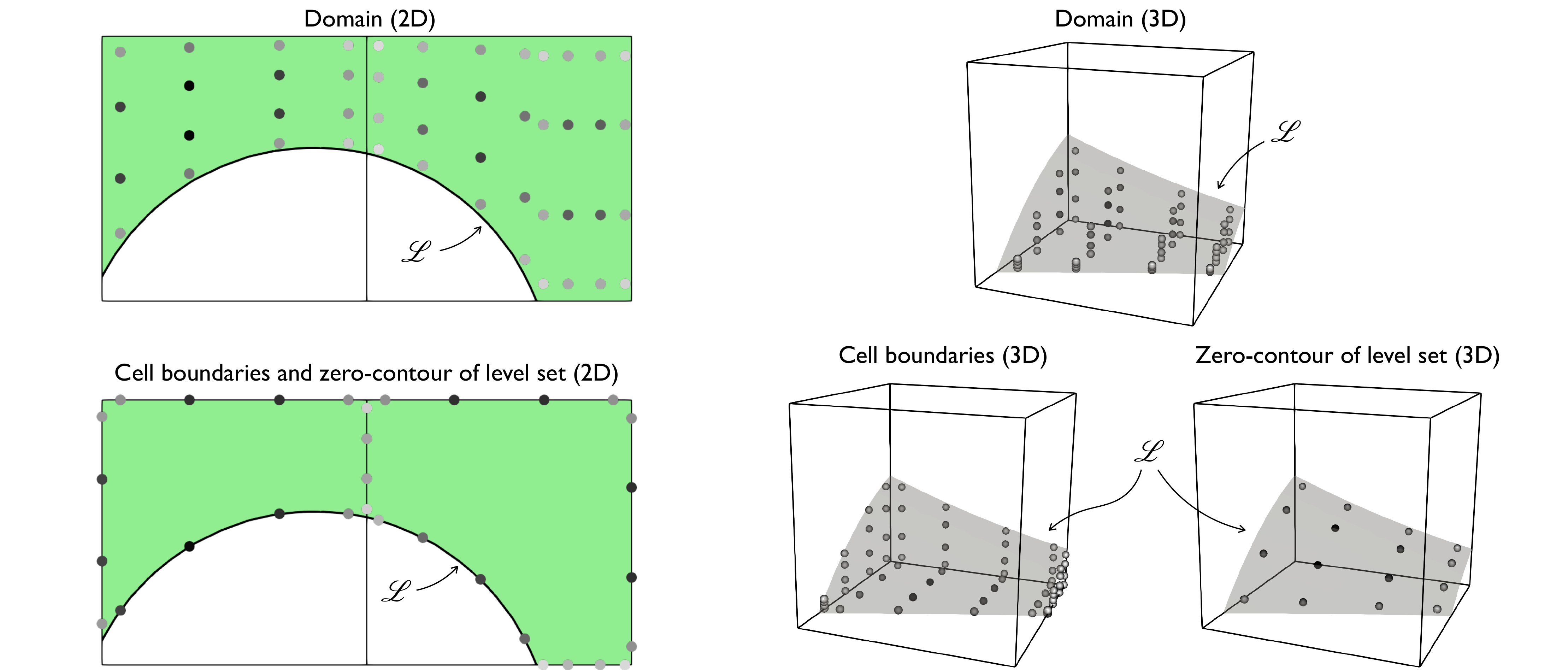}
\caption
{
    Examples of the location of the quadrature points obtained using the algorithm of Ref.\cite{saye2015high};
    in the images, the quadrature points are coloured according to their weight: light grays denote small weights while dark grays denotes large weights.
}
\label{fig:IBVP - QUADRATURE RULES}
\end{figure}

In Eqs.(\ref{eq:governing equations - semidiscrete form - bi-phase solid}) and (\ref{eq:governing equations - mass matrix}), several volumetric and boundary integrals need to be evaluated.
For the entire elements and those cell boundaries that are not cut by the embedded boundaries, see e.g.~Fig.(\ref{fig:IBVP - GRID MESH}), these integrals are evaluated with high-order accuracy using tensor-product Gauss-Legendre quadrature rules.
Meanwhile, to retain the high-order accuracy of the method in proximity of the embedded boundaries, suitable high-order integration schemes should be employed to evaluate the domains and boundary integrals of the cells cut by the zero-contour of the level set function.
Here, we make use of the high-order accurate quadrature algorithms developed in \cite{saye2015high}; an open-source implementation of these algorithms is also available \cite{saye2019algoim}.
A few examples of the kinds of quadrature schemes produced by these algorithms is shown in Fig.(\ref{fig:IBVP - QUADRATURE RULES}); it is worth stressing that the quadrature points are always inside the domain of integration and the quadrature weights are always strictly positive.
The interested reader is referred to Ref.\cite{saye2015high} for a detailed description of the algorithms generating the quadrature rules.

\subsection{$hp$ adaptive mesh refinement}\label{ssec:hpAMR}

\begin{figure}
\centering
\includegraphics[width = \textwidth]{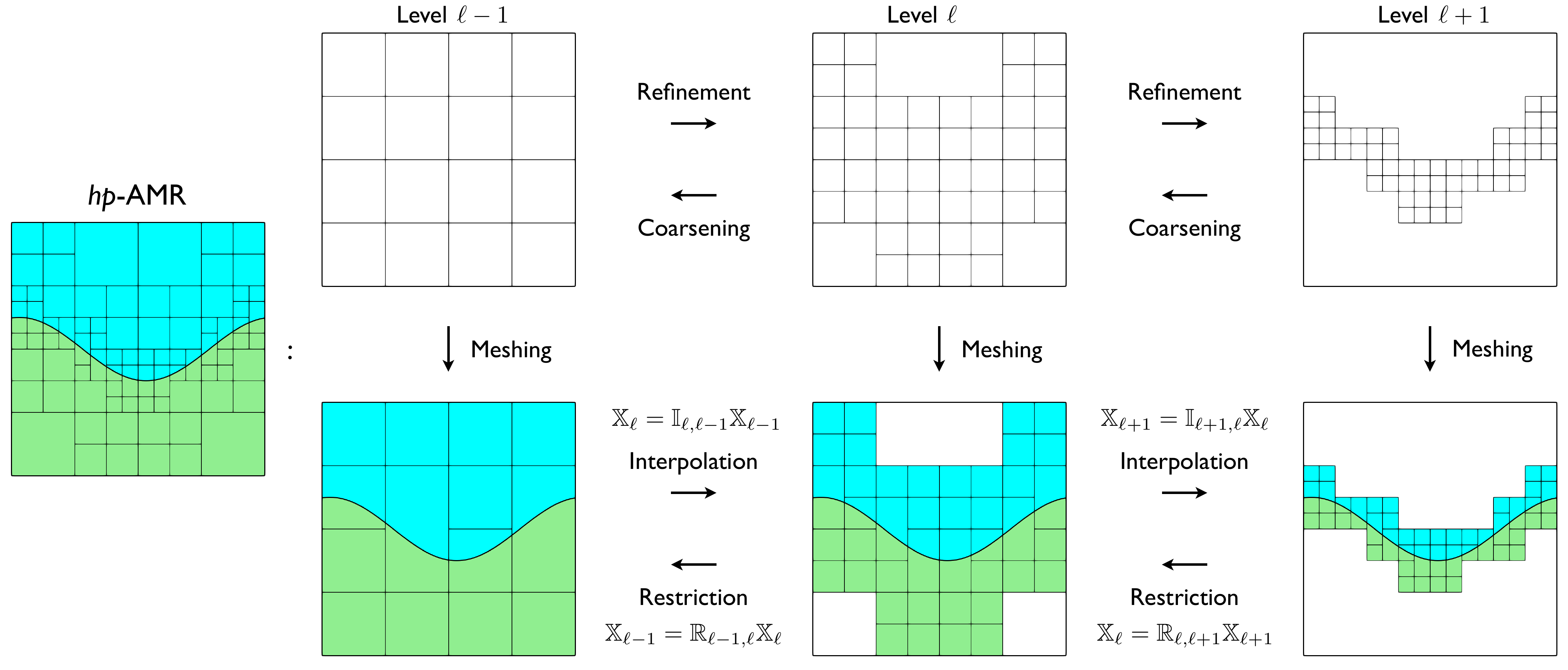}
\caption
{
    Illustration of an $hp$-AMR configuration where the implicitly-defined meshes are generated from a three-level hierarchy of structured grids.
    The figure also sketches the operations performed by the interpolation operator $\mathbb{I}_{\ell+1,\ell}$ and the restriction operator $\mathbb{R}_{\ell,\ell+1}$.
}
\label{fig:IBVP - AMR IR OPERATIONS}
\end{figure}

\begin{figure}
\centering
\includegraphics[width = \textwidth]{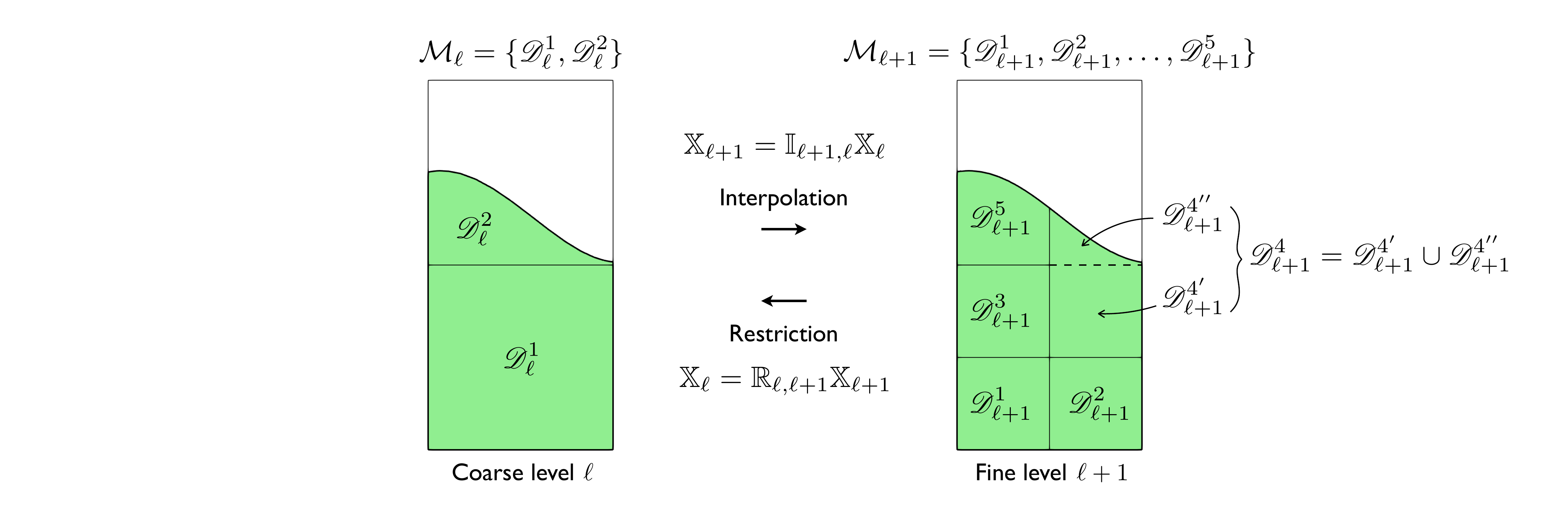}
\caption
{
    Example of an $hp$-AMR configuration where the fine element $\mathscr{D}_{\ell+1}^4$ partially covers the two coarse elements $\mathscr{D}_{\ell}^1$ and $\mathscr{D}_{\ell}^2$.
}
\label{fig:IBVP - AMR DETAIL}
\end{figure}

Owing to the discontinuous nature of DG methods, the present implicit-mesh DG framework can naturally be coupled to an $hp$ adaptive mesh refinement strategy.
To enable AMR capabilities, the framework is integrated into \texttt{AMReX} \cite{zhang2019amrex}, an open-source software library (\url{https://amrex-codes.github.io}) with functionalities for writing massively parallel applications based on adaptive structured Cartesian grids.
In \texttt{AMReX}, the AMR structure is represented as a hierarchy of overlapping levels of refinement ranging from the level identified by $\ell = 0$, which contains the coarsest grid, to the level identified by $\ell = \ell_{\mathrm{finest}}$, which contains the finest grid.
The grid at the level $\ell = 0$ is created at the beginning of the simulations and statically covers the entire domain of analysis; the grids at the levels $\ell > 0$ are created and destroyed dynamically based on user-defined refinement and coarsening criteria.
\texttt{AMReX} has been mainly employed in combination with finite volume schemes, see Ref.\cite{zhang2019amrex} and the references therein, but recently has been extended to include the implicit-mesh DG framework, see Ref.\cite{gulizzi2021coupled}.
Figure (\ref{fig:IBVP - AMR IR OPERATIONS}) illustrates a typical multi-level mesh supported by the present implementation, where the final mesh shown on the leftmost end of the figure is obtained by overlapping three implicitly-defined meshes (the bottom row of images), which in turn are generated from a three-level hierarchy of structured grids (the top row of images) according to the procedure presented in Sec.(\ref{ssec:IM}).
In the figure and in the remainder of the paper, the present AMR scheme is referred to as $hp$-AMR since it allows the use of different polynomial orders at different levels.

It is worth recalling some details about the operations associated with the use of adaptive mesh refinement, even though a thorough discussion regarding adaptive implicitly-defined meshes can be found in Refs.\cite{saye2017implicitI,fortunato2019efficient} as well as Ref.\cite{gulizzi2021coupled}.

The evolution of the $hp$-AMR levels is governed by tagging and un-tagging operations, whereby a two-value variable, or tag, is assigned to each cell at each level.
The tag determines whether the cell should be replaced by a set of finer cells or whether the existing finer resolution is not further required at that location.
Various criteria might be considered to assign and evolve the tags of the cells, with examples ranging from static manual tagging to dynamic solution-dependent tagging.
In this work, we evolve the tag of a generic cell on the basis of the DG solution at that cell.
More specifically, with reference to the cell classification introduced in Sec.(\ref{ssec:IM}), a primary cell is tagged for refinement if the DG solution of the associated implicitly-defined element $\mathscr{D}^e$, which is identified by the coefficients of the element's basis functions contained in $\mathbb{X}^e$, satisfies a condition of the form
\begin{equation}\label{eq:tagging}
f_{\mathrm{tag}}(\mathbb{X}^e) > 0,
\end{equation}
where $f_{\mathrm{tag}}$ is a user-defined function.
Small cells do not have a DG solution directly associated with them but inherit the DG solution of the primary cells with which they are merged; therefore, if a primary cell is tagged for refinement so are all small cells that are merged with it.
Empty cells are never tagged for refinement.
In case of a two-phase solid, oftentimes (typically in proximity of the embedded boundary) an element $\mathscr{D}_\alpha^e$ of the phase $\alpha$ with DG solution $\mathbb{X}_\alpha^e$ and an element $\mathscr{D}_\beta^{e'}$ of the phase $\beta$ with DG solution $\mathbb{X}_\beta^{e'}$ are associated with the same cell $\mathscr{C}$; in such a case, Eq.(\ref{eq:tagging}) will involve both $\mathbb{X}_\alpha^e$ and $\mathbb{X}_\beta^{e'}$.
We note that Eq.(\ref{eq:tagging}) is used not only for refinement but also for coarsening.
In fact, during the course of a numerical simulation, Eq.(\ref{eq:tagging}) is evaluated also for the cells covered by a finer grid and, if the DG solution of an element associated with a primary cell $\mathscr{C}$ ceases to fulfil Eq.(\ref{eq:tagging}), then all the finer cells covering the cell $\mathscr{C}$ and the merged small cells are removed.

When cells are refined or coarsened, the DG solution of the associated implicitly-defined elements must be suitably interpolated at the finer level or restricted at the coarser level, respectively.
These operations are linear, local to the elements and can be implemented using matrix-vector products involving block-sparse operators. 
In particular, consider a vector $\mathbb{X}_\ell$ containing the DG solution of all the implicitly-defined elements at level $\ell$.
Then, as sketched in Fig.(\ref{fig:IBVP - AMR IR OPERATIONS}), the interpolation operator $\mathbb{I}_{\ell+1,\ell}$ transfers the DG solution $\mathbb{X}_\ell$ of the coarse level $\ell$ to the DG solution $\mathbb{X}_{\ell+1}$ of the fine level $\ell+1$.
Formally, the interpolation operation can be written as
\begin{equation}\label{eq:interpolation}
\mathbb{X}_{\ell+1} = \mathbb{I}_{\ell+1,\ell}\mathbb{X}_{\ell},
\end{equation}
where $\mathbb{I}_{\ell+1,\ell}$ is a block-structured matrix computed via the Galerkin projection, so that, given a fine element $\mathscr{D}_{\ell+1}^{e'}$ of level $\ell+1$ and a coarse element $\mathscr{D}_{\ell}^{e}$ of level $\ell$, the block $\mathbb{I}_{\ell+1,\ell}^{e',e}$ is given by
\begin{equation}\label{eq:interpolation - Galerkin projection}
\mathbb{I}_{\ell+1,\ell}^{e',e} \equiv (\mathbb{M}_{\ell+1}^{e'})^{-1}\int_{\mathscr{D}_{\ell+1}^{e'}\cap\mathscr{D}_{\ell}^{e}}\mathbb{B}_{\ell+1}^{e'\Tr}\mathbb{B}_{\ell}^{e}~\dd{V}.
\end{equation}
To illustrate how $\mathbb{I}_{\ell+1,\ell}$ acts on $\mathbb{X}_{\ell}$, consider Fig.(\ref{fig:IBVP - AMR DETAIL}), which shows a two-element coarse mesh $\mathcal{M}_{\ell} \equiv \{\mathscr{D}_{\ell}^1,\mathscr{D}_{\ell}^2\}$ generated from a $1\times 2$ coarse grid and a five-element fine mesh $\mathcal{M}_{{\ell+1}} \equiv \{\mathscr{D}_{\ell+1}^1,\mathscr{D}_{\ell+1}^2,\dots,\mathscr{D}_{\ell+1}^5\}$ generated from a $2\times 4$ grid that is obtained by refining the coarse grid with a refinement ratio of 2.
Using Eq.(\ref{eq:interpolation - Galerkin projection}) for the AMR configuration of Fig.(\ref{fig:IBVP - AMR DETAIL}), Eq.(\ref{eq:interpolation}) becomes
\begin{equation}\label{eq:interpolation - example}
\left(\begin{array}{c}
\mathbb{X}_{\ell+1}^{1}\\
\mathbb{X}_{\ell+1}^{2}\\
\mathbb{X}_{\ell+1}^{3}\\
\mathbb{X}_{\ell+1}^{4}\\
\mathbb{X}_{\ell+1}^{5}
\end{array}\right)
=
\left[\begin{array}{cc}
\mathbb{I}_{\ell+1,\ell}^{1,1}&\bm{0}\\
\mathbb{I}_{\ell+1,\ell}^{2,1}&\bm{0}\\
\mathbb{I}_{\ell+1,\ell}^{3,1}&\bm{0}\\
\mathbb{I}_{\ell+1,\ell}^{4,1}&\mathbb{I}_{\ell+1,\ell}^{4,2}\\
\bm{0}&\mathbb{I}_{\ell+1,\ell}^{5,2}
\end{array}\right]
\left(\begin{array}{c}
\mathbb{X}_{\ell}^{1}\\
\mathbb{X}_{\ell}^{2}
\end{array}\right).
\end{equation}
It is interesting to note that Eq.(\ref{eq:interpolation - example}) reflects the configurations of the meshes of Fig.(\ref{fig:IBVP - AMR DETAIL}), including the case of the fine element $\mathscr{D}_{\ell+1}^4$ that partially covers the two distinct coarse elements $\mathscr{D}_{\ell}^1$ and $\mathscr{D}_{\ell}^2$.

The counterpart of the interpolation operation is the restriction operation, whereby the DG solution $\mathbb{X}_{\ell+1}$ of the fine level $\ell+1$ is transferred to the DG solution $\mathbb{X}_\ell$ of the coarse level $\ell$.
Similar to Eq.(\ref{eq:interpolation}), this operation can formally be written as
\begin{equation}\label{eq:restriction}
\mathbb{X}_{\ell} = \mathbb{R}_{\ell,\ell+1}\mathbb{X}_{\ell+1},
\end{equation}
where, using the Galerkin projection, the restriction operator $\mathbb{R}_{\ell,\ell+1}$ is related to the interpolation operator such that
\begin{equation}\label{eq:restriction - from interpolation}
\mathbb{R}_{\ell,\ell+1} = \mathbb{M}_{\ell}^{-1}\mathbb{I}_{\ell+1,\ell}^{\Tr}\mathbb{M}_{\ell+1},
\end{equation}
where $\mathbb{M}_{\ell}$ is the block-diagonal mass matrix of the implicitly-defined elements of the level $\ell$.

It is worth noting that the interpolation and restriction operators defined via Eqs.(\ref{eq:interpolation - Galerkin projection}) and (\ref{eq:restriction}), respectively, are valid regardless of the choice of basis functions and, therefore, naturally enable the use of different polynomial orders at different AMR levels.
From an implementation viewpoint, the present DG scheme requires the evaluation of the interpolation operators and the Cholesky decomposition of the mass matrices, while the restriction operator can be applied on-the-fly using Eq.(\ref{eq:restriction - from interpolation}).
Moreover, all standard (hyper)rectangular elements (which represent the majority of the mesh elements) share the same mass matrix, which can precomputed and stored at the beginning of the simulations; the same applies to the interpolation operator between two standard elements of two different AMR levels.
Conversely, the mass matrices and the interpolation operators of the cut elements are in general unique and are computed via Eq.(\ref{eq:governing equations - mass matrix}) and Eq.(\ref{eq:interpolation - Galerkin projection}), respectively, using high-order quadrature rules.

\subsection{Time-stepping}
The last aspect of the numerical framework regards the time-stepping, i.e.~the update in time of the coefficients of the spatial basis functions.
Whether a single-level or an $hp$-AMR scheme is considered, the time-evolution equation for the coefficients $\mathbb{X}_\alpha^e$ of a generic $e$-th element $\mathscr{D}_\alpha^e$ belonging to the phase $\alpha$ can be written as
\begin{equation}\label{eq:ode}
\mathbb{M}_\alpha^e\dot{\mathbb{X}}{}_\alpha^e = \mathbb{A}_\alpha^e(t,\mathbb{X}),
\end{equation}
where $\mathbb{A}_\alpha^e(t,\mathbb{X})$ stems from the evaluation of the right-hand side of Eq.(\ref{eq:governing equations - semidiscrete form - bi-phase solid}) and $\mathbb{X}$ formally contains the coefficients of all the mesh elements;
note however that only the DG solution from the neighboring elements of $\mathscr{D}_\alpha^e$ is required to compute $\mathbb{A}_\alpha^e(t,\mathbb{X})$.
Integration in time of Eq.(\ref{eq:ode}) is performed via an explicit high-order Runge-Kutta algorithm \cite{cockburn1998runge} matching the order of the highest spatial discretization among the mesh levels.
As explicit time-integration schemes are conditionally stable, at a generic level $\ell$ with mesh size $h_\ell$ and using a DG scheme with polynomial degree $p_\ell$, the maximum time step $\tau_\ell$ is subject to the following CFL condition
\begin{equation}\label{eq:CFL}
\frac{\tau_\ell}{h_\ell} < \frac{C_\ell \overline{f}_\ell}{c(1+2p_\ell)},
\end{equation}
where $\overline{f}_\ell$ is the volume fraction threshold triggering the cell-merging at the level $\ell$ and $C_\ell$ is a constant smaller than $1$ that does not depend on $h_\ell$ or $p_\ell$; in all simulations presented in Sec.(\ref{sec:RESULTS}), $\overline{f}_0 = \overline{f}_1 = 0.3$ and ${C}_0 = {C}_1 = 0.833$.
Moreover, in Eq.(\ref{eq:CFL}), $c \equiv c_\alpha$ for single-phase solids or $c \equiv \max\{c_\alpha,c_\beta\}$ for bi-phase solids, where $c_\alpha$ and $c_\beta$ are the maximal speeds of the elastic waves in the phases $\alpha$ and $\beta$, respectively.
For isotropic solids, $c_\alpha$ coincides with the speed of the P-waves in the phase $\alpha$.
For general anisotropic solids, the wave speed depends on the direction of propagation and, therefore, the maximum wave speed is evaluated as \cite{carcione2007wave}
\begin{equation}\label{eq:maximum wave speed}
c_\alpha^2 \equiv \max_{\bm{n}} \{\mathrm{eig}(\bm{\sGamma}_\alpha(\bm{n}))\},
\end{equation}
where $\bm{\sGamma}_\alpha(\bm{n}) \equiv \rho_\alpha^{-1}\bm{I}_{n}^\Tr\bm{c}_\alpha\bm{I}_{n}$, $\bm{I}_n \equiv n_i\bm{I}_i$ and $\mathrm{eig}(\bullet)$ returns the eigenvalues of $\bullet$.
Finally, the time step $\tau$ of the Runge-Kutta algorithm is $\tau \equiv \min_{\ell}\tau_\ell$.
\section{Results}\label{sec:RESULTS}
In this section, the capabilities of the presented implicit-mesh DG framework are assessed for two- and three-dimensional test cases involving wave propagation in single- and bi-phase elastic solids.

The numerical simulations use implicitly-defined meshes generated either from uniform grids with mesh size $h$ or from a two-level $hp$-AMR, where level 0 and level 1 have mesh size $h_0$ and $h_1$, respectively.
For the simulations using the $hp$-AMR, the coarse level $0$ is generated at the beginning of the simulation and is kept fixed, while the fine level $1$ is dynamically updated during the time evolution by refining the coarse cells with a refinement ratio $r$ such that $h_1 = h_0/r$ and the number of fine cells replacing a coarse cell is $r^d$.
Tensor-product Legendre polynomials of degree $p$ are employed to define the space $\mathcal{P}_{hp}^e$, and thus the space $\mathcal{V}_{hp}$ introduced in Eq.(\ref{eq:space}); the corresponding DG scheme is denoted by DG$_p$.
We recall that the present $hp$-AMR strategy allows the use of different polynomial orders for different AMR levels.

\subsection{Convergence analysis}
Reported here are the results of several convergence tests on single- and bi-phase solids with isotropic, orthotropic and anisotropic constitutive behavior, in two and three-dimensions.
For both the single-phase solid simulations and the bi-phase solid simulations, we assume that the phases have density $\rho_\alpha = \rho_\beta = \rho$ and stiffness $\bm{c}_\alpha = \bm{c}_\beta = \bm{c}$.
In this section, all quantities are assumed non-dimensional.
The considered isotropic solid has density $\rho = 1$, Young's modulus $Y = 1$ and Poisson's ratio $\nu = 0.3$.
The considered orthotropic solid is a FCC Copper crystal \cite{benedetti2013three} with density $\rho = 8.92$ and non-zero elastic constants $c_{11} = 168$, $c_{12} = 121$, $c_{44} = 75$.
In 2D, the constitutive behavior of the considered anisotropic solid represents the in-plane behavior of a multilayered composite material \cite{gulizzi2019novel} with density $\rho = 1.6$ and stiffness matrix
\begin{equation}
\bm{c} = \left[\begin{array}{ccc}
0.5637 & 0.2963 & 0.3158 \\
       & 0.5637 & 0.3158 \\
\mathrm{Sym.}&  & 0.3111
\end{array}\right].
\end{equation}
In 3D, the considered anisotropic solid is an Olivine crystal \cite{browaeys2004decomposition,de2007arbitrary} with density $\rho = 1.0$ and whose orthorhombic axes are tilted and aligned with the directions $[1,1,1]$, $[-1,1,0]$ and $[-1,-1,2]$ such that the stiffness matrix in the global reference system is
\begin{equation}
\bm{c} = \left[\begin{array}{cccccc}
185.8 &  67.3 &  76.2 & 2.69 &   17.6 & -5.44 \\
      & 170.0 & 62.67 & 4.62 &  -6.60 & -6.53 \\
      &       & 219.8 & 3.08 &  30.48 & -2.72 \\
      &       &       & 59.0 & -1.905 &  2.83 \\
      &\mathrm{Sym.}& &      &   79.2 &  2.12 \\
      &       &       &      &        &  57.0
\end{array}\right].
\end{equation}

We start by constructing an exact solution of Eq.(\ref{eq:governing equations - vector form}) with zero source term.
Consider a plane-wave vector field of the form
\begin{equation}\label{eq:plane wave}
\bm{U}(t,\bm{x}) = \wtilde{\bm{U}}\sin(\omega t-\kappa_i x_i),
\end{equation}
where $\kappa_i$ is the $i$-th component of the wave vector $\bm{\kappa} = (\kappa_1,\dots,\kappa_d)$, $\omega$ is the angular frequency and $\wtilde{\bm{U}}$ is a constant vector.
For a given choice of $\bm{\kappa}$, by plugging Eq.(\ref{eq:plane wave}) into Eq.(\ref{eq:governing equations - vector form}) with $\bm{S}_\alpha = \bm{0}$, $\omega$ and $\wtilde{\bm{U}}$ are obtained as the eigenvalue and the eigenvector, respectively, of the eigenvalue problem
\begin{equation}\label{eq:eigenvalue problem}
\bm{A}_\kappa \wtilde{\bm{U}} = \omega \wtilde{\bm{U}},
\end{equation}
where $\bm{A}_\kappa$ is $N_U\times N_U$ matrix given by
\begin{equation}
\bm{A}_\kappa \equiv
\left[\begin{array}{cc}
\bm{0}&\bm{I}_\kappa^\Tr\bm{c}\\
\rho^{-1}\bm{I}_\kappa&\bm{0}
\end{array}\right]
\quad\mathrm{and}\quad
\bm{I}_\kappa \equiv \kappa_i\bm{I}_i.
\end{equation}
There are in general $N_U$ couples $\{\omega_v,\wtilde{\bm{U}}_v\}$, $v = 1,\dots,N_U$, that are solution of Eq.(\ref{eq:eigenvalue problem}); see Ref.\cite{zhan2018exact} for more details about the eigenvalue properties of the matrix $\bm{A}_\kappa$.
Then, an exact solution Eq.(\ref{eq:governing equations - vector form}) can be written as the following linear superposition of plane waves
\begin{equation}\label{eq:exact solution}
\bm{U}^{\mathrm{exact}}(t,\bm{x}) = \sum_{v = 1}^{N_U}\wtilde{\bm{U}}_v\sin(\omega_v t-\kappa_i x_i),
\end{equation}
where the vectors $\wtilde{\bm{U}}_v$, $v = 1,\dots,N_U$, are normalized to have unit amplitude.
In all simulations, the exact solution $\bm{U}^\mathrm{exact}$ given in Eq.(\ref{eq:exact solution}) is employed to set the initial conditions as $\bm{U}_{\alpha 0}(\bm{x}) = \bm{U}_{\beta 0}(\bm{x}) = \bm{U}^\mathrm{exact}(0,\bm{x})$, while the maximum eigenvalue $\omega_{\mathrm{max}} \equiv \max_{v}\omega_v$ determines the final time $T$ of evolution as $T = 2\pi/\omega_{\mathrm{max}}$.

As the last ingredient of this convergence analysis, we introduce two error measures between the solution $\bm{U}$ computed via the present DG scheme and the exact solution $\bm{U}^\mathrm{exact}$ given in Eq.(\ref{eq:exact solution}).
The error measures are
\begin{equation}\label{eq:error measures}
e_{L_\infty}(\bm{U},\bm{U}^\mathrm{exact}) \equiv \frac{||\bm{U}-\bm{U}^\mathrm{exact}||_{L_\infty}}{||\bm{U}^\mathrm{exact}||_{L_\infty}}
\quad\mathrm{and}\quad
e_{L_2}(\bm{U},\bm{U}^\mathrm{exact}) \equiv \left[\frac{E(\bm{U}-\bm{U}^\mathrm{exact})}{E(\bm{U}^\mathrm{exact})}\right]^{1/2},
\end{equation}
where the $L_\infty$ norm is evaluated by computing the maximum error at the quadrature points among all the components of $\bm{U}$ and the $e_{L_2}$ error is evaluated by introducing an energy norm.
In particular, let $E_\alpha(\bm{U})$ denote the energy associated with the solution $\bm{U}$ for the phase $\alpha$ given by
\begin{equation}\label{eq:energy}
E_\alpha(\bm{U}) = \frac{1}{2}\int_{\mathscr{D}_\alpha}\left(\rho_\alpha^{-1}\bm{m}_\alpha^\Tr\bm{m}_\alpha+\bm{\gamma}_{\alpha}^\Tr\bm{c}_{\alpha}\bm{\gamma}_{\alpha}\right)\dd{V};
\end{equation}
then, in Eq.(\ref{eq:error measures}), $E(\bm{U}) = E_\alpha(\bm{U})$ and $E(\bm{U}) = E_\alpha(\bm{U})+E_\beta(\bm{U})$ for single-phase and bi-phase solids, respectively.

\begin{figure}
\centering
\includegraphics[width = \textwidth]{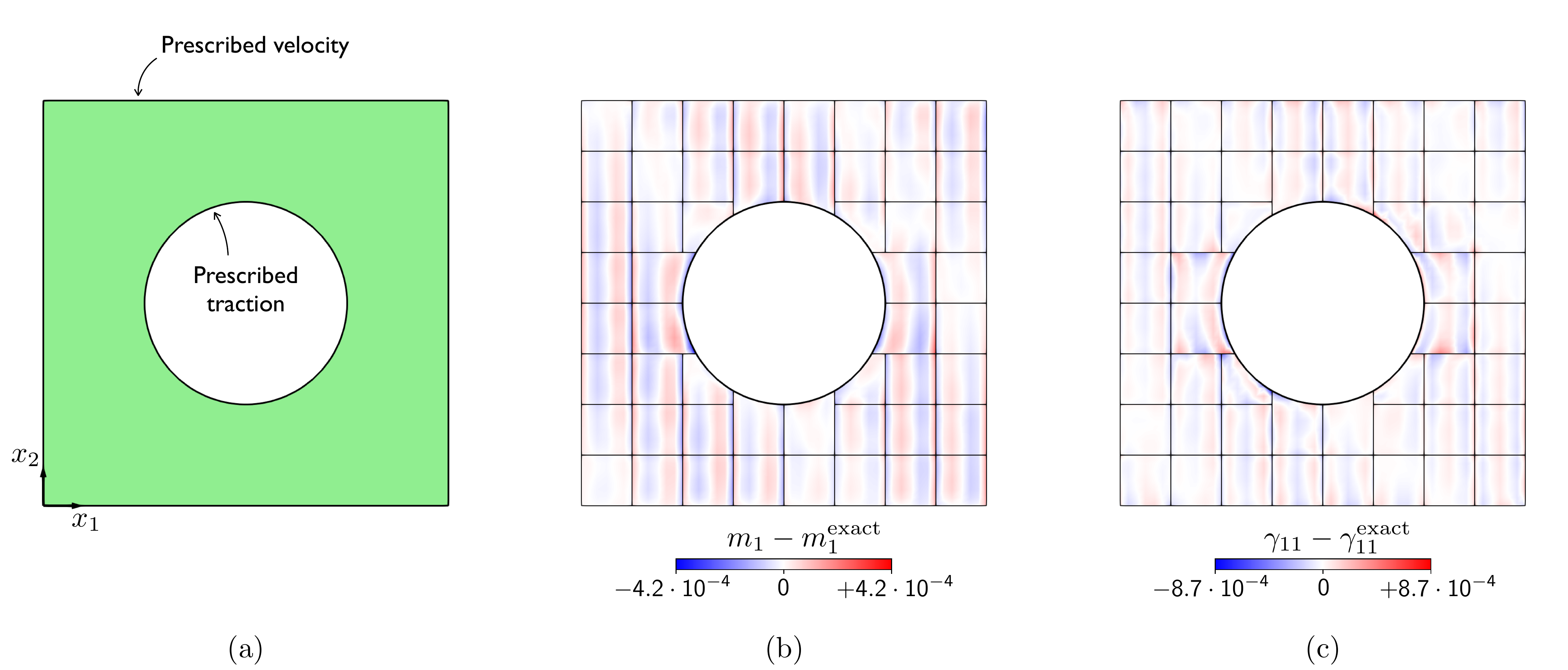}
\caption
{
    (a) Geometry and boundary conditions for the single-phase solid employed in the 2D $hp$-convergence analysis.
    Error in (b) the momentum component $m_1$ and (c) the strain component $\gamma_{11}$ obtained with an implicitly-defined mesh generated from an $8^2$ background grid and a DG$_3$ scheme for the anisotropic material response.
}
\label{fig:RESULTS - CONVERGENCE SINGLE PHASE DOMAIN - GEOMETRY ERROR 2D}
\end{figure}

\begin{figure}
\centering
\begin{subfigure}{0.32\textwidth}
\includegraphics[width=\textwidth]{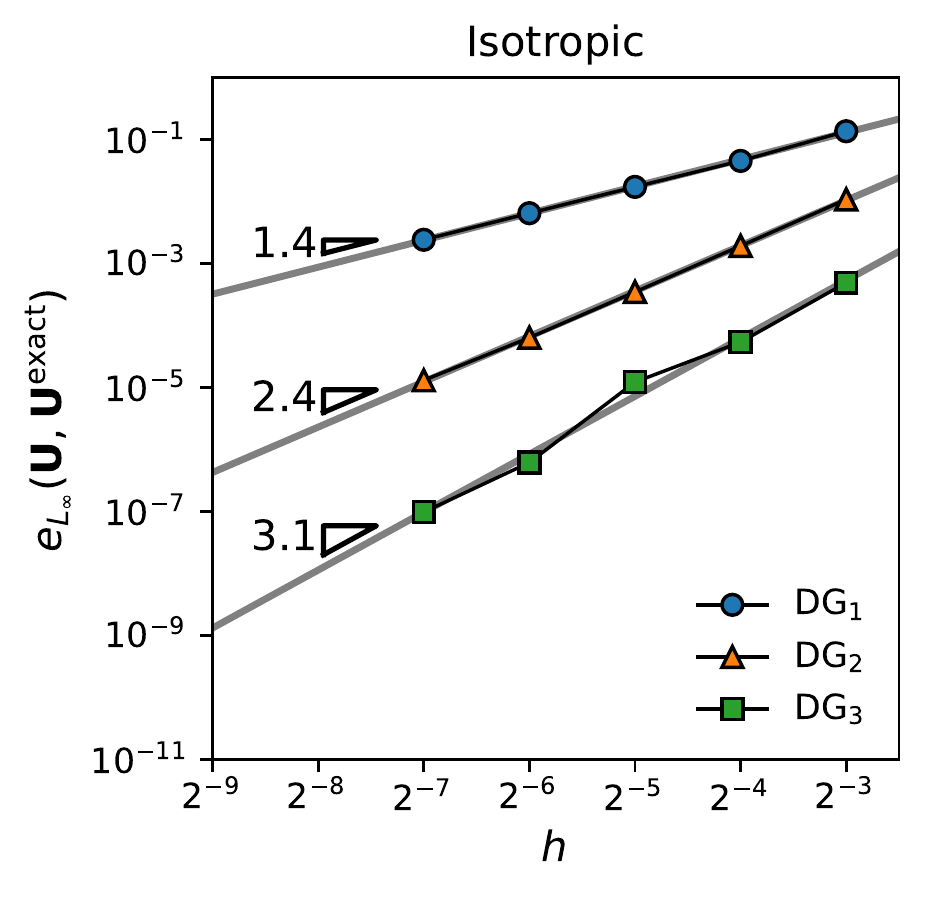}
\end{subfigure}
\
\begin{subfigure}{0.32\textwidth}
\includegraphics[width=\textwidth]{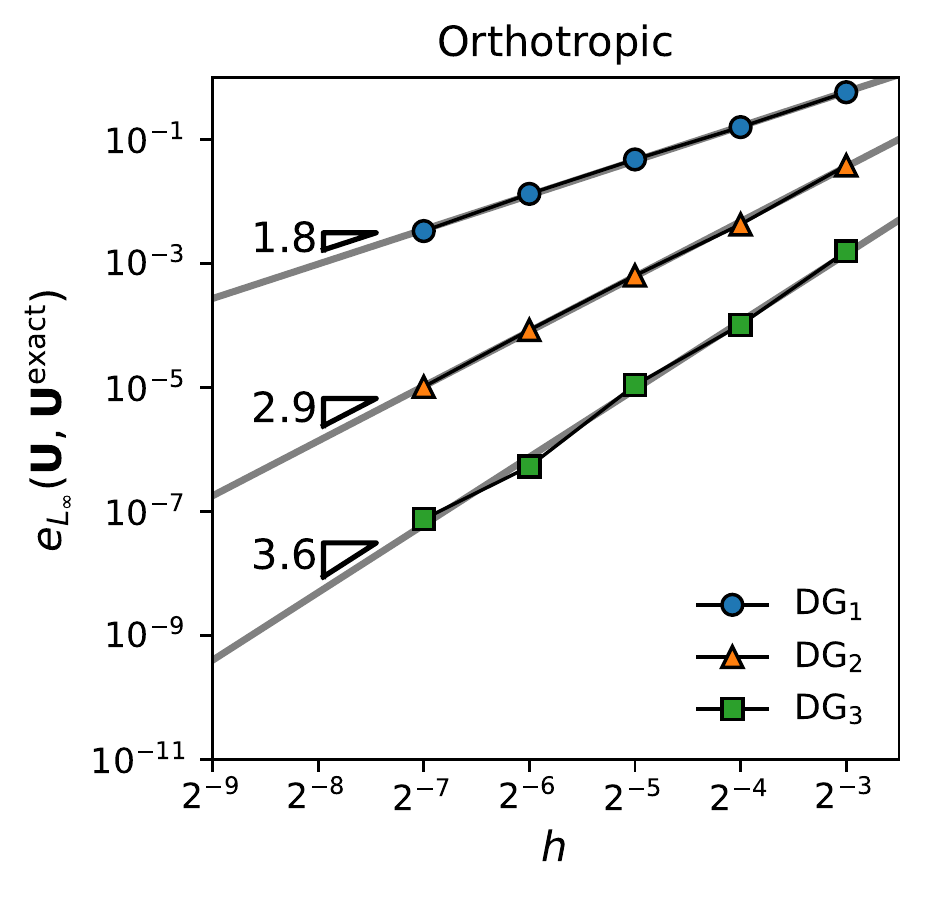}
\end{subfigure}
\
\begin{subfigure}{0.32\textwidth}
\includegraphics[width=\textwidth]{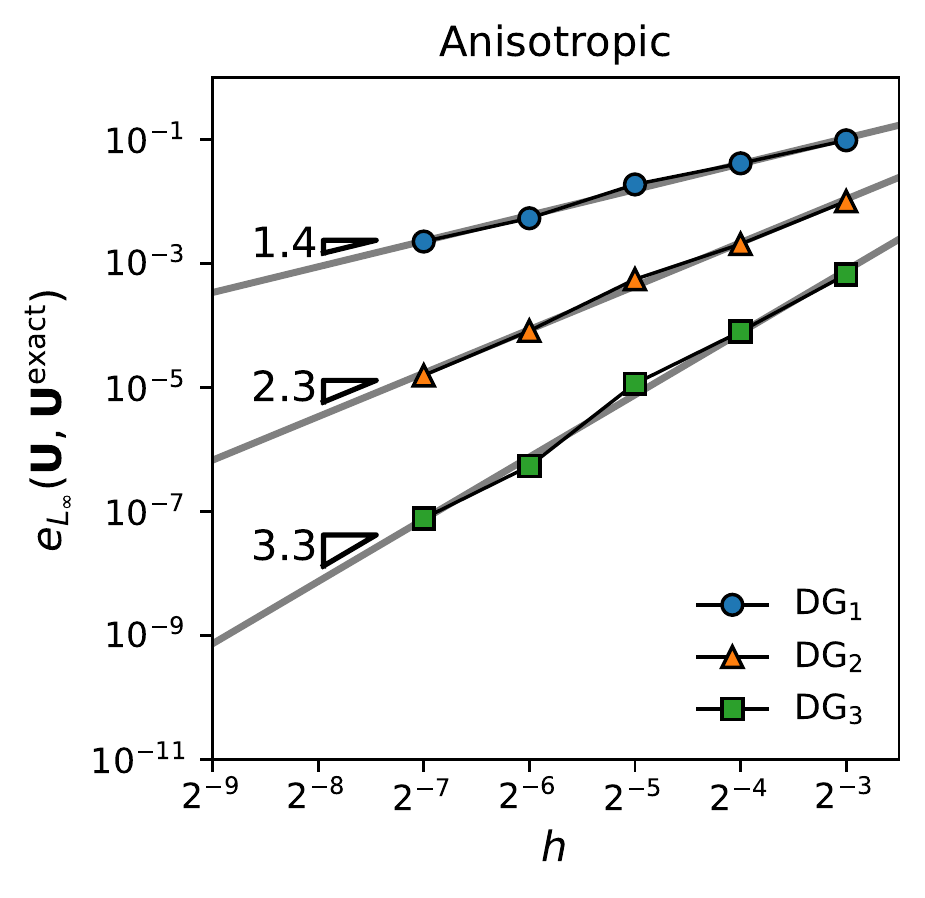}
\end{subfigure}
\\
\begin{subfigure}{0.32\textwidth}
\includegraphics[width=\textwidth]{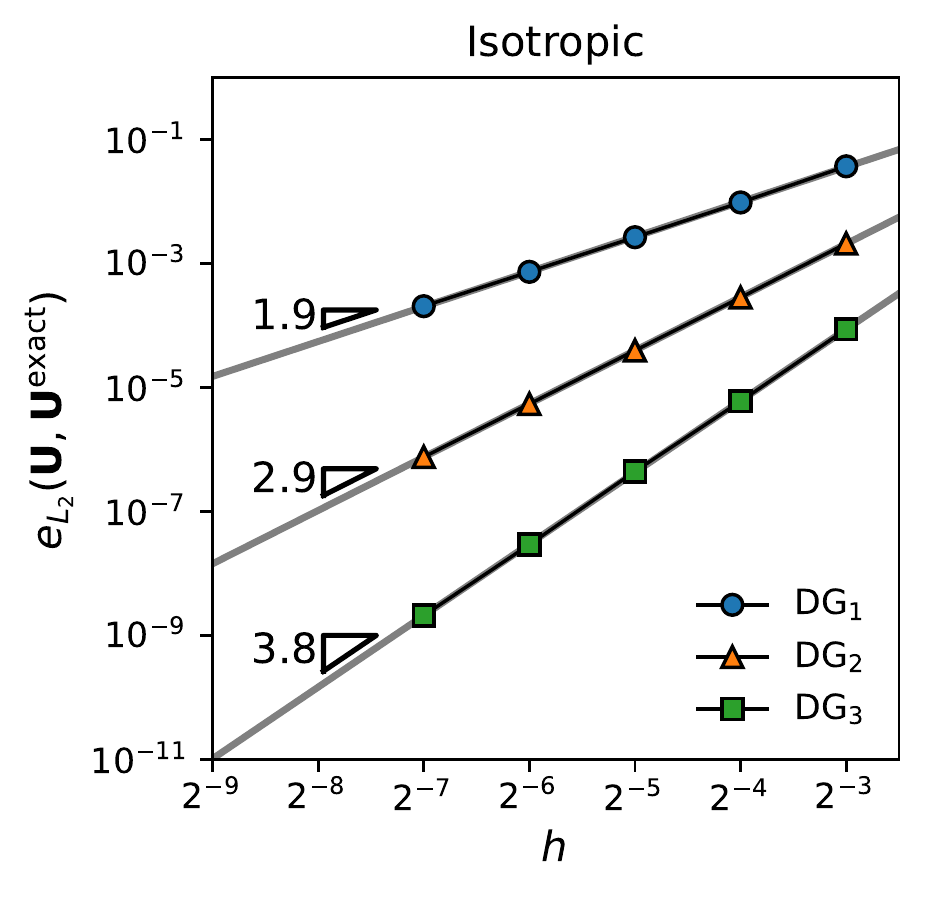}
\end{subfigure}
\
\begin{subfigure}{0.32\textwidth}
\includegraphics[width=\textwidth]{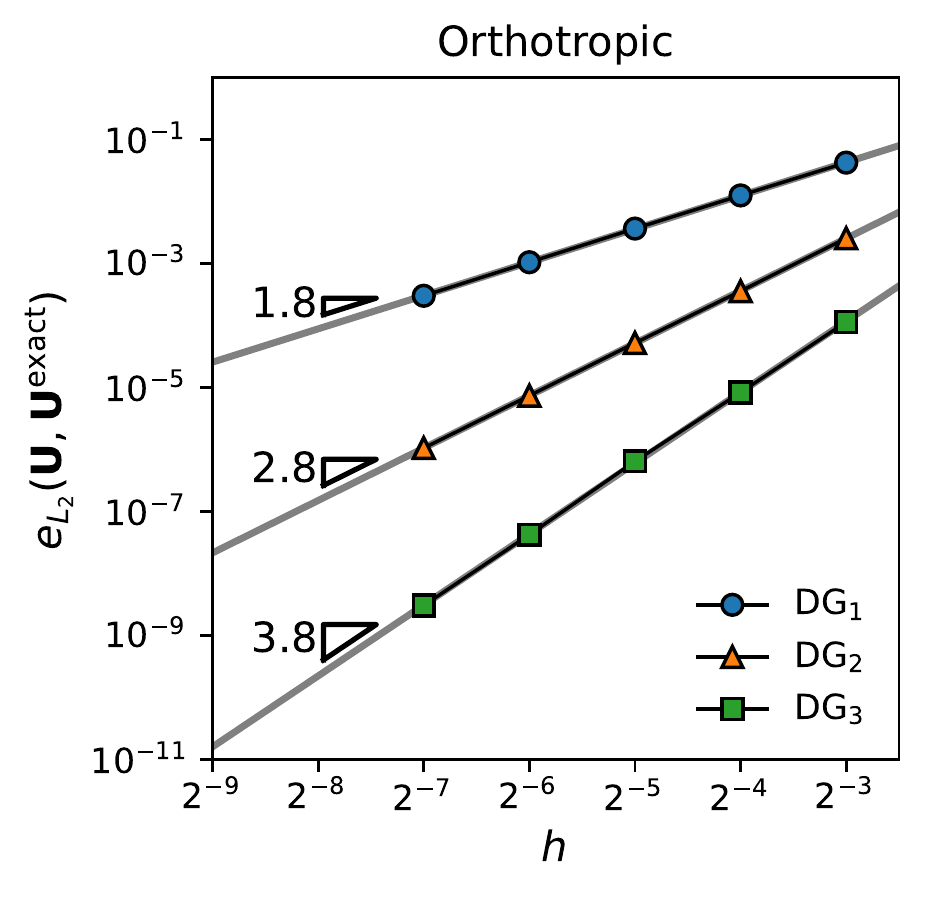}
\end{subfigure}
\
\begin{subfigure}{0.32\textwidth}
\includegraphics[width=\textwidth]{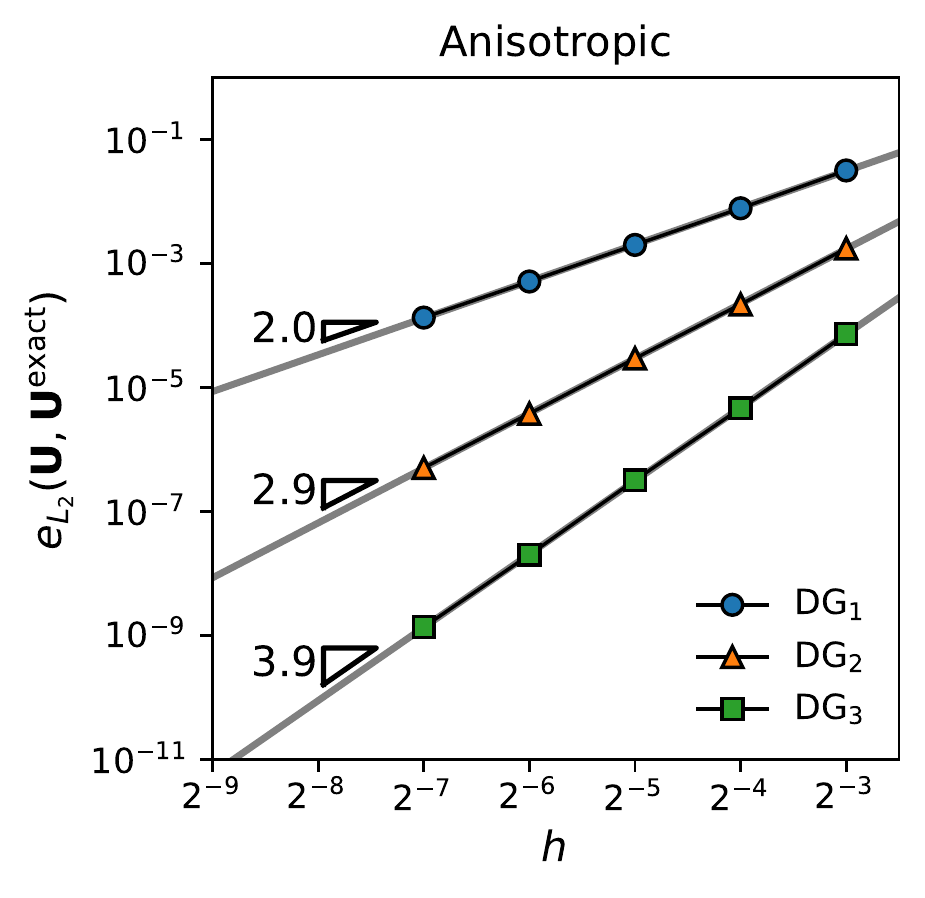}
\end{subfigure}
\caption
{
    (top row) $e_{L_\infty}$ error and (bottom row) $e_{L_2}$ error for the 2D single-phase solid of Fig.(\ref{fig:RESULTS - CONVERGENCE SINGLE PHASE DOMAIN - GEOMETRY ERROR 2D}) for different constitutive behaviors.
}
\label{fig:RESULTS - CONVERGENCE SINGLE PHASE DOMAIN - HP STUDY 2D}
\end{figure}

Figure (\ref{fig:RESULTS - CONVERGENCE SINGLE PHASE DOMAIN - GEOMETRY ERROR 2D}a) shows the geometry and the boundary conditions for the 2D single-phase solid case.
The geometry consists of a square with a circular cavity and is defined in the background unit square $[0,1]^2$ by the level set function
\begin{equation}
\varphi(\bm{x}) = R^2-(x_1-o_1)^2-(x_2-o_2)^2,
\end{equation}
where $R = 0.25$ and $o_1 = o_2 = 0.5$.
A velocity field $\overline{\bm{v}}$ is prescribed on the outer boundary $\mathscr{B}_\alpha$ of the background square, whereas a traction field $\overline{\bm{t}}$ is prescribed on the zero contour $\mathscr{L}_\alpha$ of the level set function.
The exact solution $\bm{U}^\mathrm{exact}$ is specified by $\kappa_1 = 2\pi\cos(\pi/6)$ and $\kappa_2 = 2\pi\sin(\pi/6)$ and is employed to evaluate the fields $\overline{\bm{v}}$ and $\overline{\bm{t}}$ at any $(t,\bm{x})$ on the geometry's boundaries.
Figures (\ref{fig:RESULTS - CONVERGENCE SINGLE PHASE DOMAIN - GEOMETRY ERROR 2D}b) and (\ref{fig:RESULTS - CONVERGENCE SINGLE PHASE DOMAIN - GEOMETRY ERROR 2D}c) show the error in the momentum component $m_1$ and the strain component $\gamma_{11}$, respectively, when an $8^2$ background grid and a DG$_3$ scheme are employed for the anisotropic material response case.
The figures also display the implicitly-defined mesh generated from the background grid.
$hp$-convergence plots of the two error measures given in Eq.(\ref{eq:error measures}) for the 2D single-phase solid with isotropic, orthotropic and anisotropic behavior are shown in Fig.(\ref{fig:RESULTS - CONVERGENCE SINGLE PHASE DOMAIN - HP STUDY 2D}).

\begin{figure}
\centering
\includegraphics[width = \textwidth]{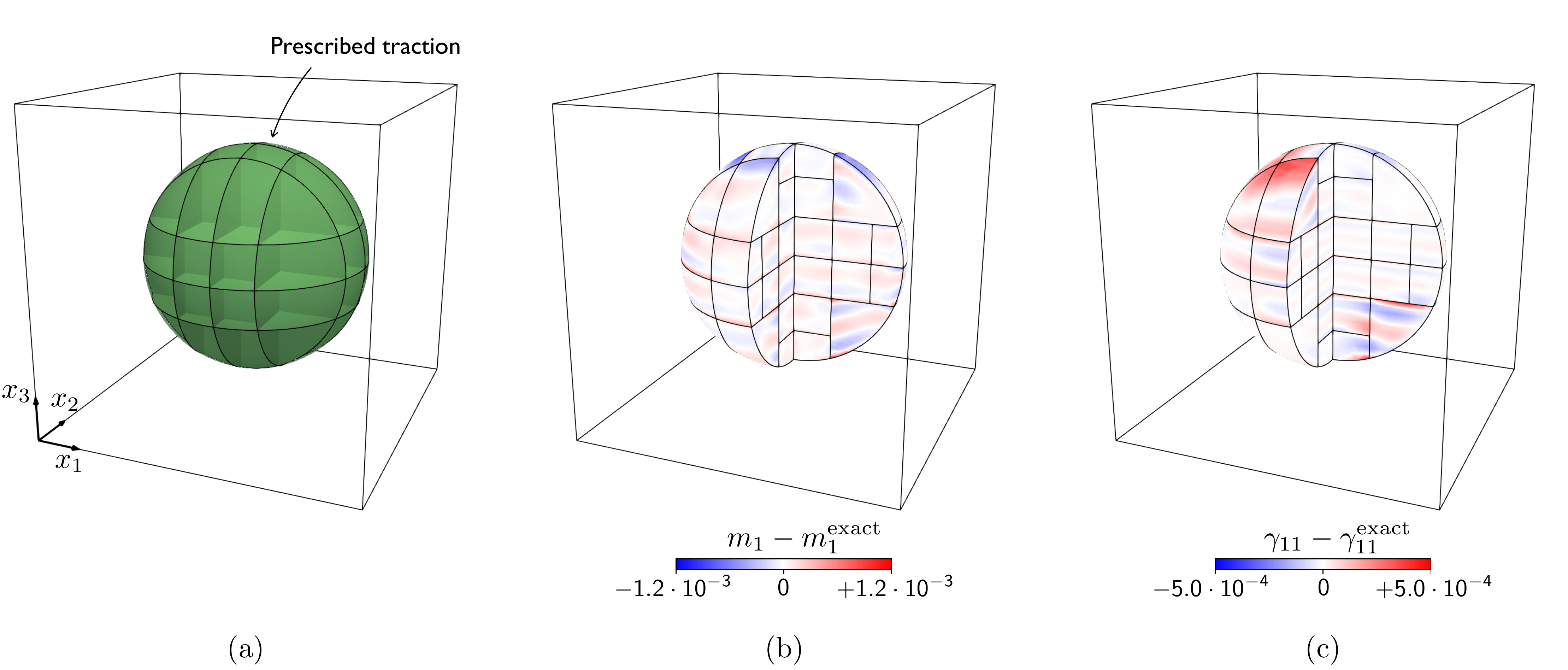}
\caption
{
    (a) Geometry, boundary conditions and implicitly-defined mesh generated from an $8^3$ background grid for the single-phase solid employed in the 3D $hp$-convergence analysis.
    Error in (b) the momentum component $m_1$ and (c) the strain component $\gamma_{11}$ obtained with the implicitly-defined mesh of figure (a) and a DG$_3$ scheme for the anisotropic material response.
}
\label{fig:RESULTS - CONVERGENCE SINGLE PHASE DOMAIN - GEOMETRY ERROR 3D}
\end{figure}

\begin{figure}
\centering
\begin{subfigure}{0.32\textwidth}
\includegraphics[width=\textwidth]{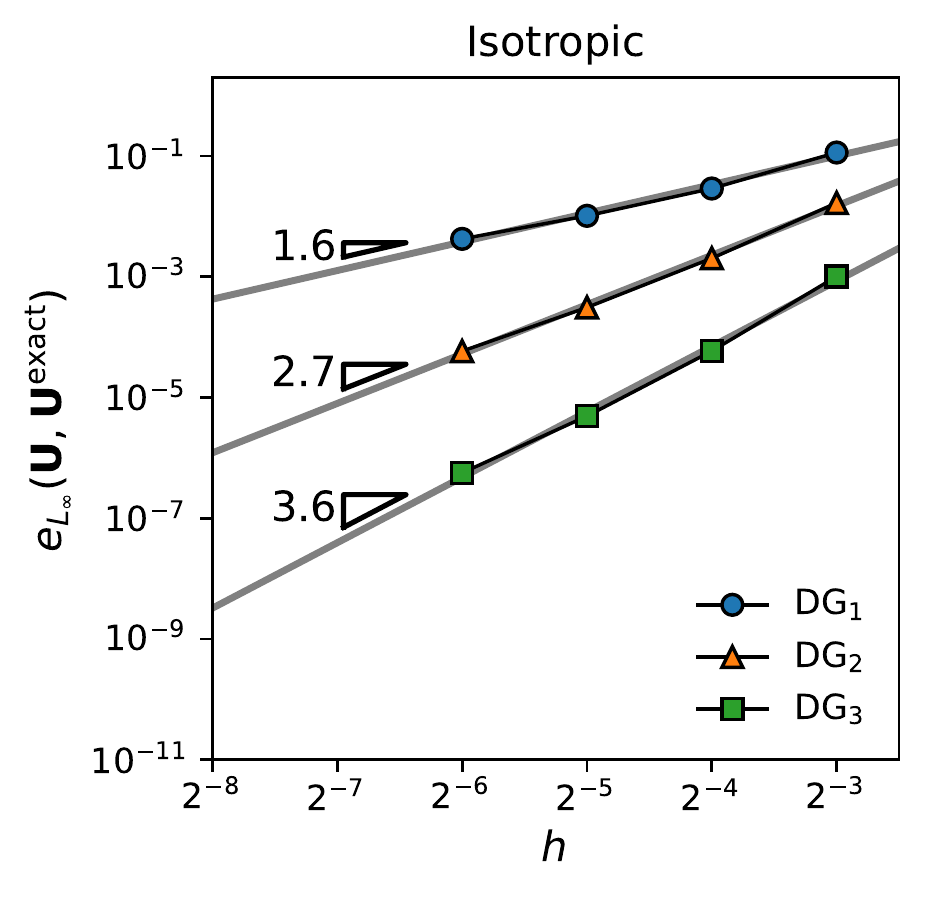}
\end{subfigure}
\
\begin{subfigure}{0.32\textwidth}
\includegraphics[width=\textwidth]{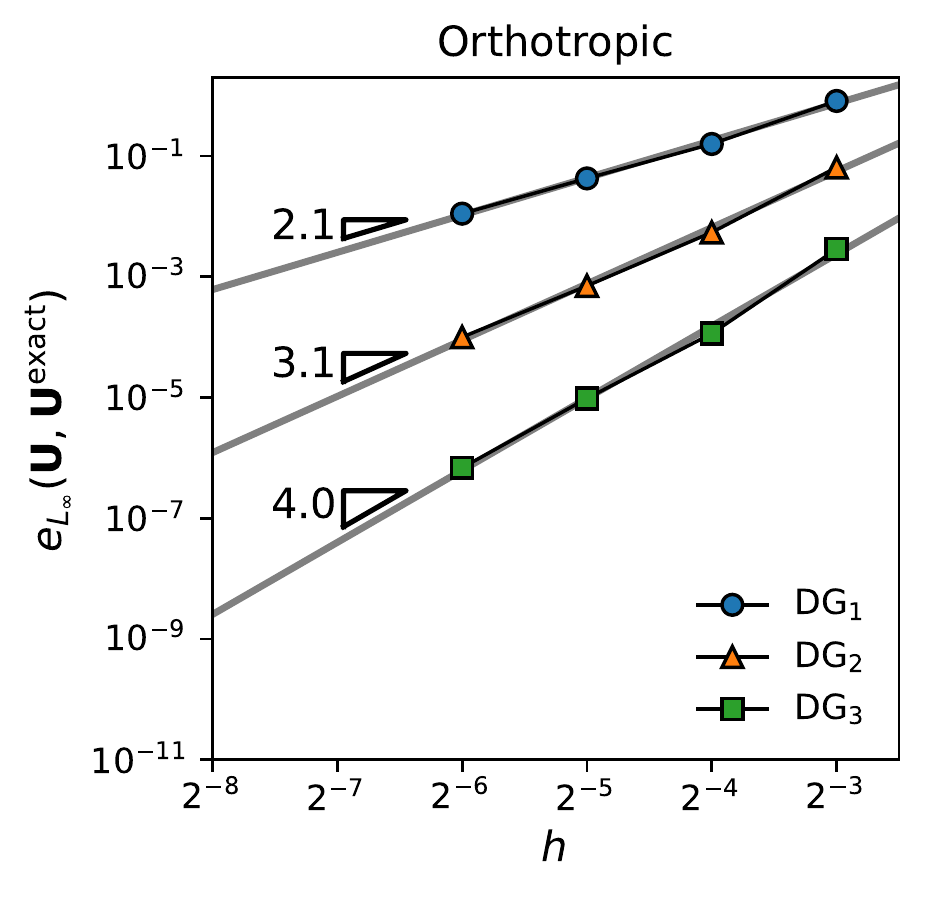}
\end{subfigure}
\
\begin{subfigure}{0.32\textwidth}
\includegraphics[width=\textwidth]{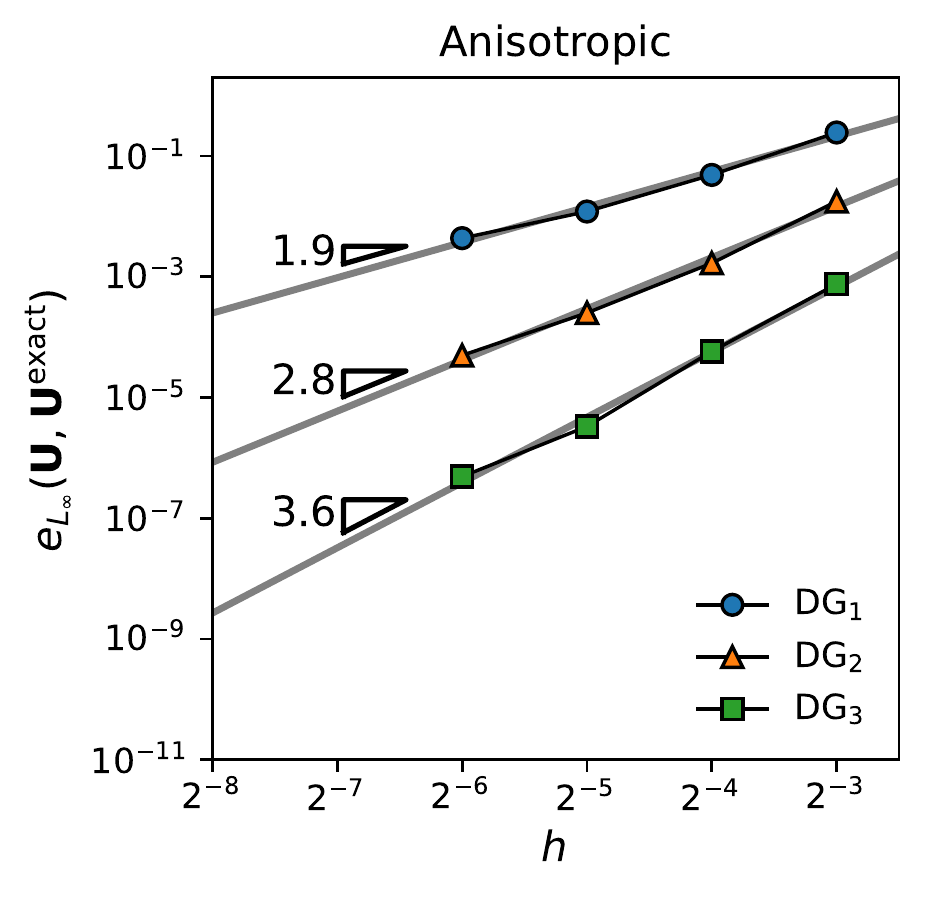}
\end{subfigure}
\\
\begin{subfigure}{0.32\textwidth}
\includegraphics[width=\textwidth]{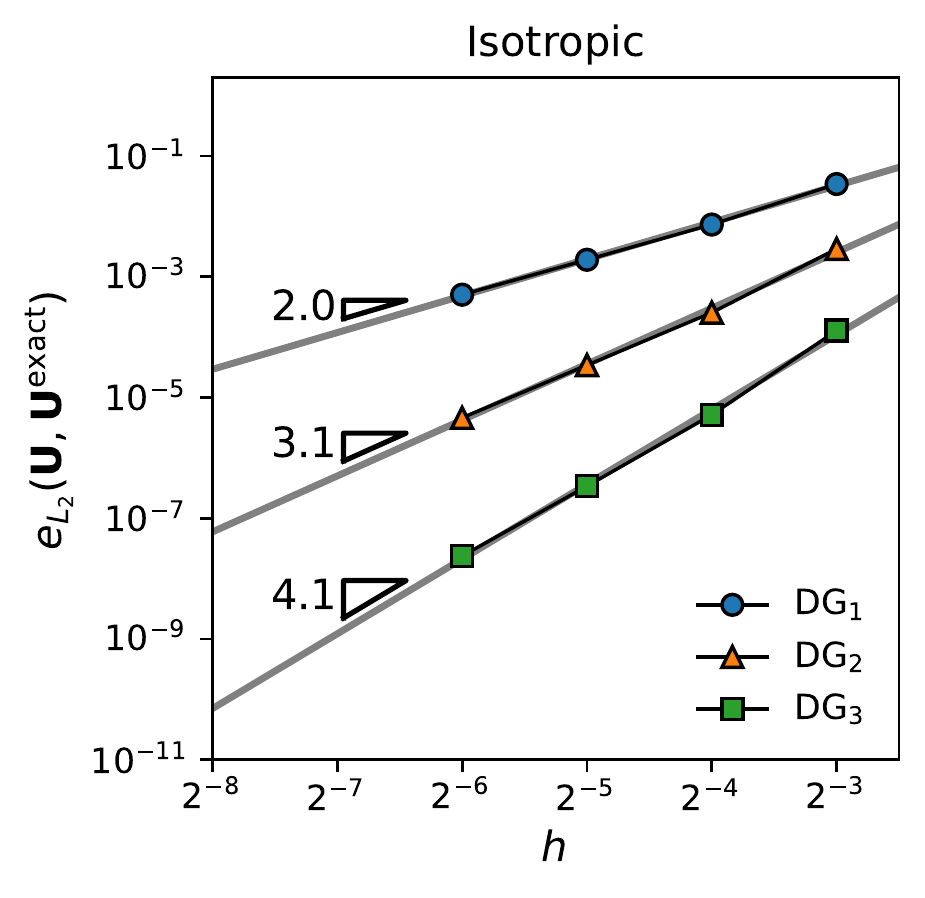}
\end{subfigure}
\
\begin{subfigure}{0.32\textwidth}
\includegraphics[width=\textwidth]{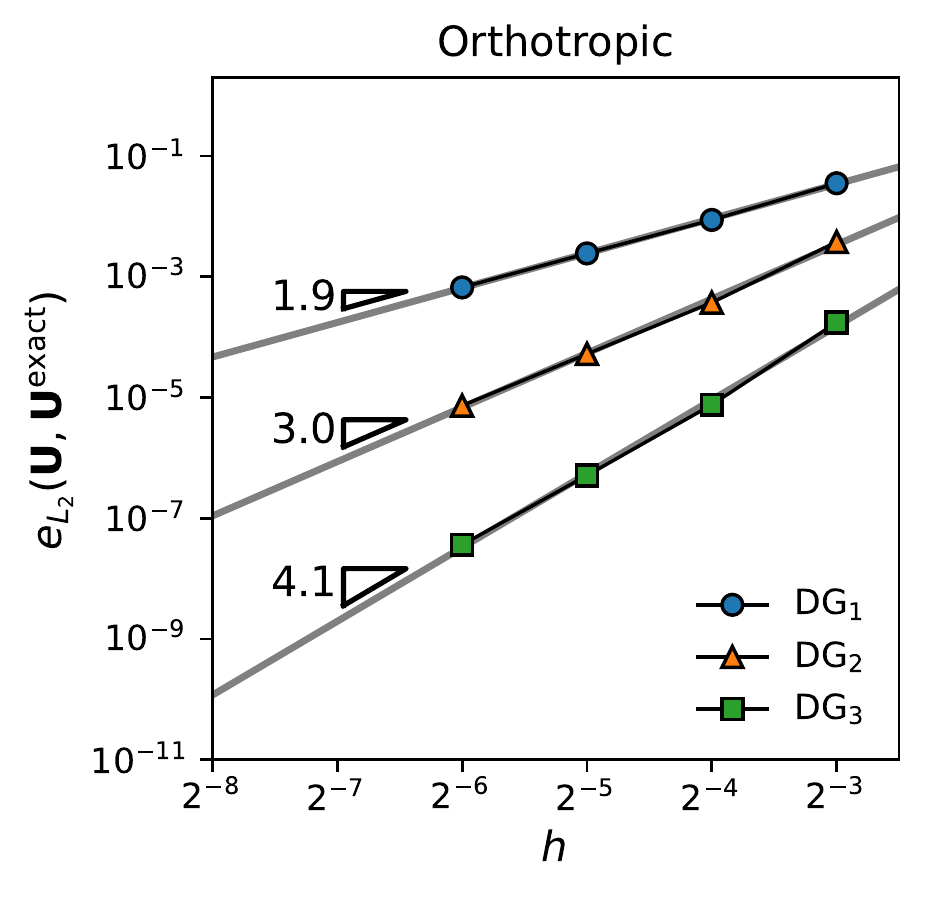}
\end{subfigure}
\
\begin{subfigure}{0.32\textwidth}
\includegraphics[width=\textwidth]{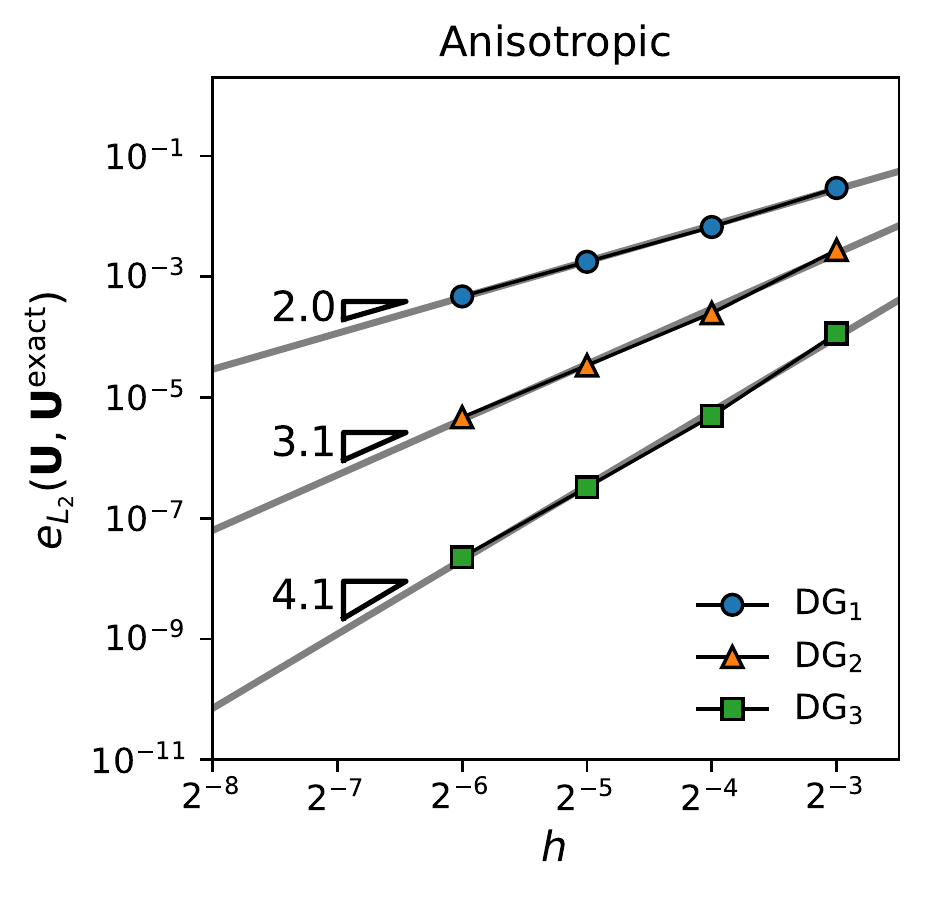}
\end{subfigure}
\caption
{
    (top row) $e_{L_\infty}$ error and (bottom row) $e_{L_2}$ error for the 3D single-phase solid of Fig.(\ref{fig:RESULTS - CONVERGENCE SINGLE PHASE DOMAIN - GEOMETRY ERROR 3D}) for different constitutive behaviors.
}
\label{fig:RESULTS - CONVERGENCE SINGLE PHASE DOMAIN - HP STUDY 3D}
\end{figure}

Figure (\ref{fig:RESULTS - CONVERGENCE SINGLE PHASE DOMAIN - GEOMETRY ERROR 3D}a) shows the geometry and the boundary conditions for the 3D single-phase solid case.
The geometry consists of a sphere implicitly-defined in the background unit cube $[0,1]^3$ by the level set function
\begin{equation}
\varphi(\bm{x}) = (x_1-o_1)^2+(x_2-o_2)^2+(x_3-o_3)^2-R^2,
\end{equation}
where $R = 0.35$ and $o_1 = o_2 = o_3 = 0.5$.
A traction field $\overline{\bm{t}}$ is prescribed on sphere's outer boundary using the exact solution $\bm{U}^\mathrm{exact}$ with $\kappa_1 = 2\pi\cos(\pi/3)\sin(\pi/6)$, $\kappa_2 = 2\pi\sin(\pi/3)\sin(\pi/6)$ and $\kappa_3 = 2\pi\cos(\pi/6)$.
Figure (\ref{fig:RESULTS - CONVERGENCE SINGLE PHASE DOMAIN - GEOMETRY ERROR 3D}a) also shows the implicitly-defined generated from an $8^3$ background grid and some of the implicitly-defined elements in proximity of the embedded boundary.
Figures (\ref{fig:RESULTS - CONVERGENCE SINGLE PHASE DOMAIN - GEOMETRY ERROR 3D}b) and (\ref{fig:RESULTS - CONVERGENCE SINGLE PHASE DOMAIN - GEOMETRY ERROR 3D}c) show the error in the momentum component $m_1$ and the strain component $\gamma_{11}$, respectively, when the implicit-mesh of Fig.(\ref{fig:RESULTS - CONVERGENCE SINGLE PHASE DOMAIN - GEOMETRY ERROR 3D}a) and a DG$_3$ scheme are employed for the anisotropic material response case.
The $hp$-convergence plots of the two error measures given in Eq.(\ref{eq:error measures}) for the 3D single-phase solid with isotropic, orthotropic and anisotropic behavior are then reported in Fig.(\ref{fig:RESULTS - CONVERGENCE SINGLE PHASE DOMAIN - HP STUDY 3D}).

\begin{figure}
\centering
\includegraphics[width = \textwidth]{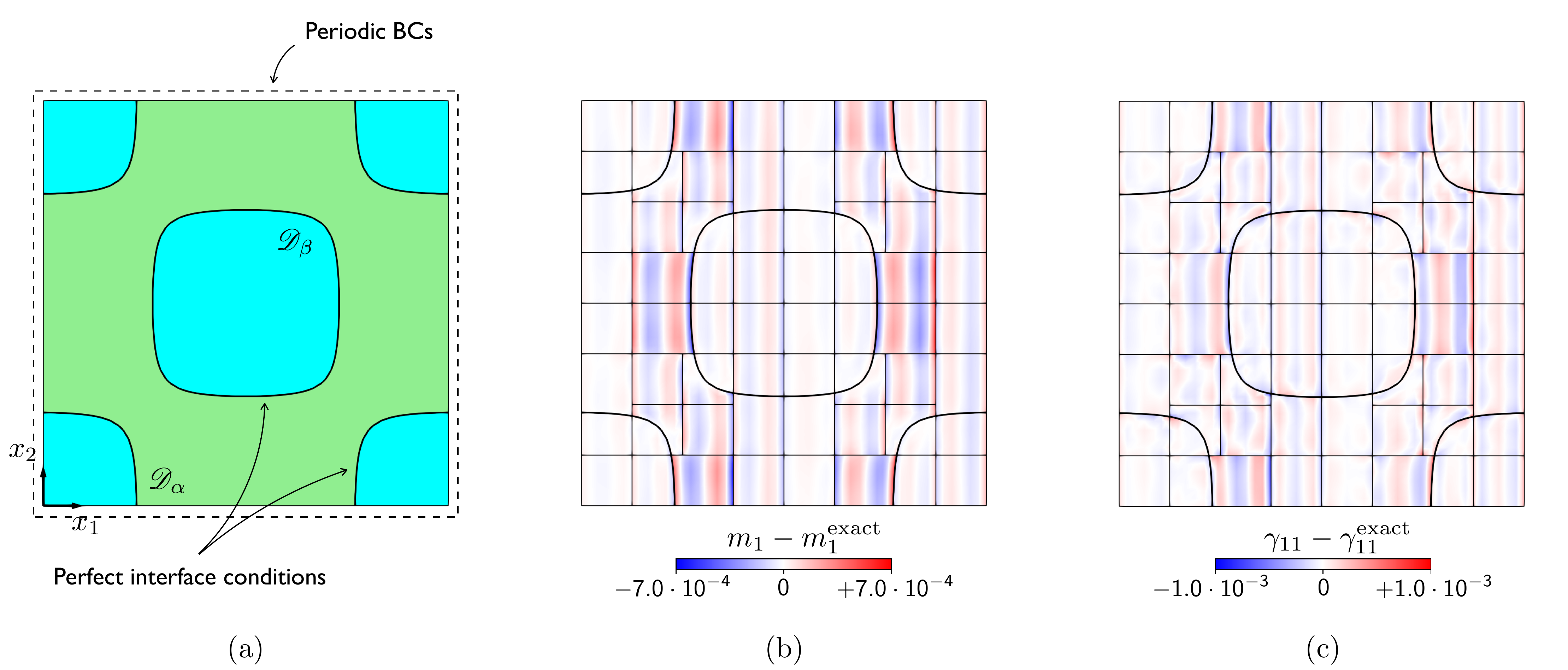}
\caption
{
    (a) Geometry and boundary conditions for the bi-phase solid employed in the 2D $hp$-convergence analysis.
    Error in (b) the momentum component $m_1$ and (c) the strain component $\gamma_{11}$ obtained with an implicitly-defined mesh generated from an 8$\times$8 background grid and a DG$_3$ scheme for the anisotropic material response.
}
\label{fig:RESULTS - CONVERGENCE TWO PHASE DOMAIN - GEOMETRY ERROR 2D}
\end{figure}

\begin{figure}
\centering
\begin{subfigure}{0.32\textwidth}
\includegraphics[width=\textwidth]{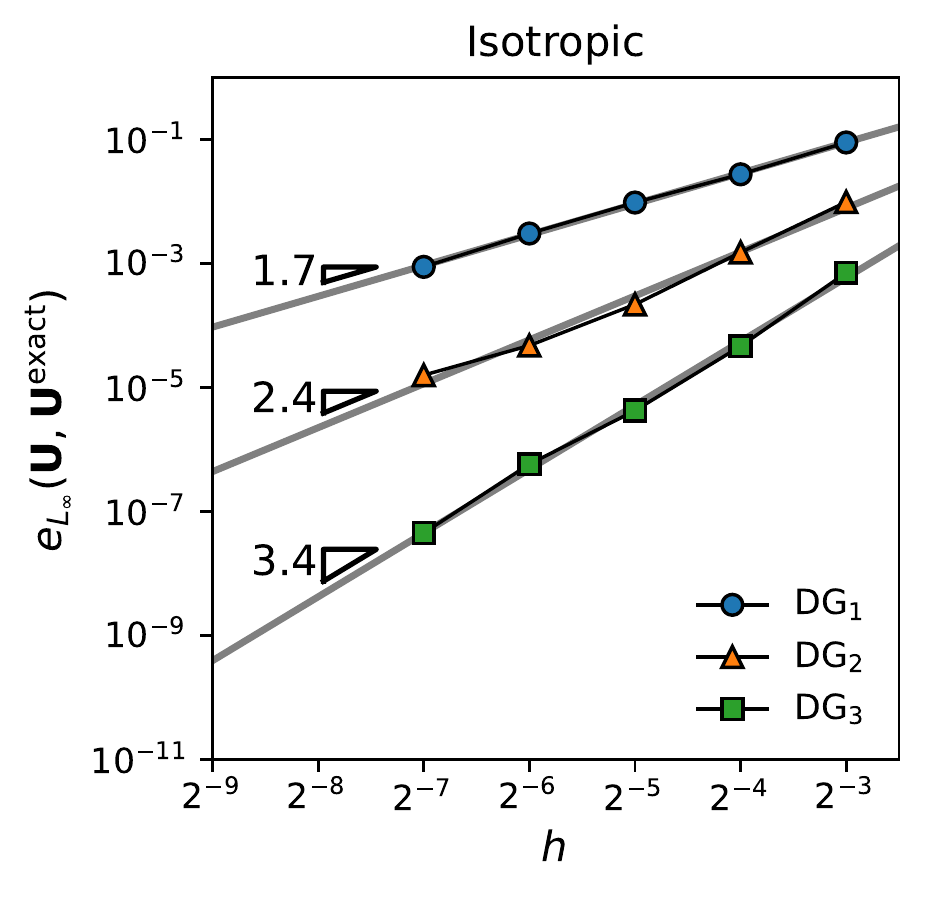}
\end{subfigure}
\
\begin{subfigure}{0.32\textwidth}
\includegraphics[width=\textwidth]{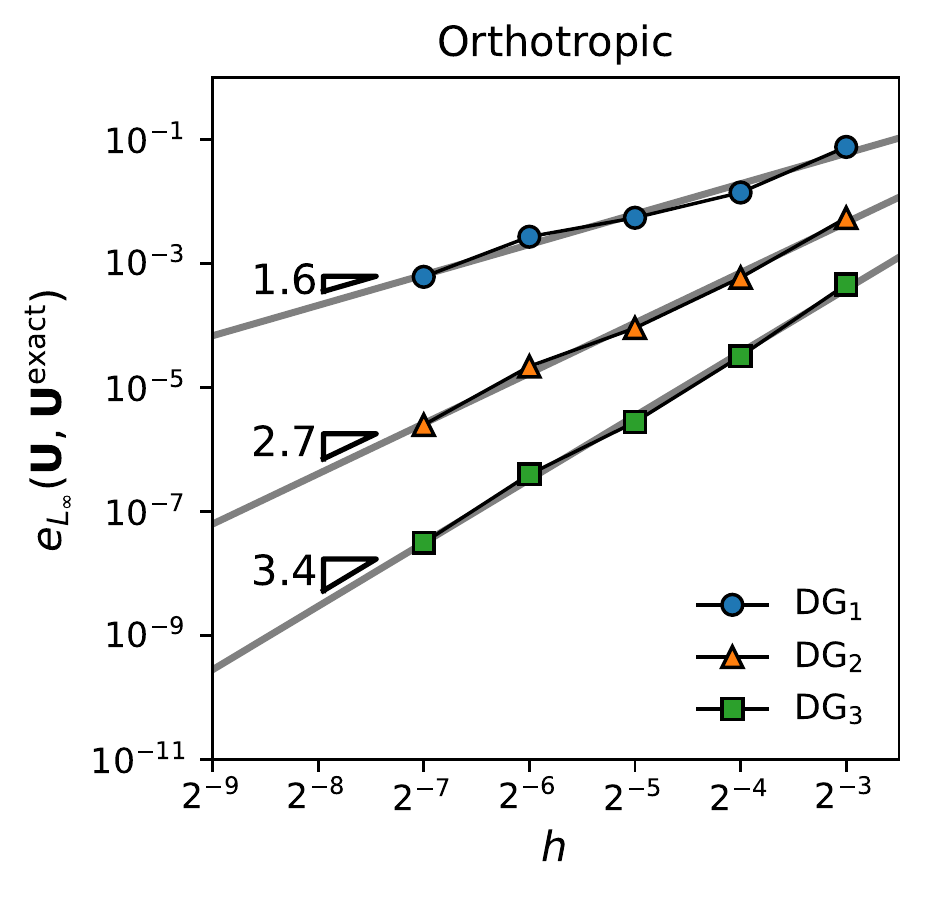}
\end{subfigure}
\
\begin{subfigure}{0.32\textwidth}
\includegraphics[width=\textwidth]{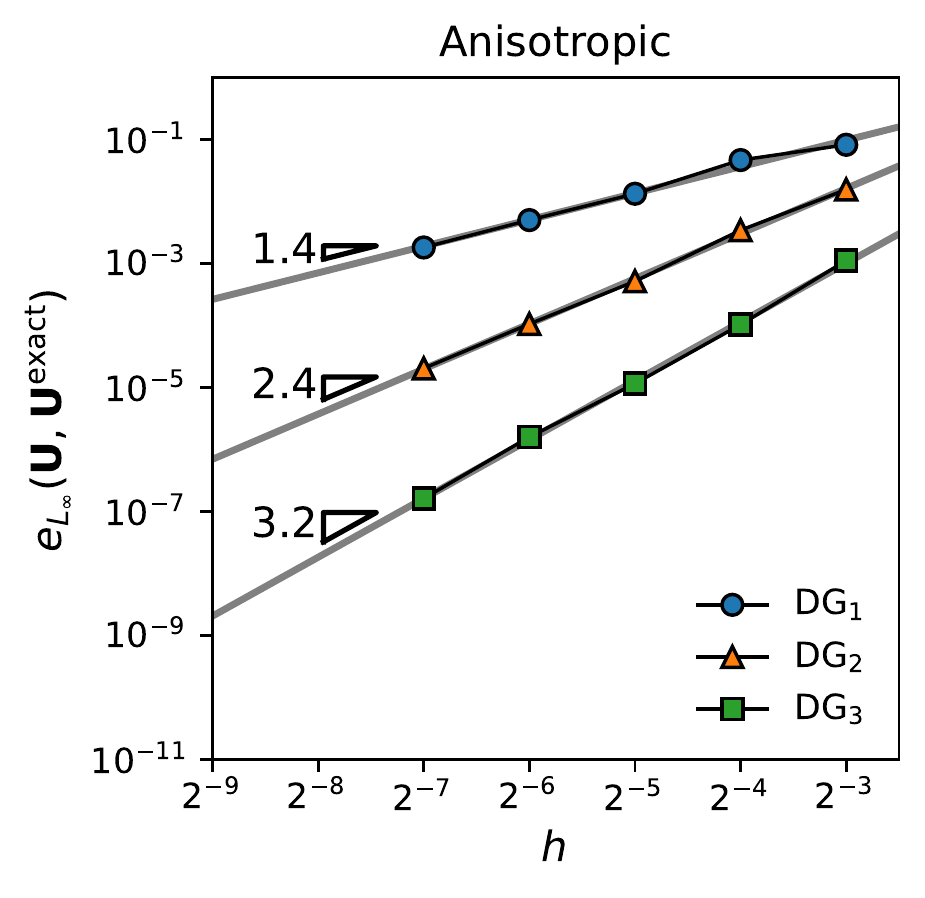}
\end{subfigure}
\\
\begin{subfigure}{0.32\textwidth}
\includegraphics[width=\textwidth]{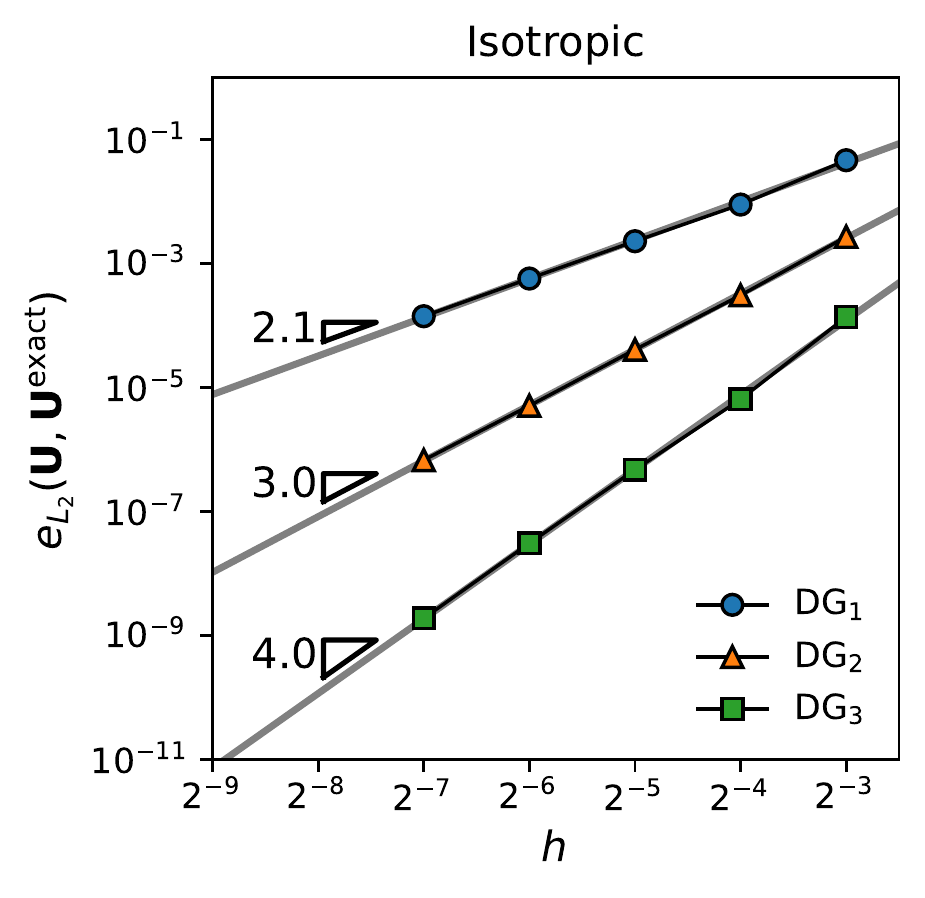}
\end{subfigure}
\
\begin{subfigure}{0.32\textwidth}
\includegraphics[width=\textwidth]{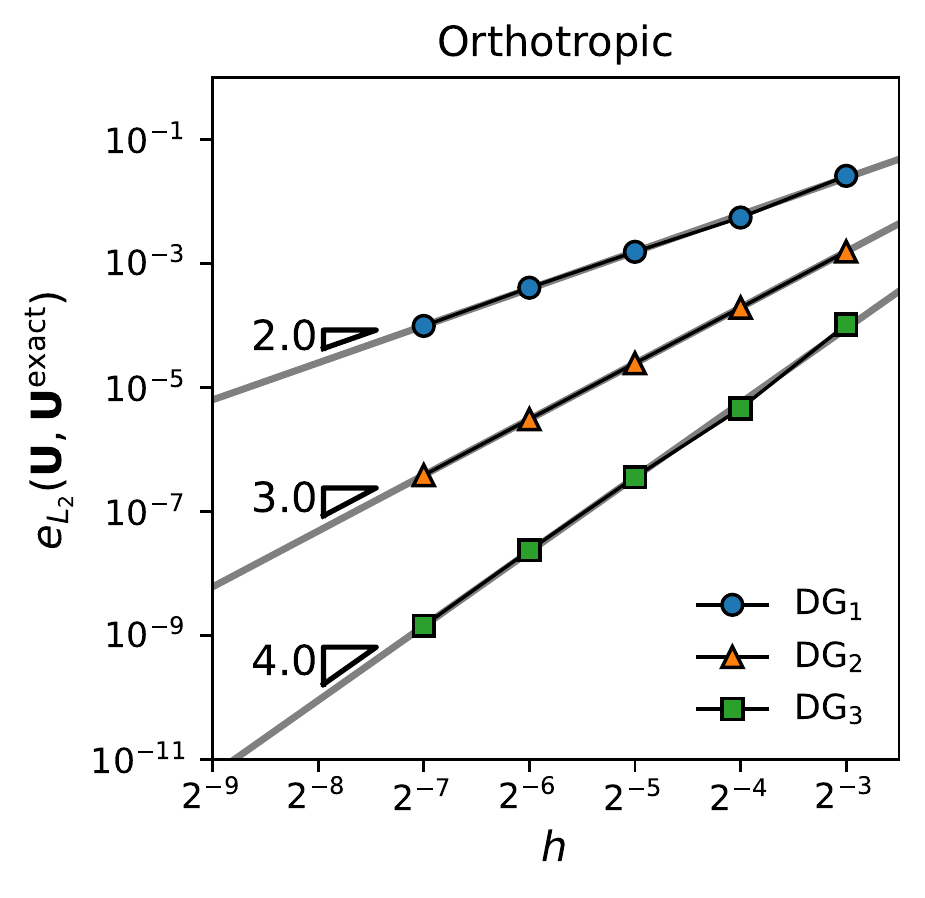}
\end{subfigure}
\
\begin{subfigure}{0.32\textwidth}
\includegraphics[width=\textwidth]{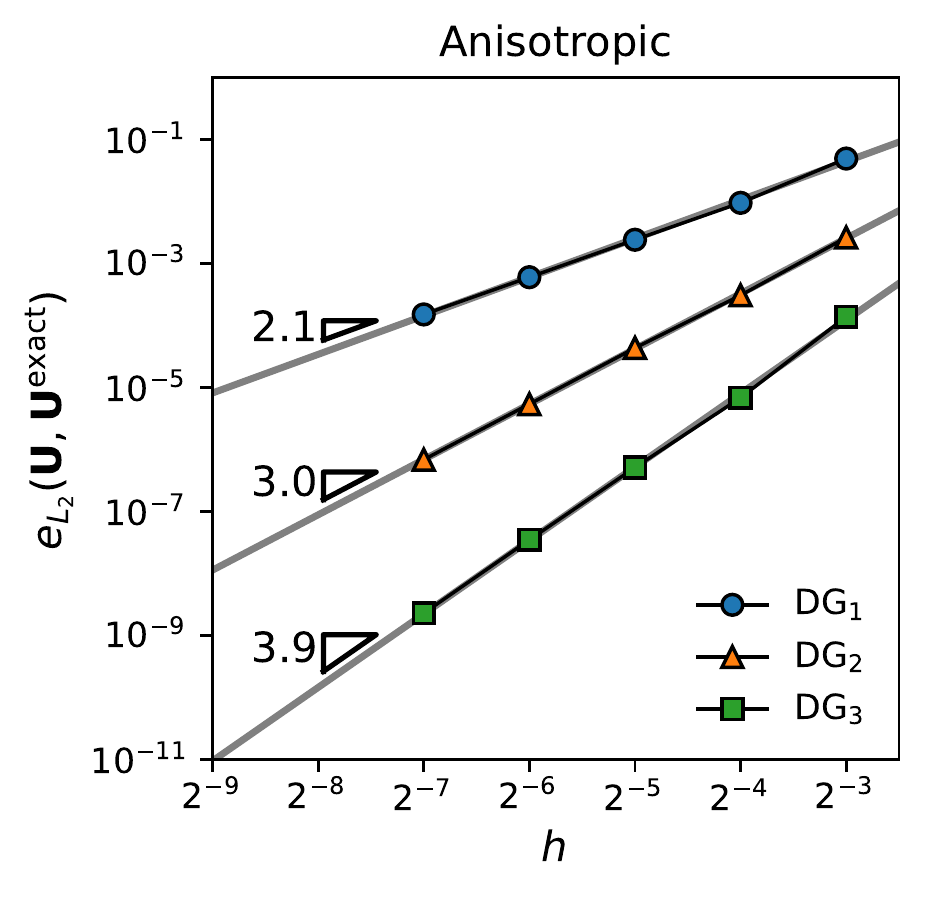}
\end{subfigure}
\caption
{
    (top row) $e_{L_\infty}$ error and (bottom row) $e_{L_2}$ error for the 2D single-phase solid of Fig.(\ref{fig:RESULTS - CONVERGENCE TWO PHASE DOMAIN - GEOMETRY ERROR 2D}) for different constitutive behaviors.
}
\label{fig:RESULTS - CONVERGENCE TWO PHASE DOMAIN - HP STUDY 2D}
\end{figure}

Figure (\ref{fig:RESULTS - CONVERGENCE TWO PHASE DOMAIN - GEOMETRY ERROR 2D}a) shows the geometry and the boundary conditions for the 2D two-phase solid case.
The geometry is periodic and is defined in the background unit square $[0,1]^2$ by the level set function
\begin{equation}
\varphi(\bm{x}) = \cos(2 \pi x_1)\cos(2 \pi x_2)-1/8.
\end{equation}
Periodic boundary conditions are prescribed on the outer boundaries $\mathscr{B}_\alpha$ and $\mathscr{B}_\beta$ of the background square, whereas perfect interface conditions as given in Eq.(\ref{eq:interface conditions}) are prescribed on $\mathscr{L}_{\alpha,\beta}$.
The exact solution $\bm{U}^\mathrm{exact}$ is specified by $\kappa_1 = 2\pi$ and $\kappa_2 = 0$.
Figures (\ref{fig:RESULTS - CONVERGENCE TWO PHASE DOMAIN - GEOMETRY ERROR 2D}b) and (\ref{fig:RESULTS - CONVERGENCE TWO PHASE DOMAIN - GEOMETRY ERROR 2D}c) show the error in the momentum component $m_1$ and the strain component $\gamma_{11}$, respectively, when an $8^2$ background grid and a DG$_3$ scheme are employed for the anisotropic material response case.
The obtained $hp$-convergence plots for the isotropic, orthotropic and anisotropic behavior are reported in Fig.(\ref{fig:RESULTS - CONVERGENCE TWO PHASE DOMAIN - HP STUDY 2D}).

\begin{figure}
\centering
\includegraphics[width = \textwidth]{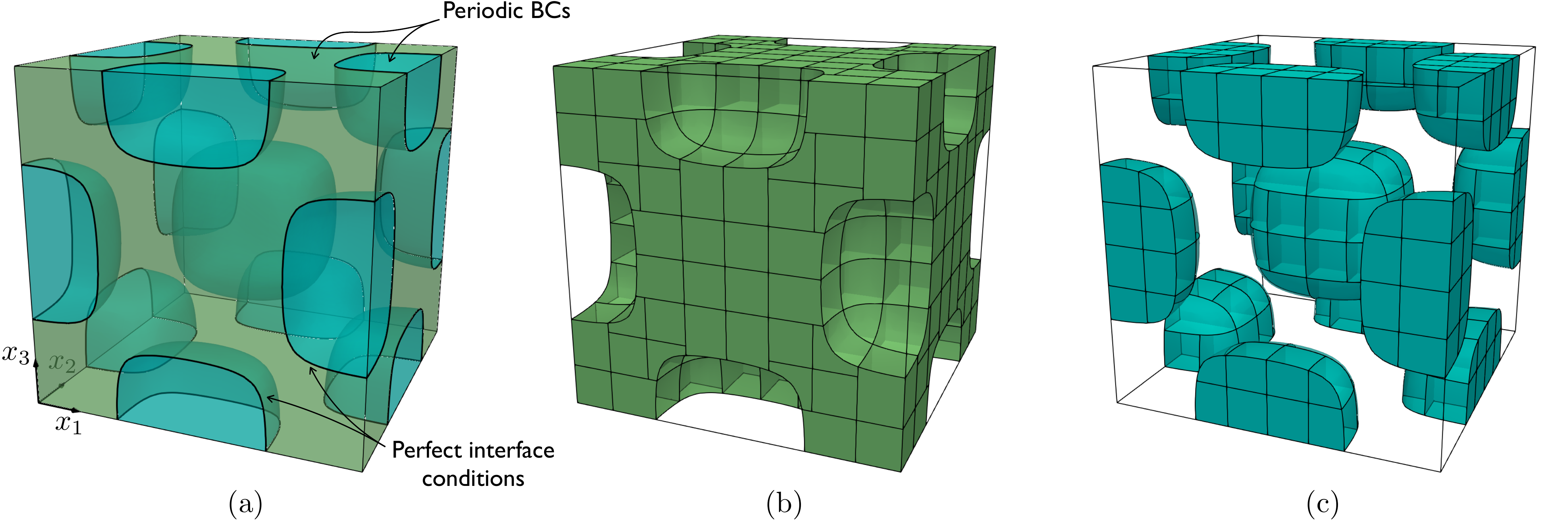}
\caption
{
    (a) Geometry and boundary conditions for the bi-phase solid employed in the 3D $hp$-convergence analysis.
    Implicitly-defined mesh of (b) the phase $\alpha$ and (c) the phase $\beta$ generated from an 8$\times$8$\times$8 background grid for the solid of figure (a). 
}
\label{fig:RESULTS - CONVERGENCE TWO PHASE DOMAIN - GEOMETRY MESH 3D}
\end{figure}

\begin{figure}
\centering
\includegraphics[width = \textwidth]{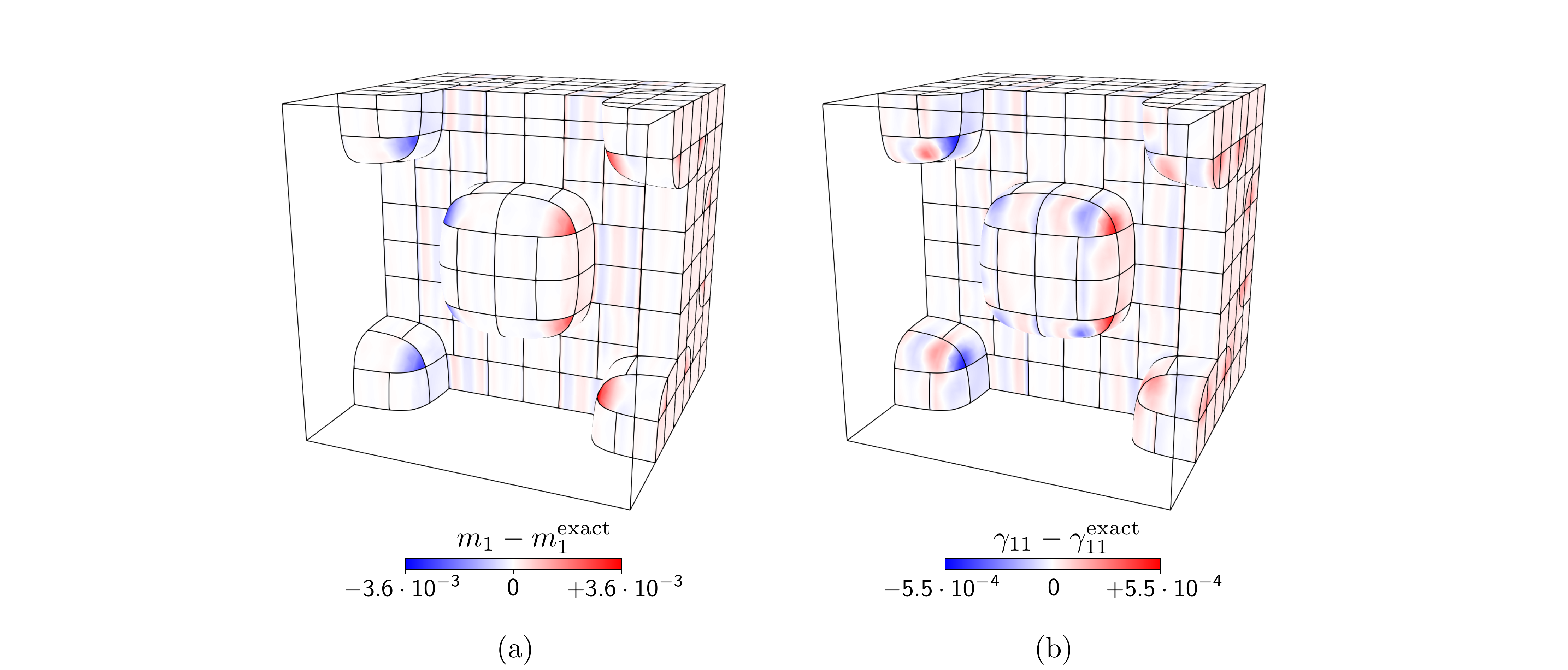}
\caption
{
    Error in (a) the momentum component $m_1$ and (b) the strain component $\gamma_{11}$ obtained with the implicitly-defined mesh of Fig.(\ref{fig:RESULTS - CONVERGENCE TWO PHASE DOMAIN - GEOMETRY MESH 3D}) and a DG$_3$ scheme for the anisotropic material response.
}
\label{fig:RESULTS - CONVERGENCE TWO PHASE DOMAIN - ERROR 3D}
\end{figure}

\begin{figure}
\centering
\begin{subfigure}{0.32\textwidth}
\includegraphics[width=\textwidth]{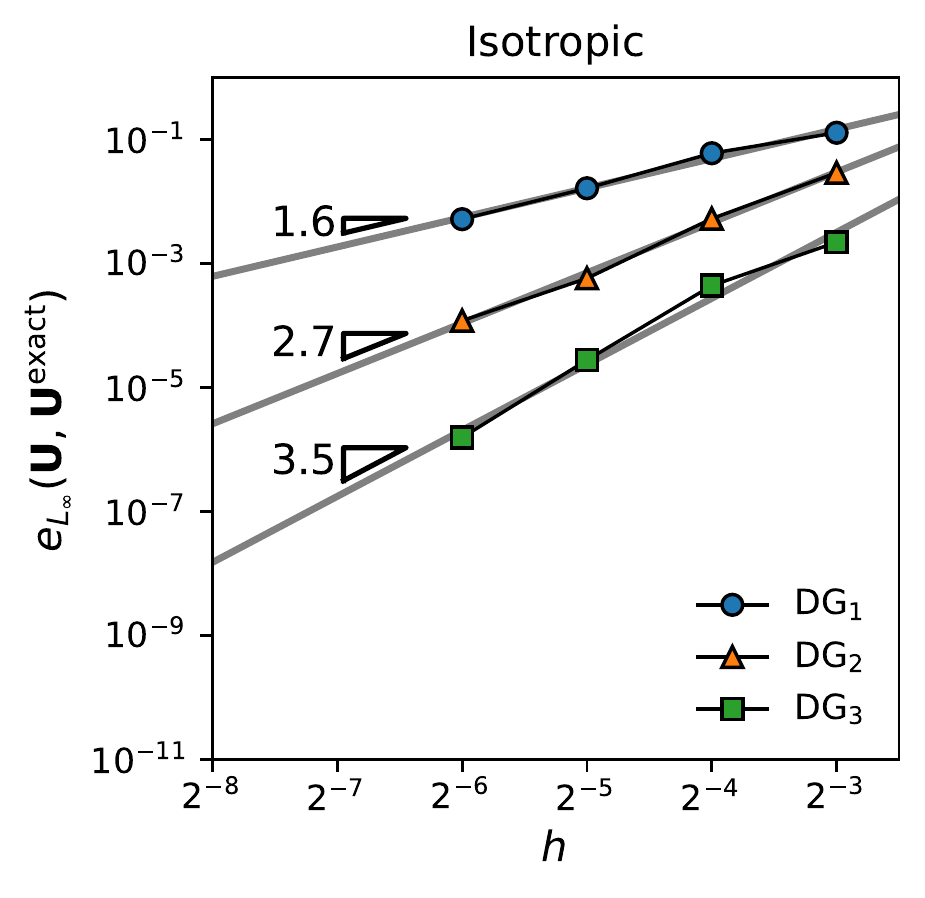}
\end{subfigure}
\
\begin{subfigure}{0.32\textwidth}
\includegraphics[width=\textwidth]{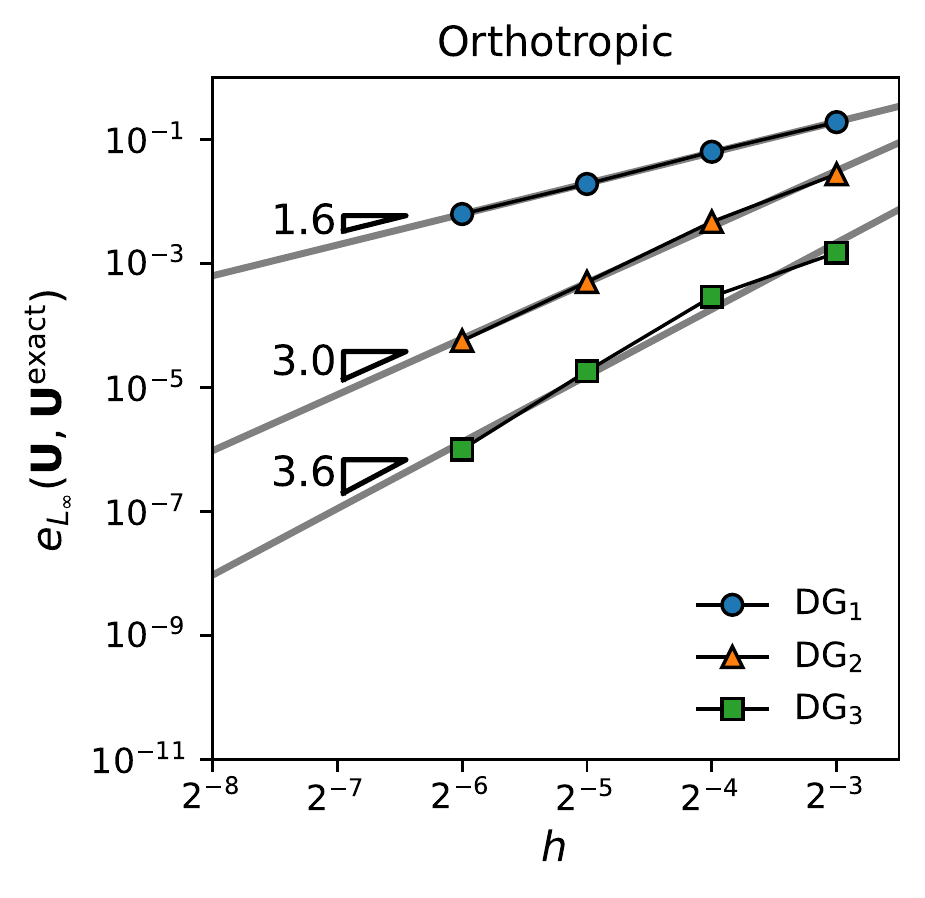}
\end{subfigure}
\
\begin{subfigure}{0.32\textwidth}
\includegraphics[width=\textwidth]{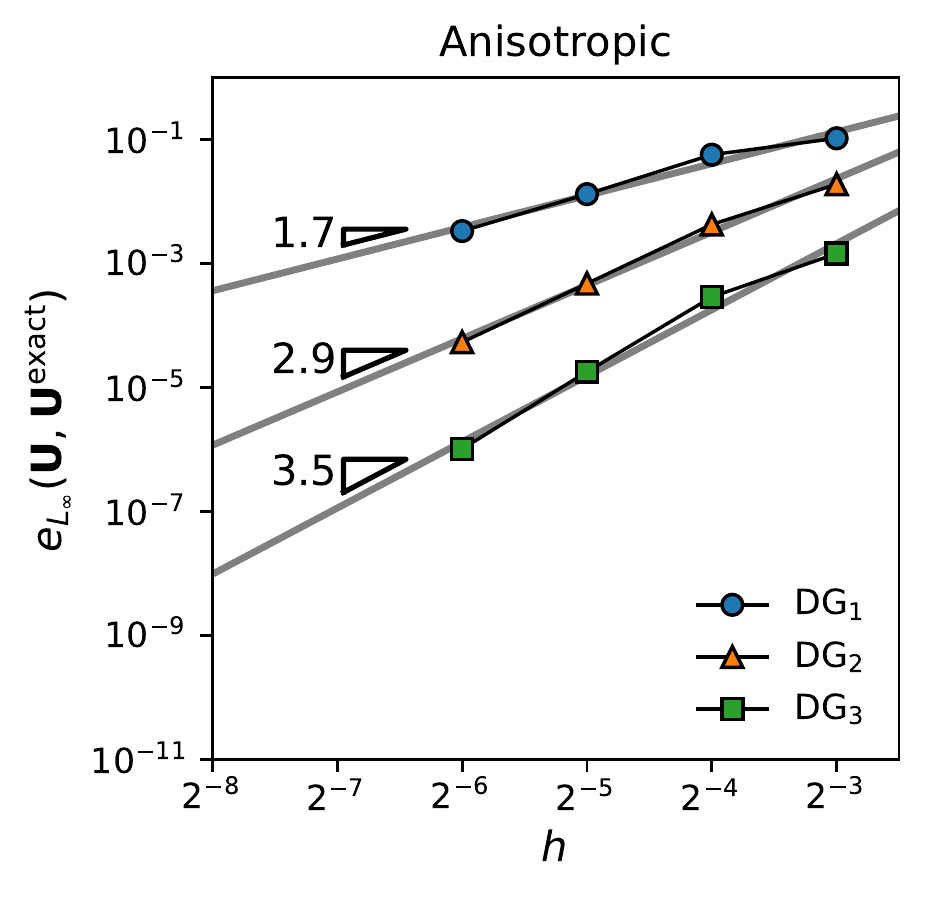}
\end{subfigure}
\\
\begin{subfigure}{0.32\textwidth}
\includegraphics[width=\textwidth]{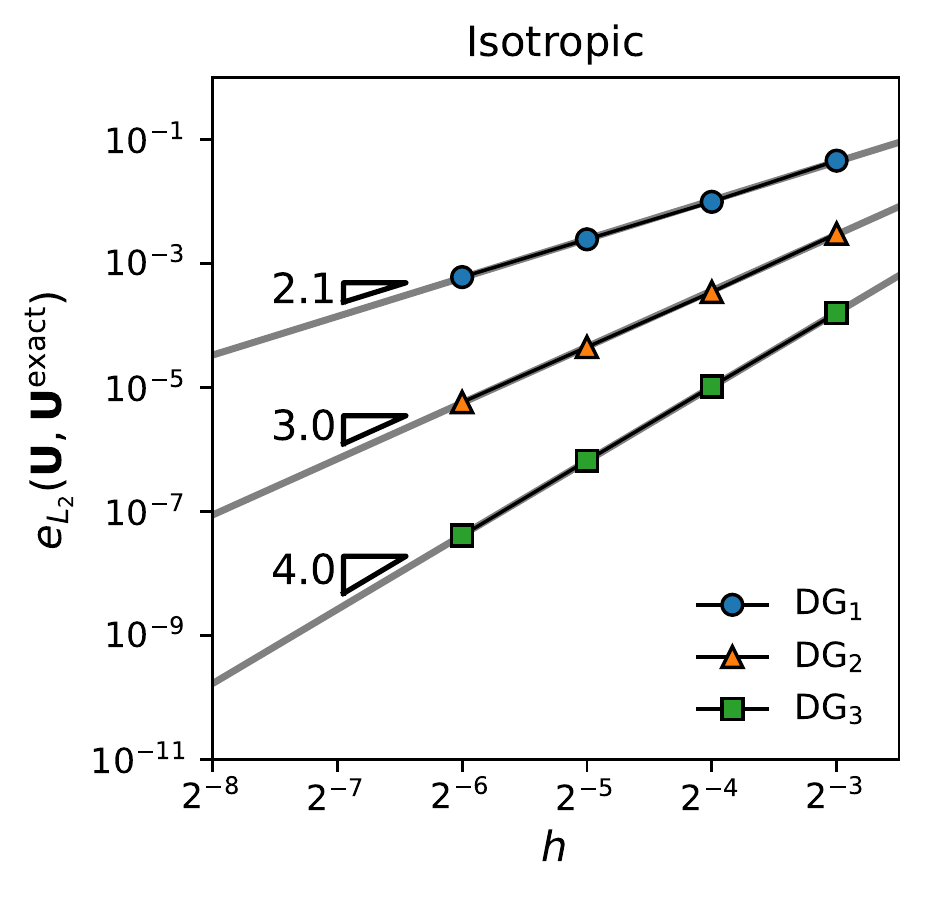}
\end{subfigure}
\
\begin{subfigure}{0.32\textwidth}
\includegraphics[width=\textwidth]{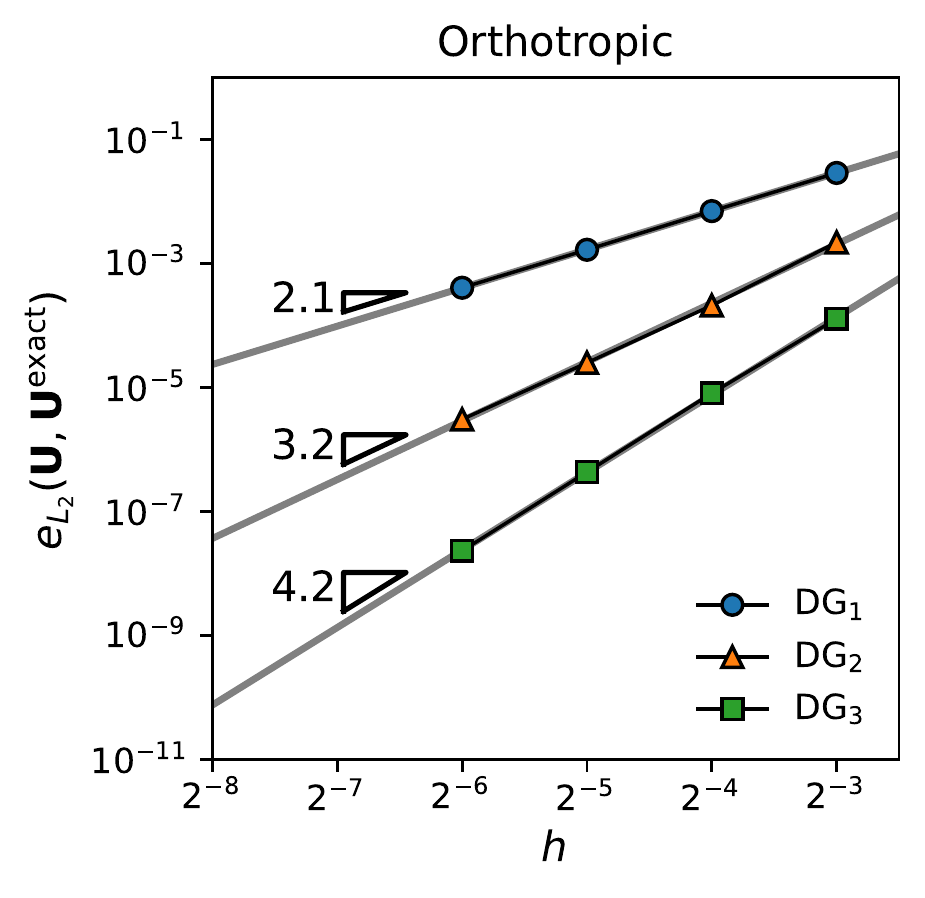}
\end{subfigure}
\
\begin{subfigure}{0.32\textwidth}
\includegraphics[width=\textwidth]{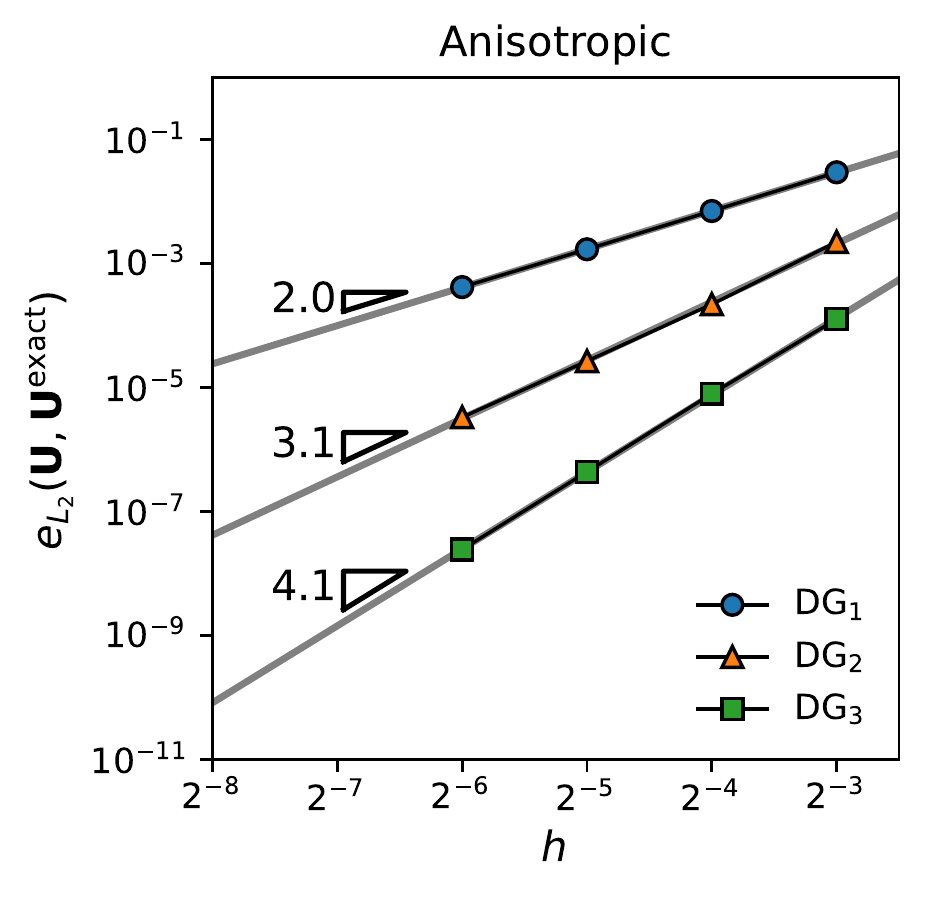}
\end{subfigure}
\caption
{
    (top row) $e_{L_\infty}$ error and (bottom row) $e_{L_2}$ error for the 3D single-phase solid of Fig.(\ref{fig:RESULTS - CONVERGENCE TWO PHASE DOMAIN - GEOMETRY MESH 3D}) for different constitutive behaviors.
}
\label{fig:RESULTS - CONVERGENCE TWO PHASE DOMAIN - HP STUDY 3D}
\end{figure}

Figure (\ref{fig:RESULTS - CONVERGENCE TWO PHASE DOMAIN - GEOMETRY MESH 3D}a) shows the geometry and the boundary conditions for the 3D two-phase solid case.
The geometry is periodic and is defined in the background unit cube $[0,1]^3$ by the level set function
\begin{equation}
\varphi(\bm{x}) = -\cos(2 \pi x_1)\cos(2 \pi x_2)\cos(2 \pi x_3)-1/8.
\end{equation}
Periodic boundary conditions are prescribed on the outer boundaries $\mathscr{B}_\alpha$ and $\mathscr{B}_\beta$ of the background square, whereas perfect interface conditions are prescribed on $\mathscr{L}_{\alpha,\beta}$.
The exact solution $\bm{U}^\mathrm{exact}$ is specified by $\kappa_1 = 2\pi$ and $\kappa_2 = \kappa_3 = 0$.
Figures (\ref{fig:RESULTS - CONVERGENCE TWO PHASE DOMAIN - GEOMETRY MESH 3D}b) and (\ref{fig:RESULTS - CONVERGENCE TWO PHASE DOMAIN - GEOMETRY MESH 3D}c) display the implicitly-defined mesh of the phase $\alpha$ and the phase $\beta$, respectively, generated from an $8^3$ background grid.
For this mesh, Figs.(\ref{fig:RESULTS - CONVERGENCE TWO PHASE DOMAIN - ERROR 3D}a) and (\ref{fig:RESULTS - CONVERGENCE TWO PHASE DOMAIN - ERROR 3D}b) show the error in the momentum component $m_1$ and the strain component $\gamma_{11}$, respectively, when a DG$_3$ scheme is employed for the anisotropic material response case.
Finally, the $hp$-convergence plots for the isotropic, orthotropic and anisotropic behavior are reported in Fig.(\ref{fig:RESULTS - CONVERGENCE TWO PHASE DOMAIN - HP STUDY 3D}).

To conclude this part of the numerical results, we observe that in all simulations the error between the exact solution and the numerical solution in the extended elements is mildly larger than the error in the regular (hyper)rectangular elements, see for example Fig.(\ref{fig:RESULTS - CONVERGENCE SINGLE PHASE DOMAIN - GEOMETRY ERROR 3D}c) or Fig.(\ref{fig:RESULTS - CONVERGENCE TWO PHASE DOMAIN - ERROR 3D}c).
This has also been observed for other applications of the present implicit-mesh DG approaches \cite{saye2017implicitI,saye2017implicitII,gulizzi2021coupled} and is an expected behavior if one considers that the extended elements are in general larger and geometrically less regular than the (hyper)rectangular elements.
However, as shown by all the $hp$-convergence plots, and consistently with previous observations \cite{saye2017implicitI,saye2017implicitII,gulizzi2021coupled}, the present implicit-mesh DG method provides a high-order accurate solution of the elastodynamics problem demonstrated by a $\mathcal{O}(h^{p+a})$ convergence rate in the energy norm, with $a$ ranging from $0.8$ to $1.2$, and by a $\mathcal{O}(h^{p+a})$ convergence rate in the $L_\infty$ norm, with $a$ ranging from $0.1$ to $1.1$.

\subsection{Lamb's problem}\label{ssec:Lamb problem}

\begin{figure}
\centering
\includegraphics[width = \textwidth]{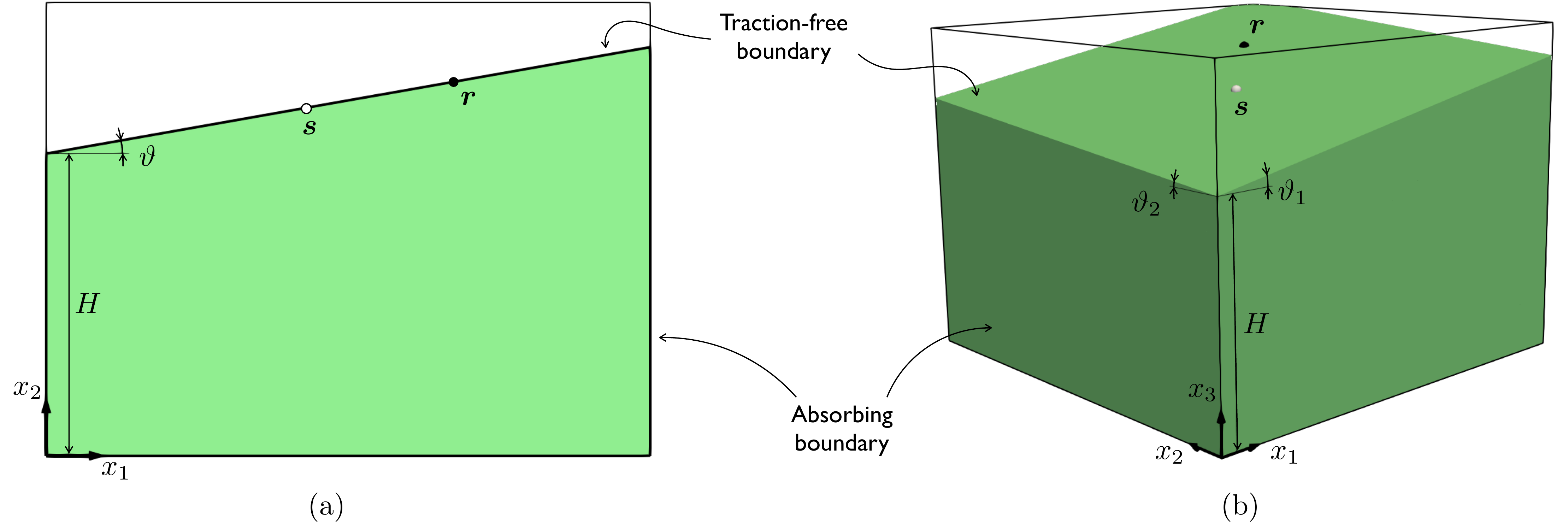}
\caption
{
    Geometry and boundary conditions for (a) the 2D and (b) the 3D Lamb's problem.
    In the figures, $\bm{s}$ denotes the location of the source point where the concentrated force is applied whereas $\bm{r}$ denotes the location of the receiver where the elastodynamic response is measured.
}
\label{fig:RESULTS - LAMB PROBLEM - GEOMETRY}
\end{figure}

In this section, we consider a classical problem in elastodynamics, namely the Lamb's problem, which admits an exact solution \cite{berg1994analytical,kausel2013lamb} and has been used to assess various numerical models, see e.g.~Refs.\cite{komatitsch1998spectral,kaser2006arbitrary,tavelli2019simple}.
The problem consists of evaluating the distribution of the mechanical fields due to a concentrated force that is applied perpendicular to an infinite free surface, i.e.~with zero-traction boundary conditions.
The well-known distinctive feature of this problem is the appearance of the Rayleigh waves, which travel along the free surface of the domain but not in its depth.

Following the problem setup of Refs.\cite{komatitsch1998spectral,kaser2006arbitrary,tavelli2019simple}, the 2D geometry is implicitly-defined in the background rectangle $\mathscr{R} = [0,4000\, \mathrm{m}]\times[0,3000\, \mathrm{m}]$ by the level set function
\begin{equation}\label{eq:Lamb problem 2D - level set}
\varphi(\bm{x}) = x_2-H-\tan(\vartheta)x_1,
\end{equation}
where $H = 2000\, \mathrm{m}$ and $\vartheta = 10^\circ$.
We note that the problem may be stated in a reference system that is aligned with the free surface and the elastic domain may be straightforwardly meshed with a simple structured grid without involving implicitly-defined elements; nevertheless, this setup is a common benchmark problem, including to assess EB approaches, see, e.g., Ref.\cite{tavelli2019simple}.
In 3D, we consider a simple extension of the 2D case, whereby the geometry is implicitly-defined in the background prism $\mathscr{R} = [0,4000\, \mathrm{m}]\times[0,4000\, \mathrm{m}]\times[0,3000\, \mathrm{m}]$ by the level set function
\begin{equation}\label{eq:Lamb problem 3D - level set}
\varphi(\bm{x}) = x_3-H-\tan(\vartheta_1)x_1-\tan(\vartheta_2)x_2,
\end{equation}
where $H = 2000\, \mathrm{m}$, $\vartheta_1 = 10^\circ$ and $\vartheta_1 = 5^\circ$.
The 2D geometry and the 3D geometry are displayed in Fig.(\ref{fig:RESULTS - LAMB PROBLEM - GEOMETRY}a) and (\ref{fig:RESULTS - LAMB PROBLEM - GEOMETRY}b), respectively.
In both figures, the point $\bm{s}$ denotes the location of the force and the point $\bm{r}$ denotes the location of a receiver where the resulting elastodynamic response is measured; in 2D, the source point and the receiver are located at $\bm{s} = (s_1,s_2)$ and $\bm{r} = (r_1,r_2)$, where $s_1 = 1720\, \mathrm{m}$, $r_1 = 2694.96\, \mathrm{m}$ and $s_2$ and $r_2$ are obtained via Eq.(\ref{eq:Lamb problem 2D - level set}); in 3D, the source point and the receiver are located at $\bm{s} = (s_1,s_2,s_3)$ and $\bm{r} = (r_1,r_2,r_4)$, where $s_1 = s_2 = 1720\, \mathrm{m}$, $r_1 = r_2 = 2694.96\, \mathrm{m}$ and $s_3$ and $r_3$ are obtained via Eq.(\ref{eq:Lamb problem 3D - level set}).
The elastic domain is an isotropic solid with density $\rho_\alpha = 2200.0\, \mathrm{kg}/\mathrm{m}^3$ and elastic properties determined by the velocity of the P-waves $c_P = 3200.0\, \mathrm{m}/\mathrm{s}$ and the velocity of the S-waves $c_S = 1847.5\, \mathrm{m}/\mathrm{s}$.
Traction-free boundary conditions, i.e. $\overline{\bm{t}} = \bm{0}$, are prescribed on $\mathscr{L}_\alpha$ whereas absorbing boundary conditions are prescribed on $\mathscr{B}_\alpha$.
The final time of the simulation is $T = 1\, \mathrm{s}$.
Finally, the concentrated force is modelled by setting $\rho_\alpha\bm{b}_\alpha$ in the source term $\bm{S}_\alpha$, see Eq.(\ref{eq:governing equations - vectors}), as 
\begin{equation}\label{eq:Lamb problem - external source}
\rho_\alpha\bm{b}_\alpha = R_w(t)\delta(\bm{x}-\bm{s})\hat{\bm{u}},
\end{equation}
where $\hat{\bm{u}}$ is the unit vector perpendicular to $\mathscr{L}_\alpha$, $\delta(\bm{x})$ is the Dirac delta function, and $R_w(t)$ is the Ricker wavelet defined as
\begin{equation}
R_w(t) \equiv a_1 \left(\frac{1}{2}+a_2(t-t_0)^2\right) e^{a_2(t-t_0)^2}
\end{equation}
being $a_1 = -2000\, \mathrm{kg}/(\mathrm{m}^2\mathrm{s}^2)$, $a_2 = -\pi^2 f_c^2$, $f_c = 14.5\, \mathrm{Hz}$, $t_0 = 0.08\, \mathrm{s}$.

In 2D, we consider three implicitly-defined meshes generated from uniform grids and an implicitly-defined mesh associated with a two-level $hp$-AMR scheme.
The three uniform meshes are generated from a $64\times 48$ grid, a $128\times 96$ grid and a $256\times 192$ grid, and use a DG$_3$ scheme.
For the AMR test problem, level $\ell = 0$ uses a $64\times 48$ grid and a DG$_1$ scheme, whereas level $\ell = 1$ uses a DG$_3$ scheme and is dynamically generated from level $\ell = 0$ via a refinement ratio $r = 4$; it follows that the $hp$-AMR has the same effective resolution of the finest uniform mesh.
To evolve the cell tagging, the function $f_\mathrm{tag}$ introduced in Eq.(\ref{eq:tagging}) implements the following energy-based threshold as
\begin{equation}\label{eq:tagging - Lamb problem}
f_\mathrm{tag} \equiv E_\alpha^e-E_0,
\end{equation}
where $E_\alpha^e$ denotes the energy associated with the solution $\bm{U}_\alpha^e$ for the element $\mathscr{D}_\alpha^e$, i.e.~it is computed via Eq.(\ref{eq:energy}) where $\mathscr{D}_\alpha$ is replaced by $\mathscr{D}_\alpha^e$, and $E_0$ is a threshold value chosen to be $E_0 = 10^{-12}\, \mathrm{J}$.
We remark that much more sophisticated refinement/coarsening criteria for DG methods exist in the literature, see for example Ref.\cite{zanotti2015space} where the authors use a criterion involving the DG solution and its derivatives in space; in this work we use a simpler criterion tuned such that the propagating waves are resolved by the DG$_3$ scheme while the remaining parts of the domain are resolved by the DG$_1$ scheme.

\begin{figure}
\centering
\begin{subfigure}{0.48\textwidth}
\includegraphics[width=\textwidth]{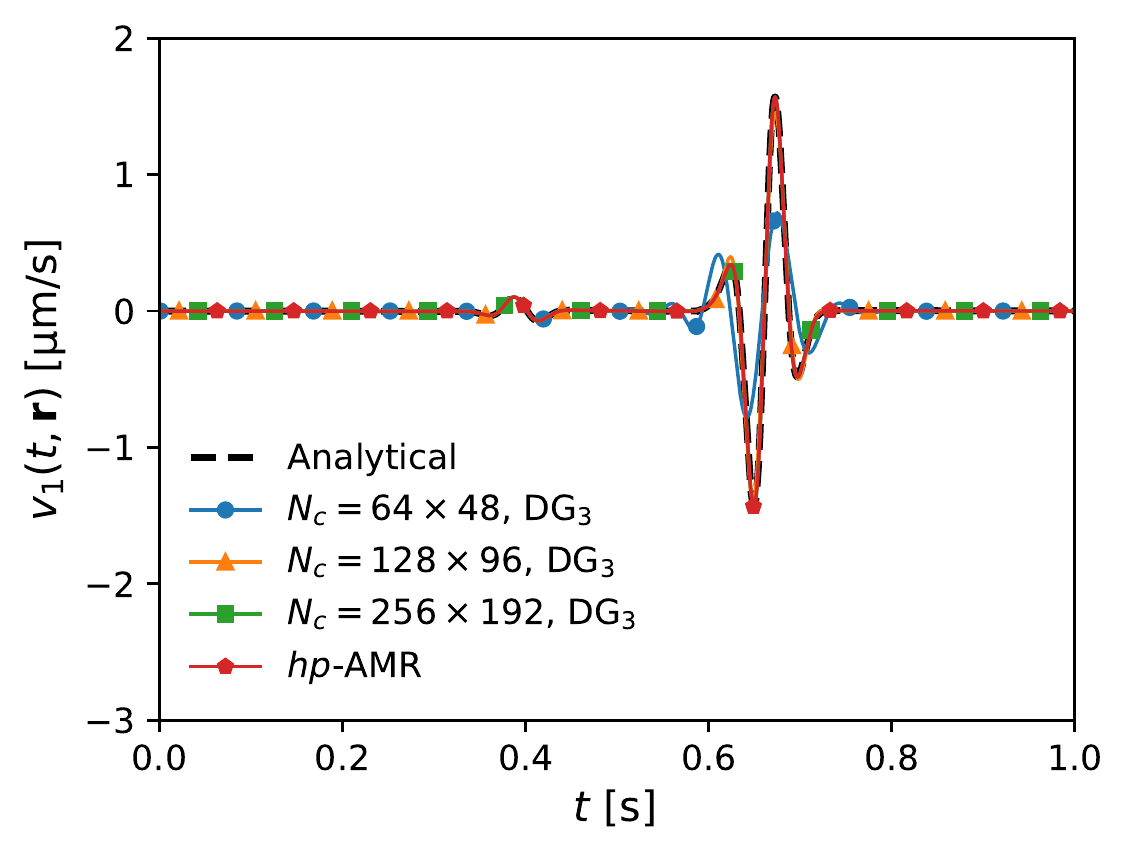}
\caption{}
\end{subfigure}
\
\begin{subfigure}{0.48\textwidth}
\includegraphics[width=\textwidth]{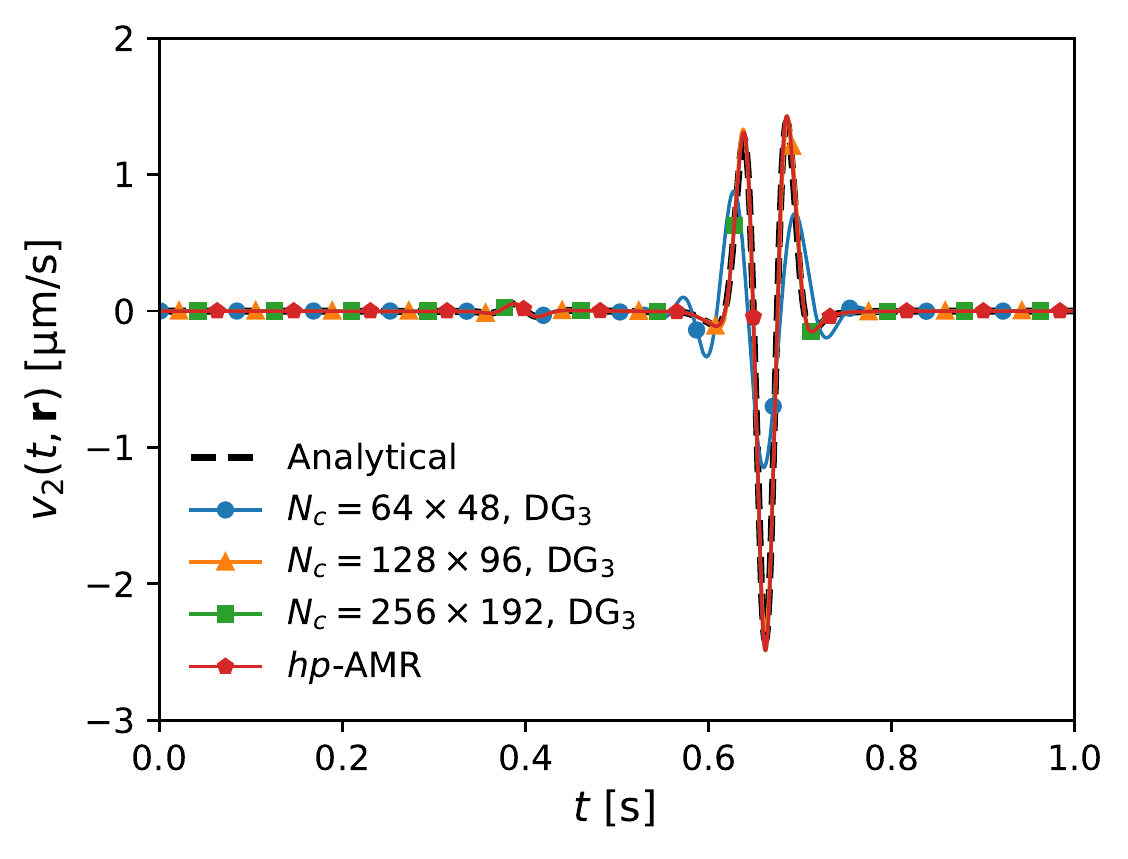}
\caption{}
\end{subfigure}
\caption
{
    Comparison between the analytical solution \cite{berg1994analytical} and the DG solution for the 2D Lamb's problem in terms of (a) horizontal velocity component $v_1$ and (b) vertical velocity component $v_2$ evaluated at the surface receiver $\bm{r}$ denoted by the black dot in Fig.(\ref{fig:RESULTS - LAMB PROBLEM - GEOMETRY}a).
}
\label{fig:RESULTS - LAMB PROBLEM - SOLUTION 2D - VELOCITY}
\end{figure}

\begin{figure}
\centering
\includegraphics[width = \textwidth]{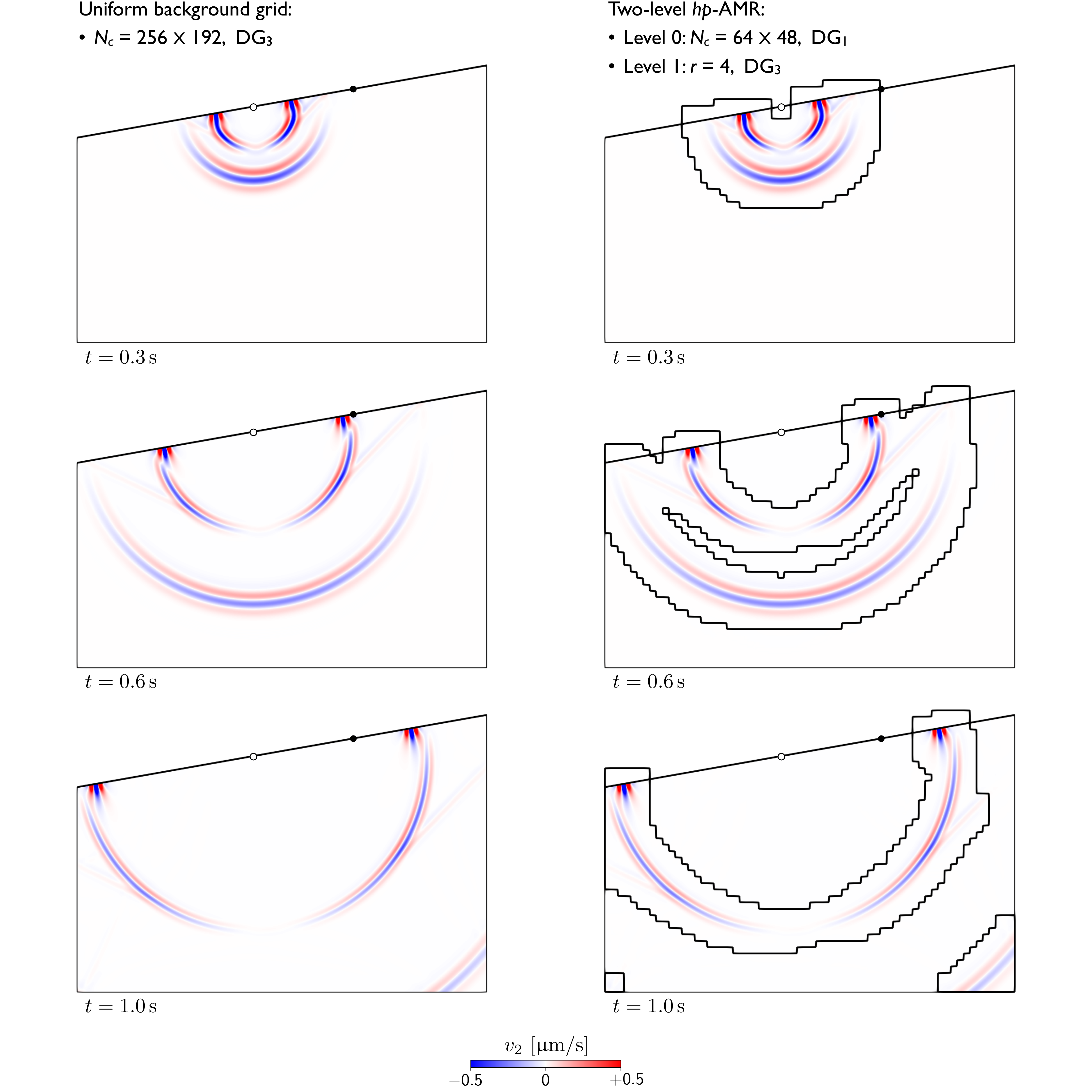}
\caption
{
    Snapshots of the vertical velocity component $v_2$ at the time instants $t = 0.3\, \mathrm{s}$, $0.6\, \mathrm{s}$ and $1.0\, \mathrm{s}$ for the 2D Lamb's problem.
    (left column) Uniform background grid.
    (right column) Two-level $hp$-AMR.
    In the right column, the stepped lines denote the boundary of the AMR level $\ell = 1$.
}
\label{fig:RESULTS - LAMB PROBLEM - SNAPSHOTS 2D}
\end{figure}

For the considered meshes, the values of the velocity components $v_1$ and $v_2$ measured at the receiver $\bm{r}$ of Fig.(\ref{fig:RESULTS - LAMB PROBLEM - GEOMETRY}a) are reported as functions of time in Figs.(\ref{fig:RESULTS - LAMB PROBLEM - SOLUTION 2D - VELOCITY}a) and (\ref{fig:RESULTS - LAMB PROBLEM - SOLUTION 2D - VELOCITY}b), respectively.
The plots show the expected convergence of the DG solution with respect to the number of mesh elements, and thus the mesh size, and the comparison between the numerical solution and the analytical solution \cite{berg1994analytical}, which is well recovered by the present scheme.

A clearer view of the wave structure generated by the concentrated force is displayed in Fig.(\ref{fig:RESULTS - LAMB PROBLEM - SNAPSHOTS 2D}), where the distribution of the velocity component $v_2$ is displayed at the time instants $t = 0.3\, \mathrm{s}$, $0.6\, \mathrm{s}$ and $1.0\, \mathrm{s}$.
The left column of the figure shows the results computed with the finest uniform grid, whereas the right column shows the results computed with the $hp$-AMR strategy; the same results are obtained with the two numerical setups.
Moreover, in either case, it is possible to distinguish the larger semicircle of the P-waves, which at $t = 1.0\, \mathrm{s}$ have almost left the domain of analysis, the smaller semicircle of the S-waves, which are travelling slower than the P-waves but have a similar spatial distribution, and the Rayleigh waves, which are travelling attached to the free surface at a speed that is slightly slower than that of the S-waves.

\begin{figure}
\centering
\begin{subfigure}{0.48\textwidth}
\includegraphics[width=\textwidth]{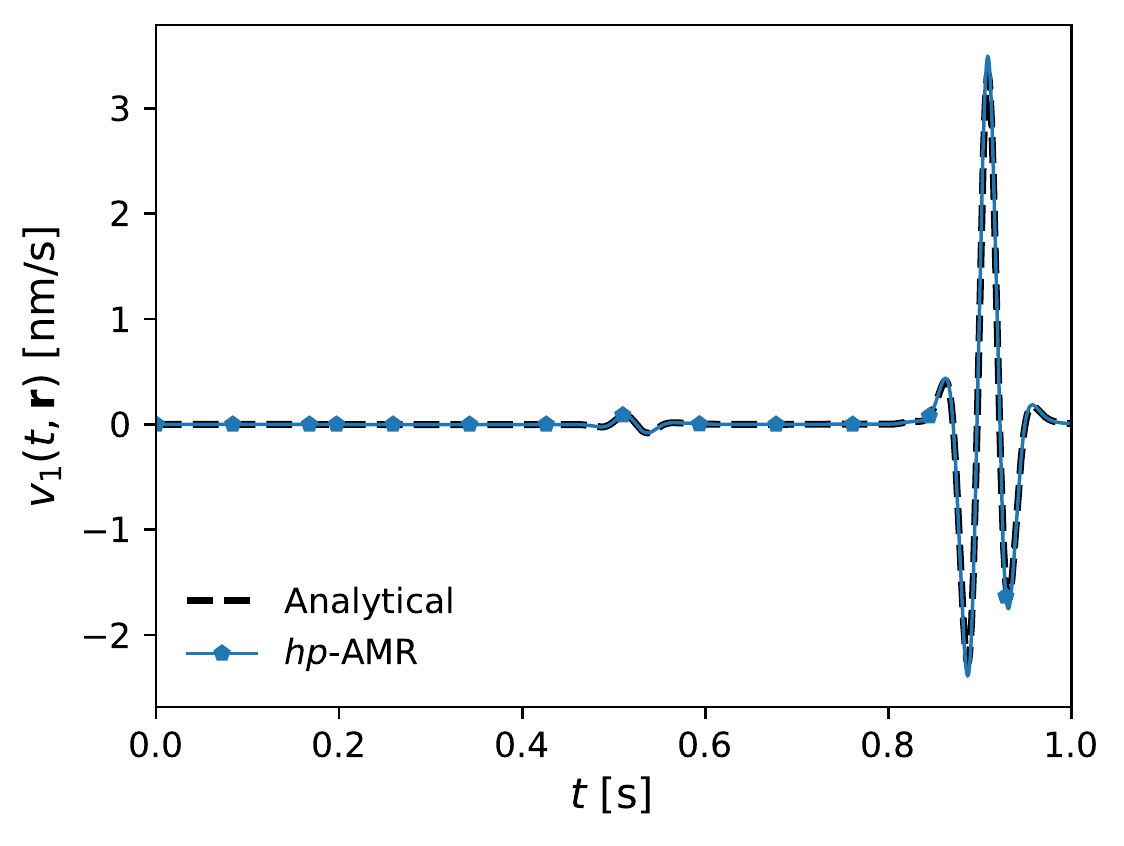}
\caption{}
\end{subfigure}
\
\begin{subfigure}{0.48\textwidth}
\includegraphics[width=\textwidth]{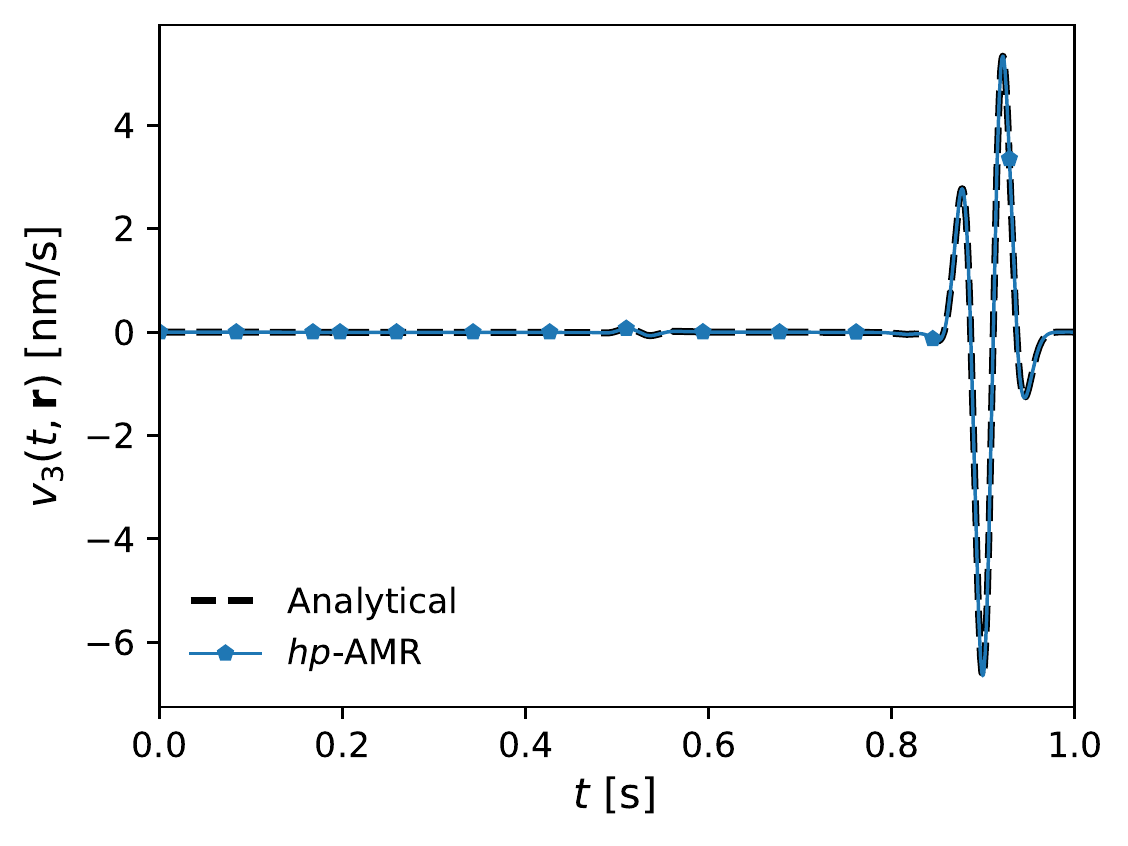}
\caption{}
\end{subfigure}
\caption
{
    Comparison between the analytical solution \cite{kausel2013lamb} and the DG solution for the 3D Lamb's problem in terms of (a) horizontal velocity component $v_1$ and (b) vertical velocity component $v_3$ evaluated at the surface receiver $\bm{r}$ denoted by the black dot in Fig.(\ref{fig:RESULTS - LAMB PROBLEM - GEOMETRY}b).
}
\label{fig:RESULTS - LAMB PROBLEM - SOLUTION 3D - VELOCITY}
\end{figure}

\begin{figure}
\centering
\includegraphics[width = \textwidth]{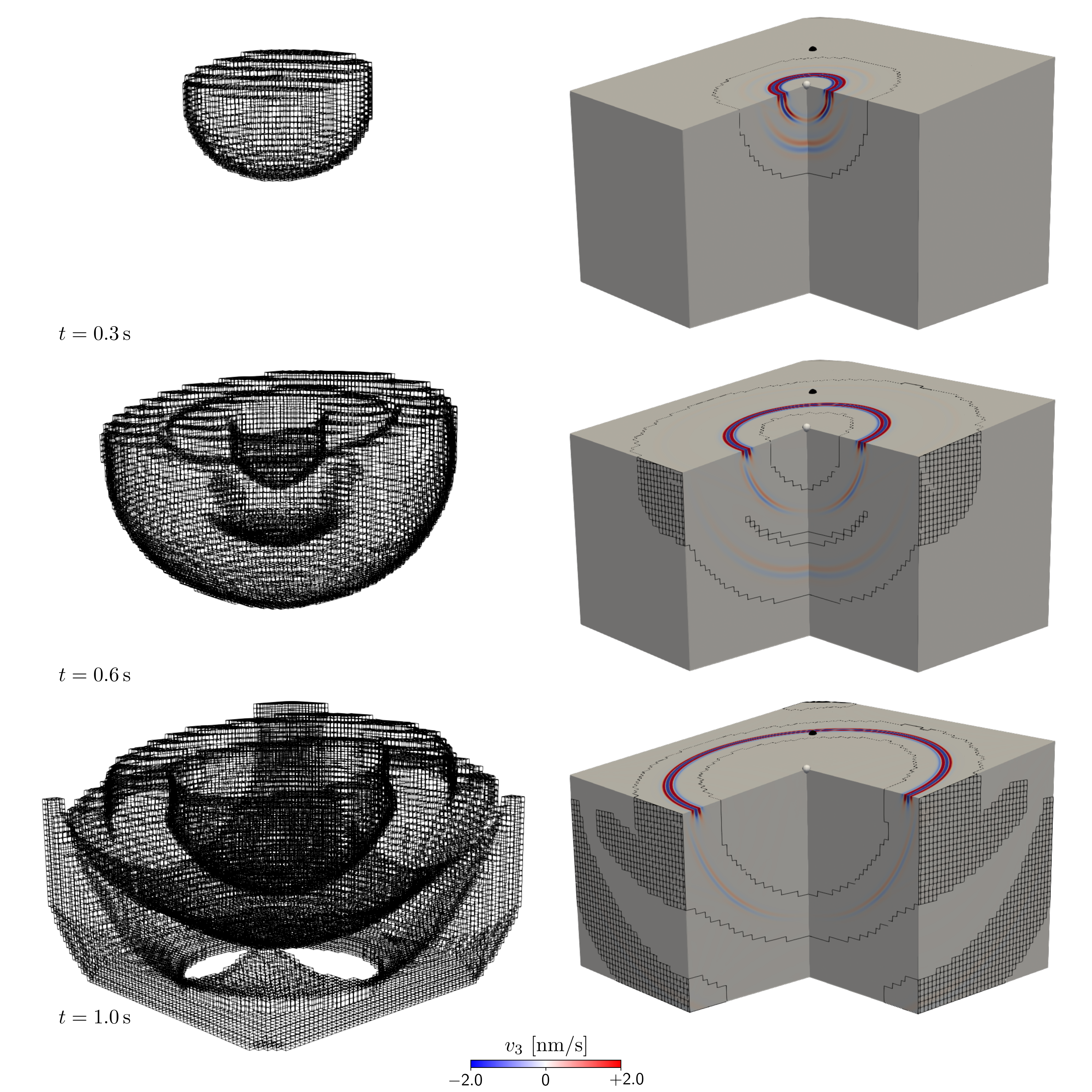}
\caption
{
    Snapshots of (left column) the boundary of the AMR level $\ell = 1$ and (right column) the vertical velocity component $v_3$ at the time instants $t = 0.3\, \mathrm{s}$, $0.6\, \mathrm{s}$ and $1.0\, \mathrm{s}$ for the 3D Lamb's problem.
    In the right column, the stepped lines denote the trace of AMR level $\ell = 1$ on the domain boundaries.
}
\label{fig:RESULTS - LAMB PROBLEM - SNAPSHOTS 3D}
\end{figure}

In 3D, guided by the 2D results, we consider only an implicitly-defined mesh associated with a two-level $hp$-AMR scheme.
The level $\ell = 0$ uses a $64\times 64\times 48$ grid and a DG$_1$ scheme, whereas level $\ell = 1$ uses a DG$_3$ scheme and is dynamically generated from level $\ell = 0$ via a refinement ratio $r = 4$.
Cell tagging is performed using Eq.(\ref{eq:tagging - Lamb problem}) where $E_0 = 10^{-18}\, \mathrm{J}$.

The values of the velocity components $v_1$ and $v_3$ measured at the receiver $\bm{r}$ of Fig.(\ref{fig:RESULTS - LAMB PROBLEM - GEOMETRY}b) are reported as functions of time in Figs.(\ref{fig:RESULTS - LAMB PROBLEM - SOLUTION 3D - VELOCITY}a) and (\ref{fig:RESULTS - LAMB PROBLEM - SOLUTION 3D - VELOCITY}b), respectively.
As in prior tests, the DG solution matches well with the exact solution.
Meanwhile, the left column of Fig.(\ref{fig:RESULTS - LAMB PROBLEM - SNAPSHOTS 3D}) shows the arrangement of the AMR level $\ell = 1$ at the time instants $t = 0.3\, \mathrm{s}$, $0.6\, \mathrm{s}$ and $1.0\, \mathrm{s}$, while the right column of Fig.(\ref{fig:RESULTS - LAMB PROBLEM - SNAPSHOTS 3D}) shows the distribution of the velocity component $v_3$ and the location of the P-, S- and Rayleigh waves at the same time instants.

\subsection{Single interface problem}

\begin{figure}
\centering
\includegraphics[width = \textwidth]{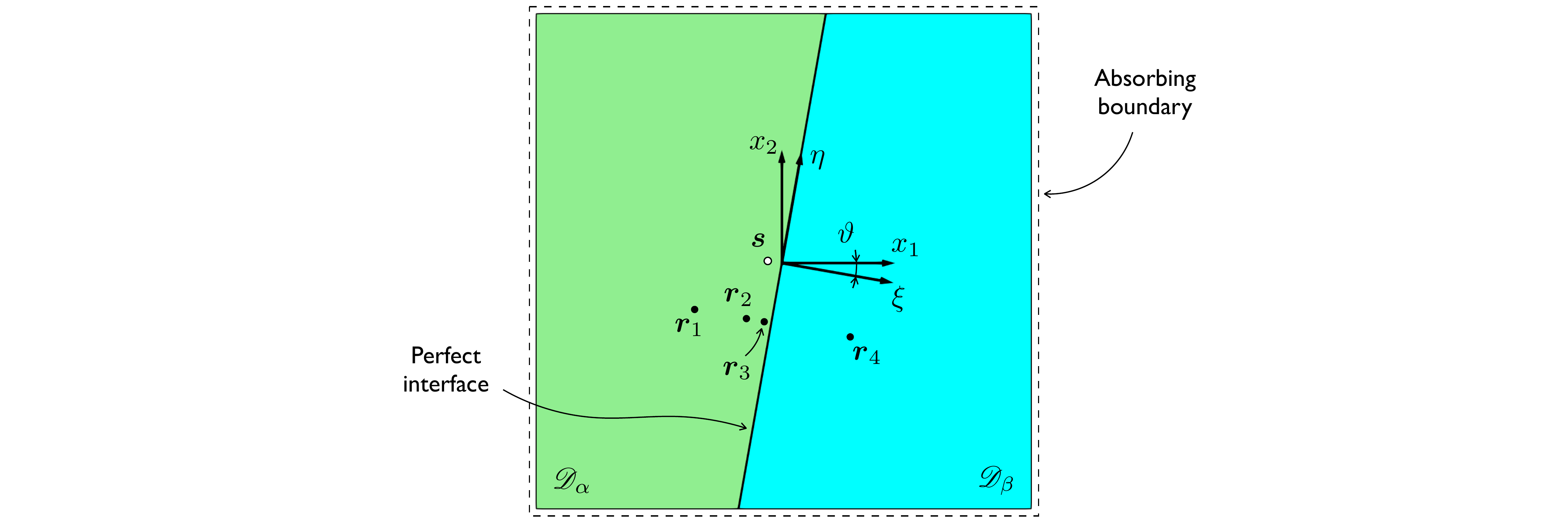}
\caption
{
    Geometry and boundary conditions for the single interface problem.
}
\label{fig:RESULTS - ONE INTERFACE PROBLEM - GEOMETRY 2D}
\end{figure}

Another classical problem in elastodynamics that has been modelled using different numerical methods, see e.g.~Refs.\cite{carcione1988wave,komatitsch2000simulation,de2007arbitrary}, regards the evaluation of the wave structure caused by a concentrated force acting in proximity of the interface between an isotropic solid and an orthotropic solid.
In the literature, the interface is typically aligned with the global reference system, i.e., the interface is perfectly horizontal or vertical; here, we instead place the interface on an angle in order to yield non-trivial implicitly-defined mesh geometry, similar to the case of the Lamb's problem discussed in Sec.(\ref{ssec:Lamb problem}).
The geometry is depicted in Fig.(\ref{fig:RESULTS - ONE INTERFACE PROBLEM - GEOMETRY 2D}) and is defined in the background square $\mathscr{R} = [-0.33\, \mathrm{m},0.33\, \mathrm{m}]^2$ by the level set function
\begin{equation}
\varphi(\bm{x}) = \cos(\vartheta)x_1+\sin(\vartheta)x_2
\end{equation}
where $\vartheta = 10^\circ$.
Figure (\ref{fig:RESULTS - ONE INTERFACE PROBLEM - GEOMETRY 2D}) also shows the location of the source point $\bm{s} = (s_1,s_2)$ where the concentrated force is applied, the location of four receiver points $\bm{r}_1$, $\bm{r}_2$, $\bm{r}_3$ and $\bm{r}_4$, where the mechanical signals are evaluated, and a local reference system that is aligned with the interface between the two phases.
In this local reference system, the coordinates of the source and receiver points are $s_\eta = -0.02\, \mathrm{m}$, $s_\xi = 0$, $r_{1\eta} = -0.105\, \mathrm{m}$, $r_{2\eta} = -0.035\, \mathrm{m}$, $r_{3\eta} = -0.01\, \mathrm{m}$, $r_{4\eta} = 0.105\, \mathrm{m}$, and $r_{1\xi} = r_{2\xi} = r_{3\xi} = r_{4\xi} = -0.08\, \mathrm{m}$.
The phase $\alpha$ and the phase $\beta$ are an orthorhombic solid and an isotropic solid, respectively, whose properties are $\rho_\alpha = \rho_\beta = 7100\, \mathrm{kg}/\mathrm{m}^3$ and
\begin{equation}
\bm{c}_\alpha = \left[\begin{array}{ccc}
165.0 & 50.0 & 0 \\
    & 62.0 & 0 \\
\mathrm{Sym.}&  & 39.6
\end{array}\right]\, \mathrm{GPa}
\quad\mathrm{and}\quad
\bm{c}_\beta = \left[\begin{array}{ccc}
165.0 & 85.8 & 0 \\
    & 165.0 & 0 \\
\mathrm{Sym.}&  & 39.6
\end{array}\right]\, \mathrm{GPa},
\end{equation}
where the elastic components are referred to the local reference system.
Perfect-interface conditions are prescribed on $\mathscr{L}_{\alpha,\beta}$ whereas absorbing boundary conditions are prescribed on $\mathscr{B}_\alpha$ and $\mathscr{B}_\beta$.
The final time of the simulation is $T = 100$ $\upmu$s.
Finally, the concentrated force is modelled by setting $\rho_\alpha\bm{b}_\alpha$ according to Eq.(\ref{eq:Lamb problem - external source}), where $\hat{\bm{u}}$ here is the unit vector parallel to $\mathscr{L}_{\alpha,\beta}$ and the parameters of the Ricker wavelet $R_w(t)$ are $a_1 = 10^{12}\, \mathrm{kg}/(\mathrm{m}^2\mathrm{s}^2)$, $a_2 = -\pi^2 f_c^2$, $f_c = 170.0\, \mathrm{kHz}$, $t_0 = 6\, {\upmu}$s.

We consider an implicitly-defined mesh generated from a $256^2$ uniform grid as well as one associated with a two-level $hp$-AMR scheme.
In the AMR case, level $\ell = 0$ uses a $64^2$ grid and a DG$_1$ scheme, whereas level $\ell = 1$ uses a DG$_3$ scheme and is dynamically generated from level $\ell = 0$ via a refinement ratio $r = 4$.
Similar to the Lamb's problem in 2D, the $hp$-AMR has the same effective resolution of the uniform mesh.
To evolve the cell tagging, the function $f_\mathrm{tag}$ introduced in Eq.(\ref{eq:tagging}) implements the following energy-based threshold
\begin{equation}\label{eq:tagging - single interface problem}
f_\mathrm{tag} \equiv \max\{E_\alpha^e,E_\beta^e\}-E_0,
\end{equation}
where $E_0 = 10^{4}\, \mathrm{J}$ and it is clear that $E_\alpha^e = 0$ if the element $\mathscr{D}_\alpha^e$ is empty.

\begin{figure}
\centering
\begin{subfigure}{0.48\textwidth}
\includegraphics[width=\textwidth]{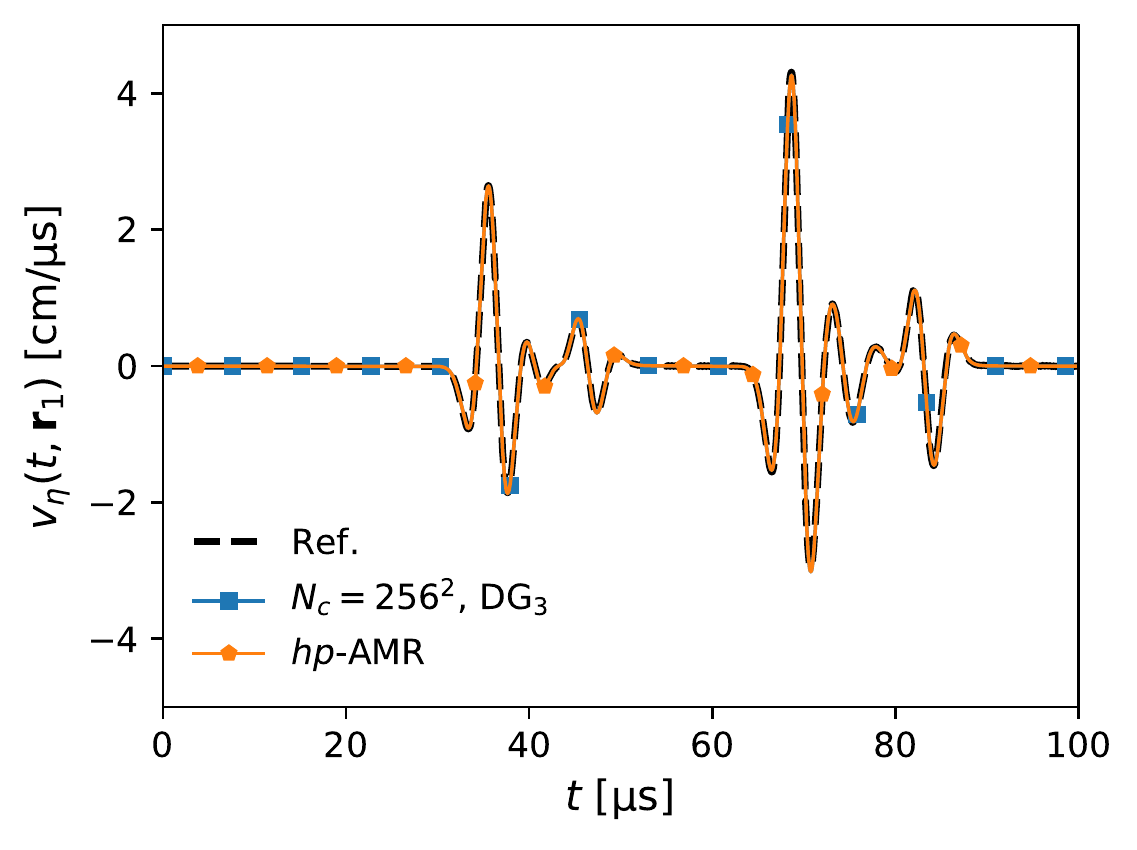}
\caption{}
\end{subfigure}
\
\begin{subfigure}{0.48\textwidth}
\includegraphics[width=\textwidth]{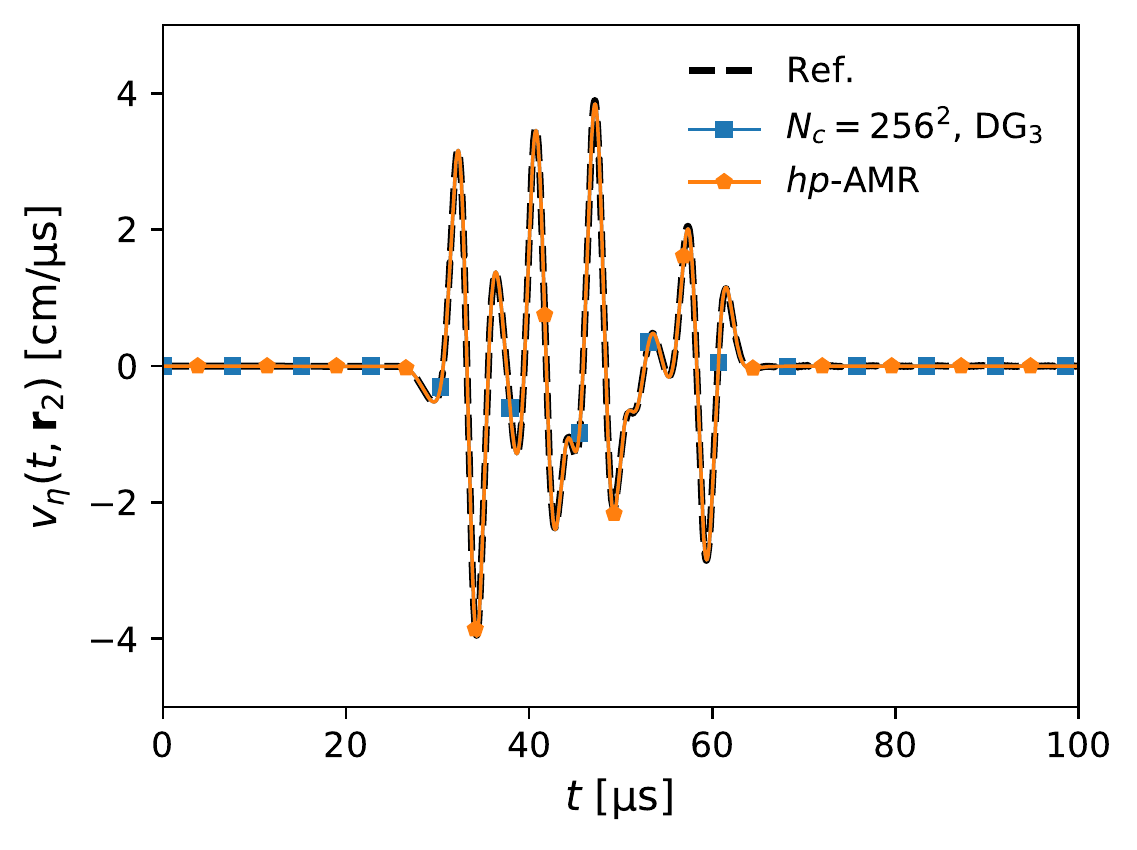}
\caption{}
\end{subfigure}
\\
\begin{subfigure}{0.48\textwidth}
\includegraphics[width=\textwidth]{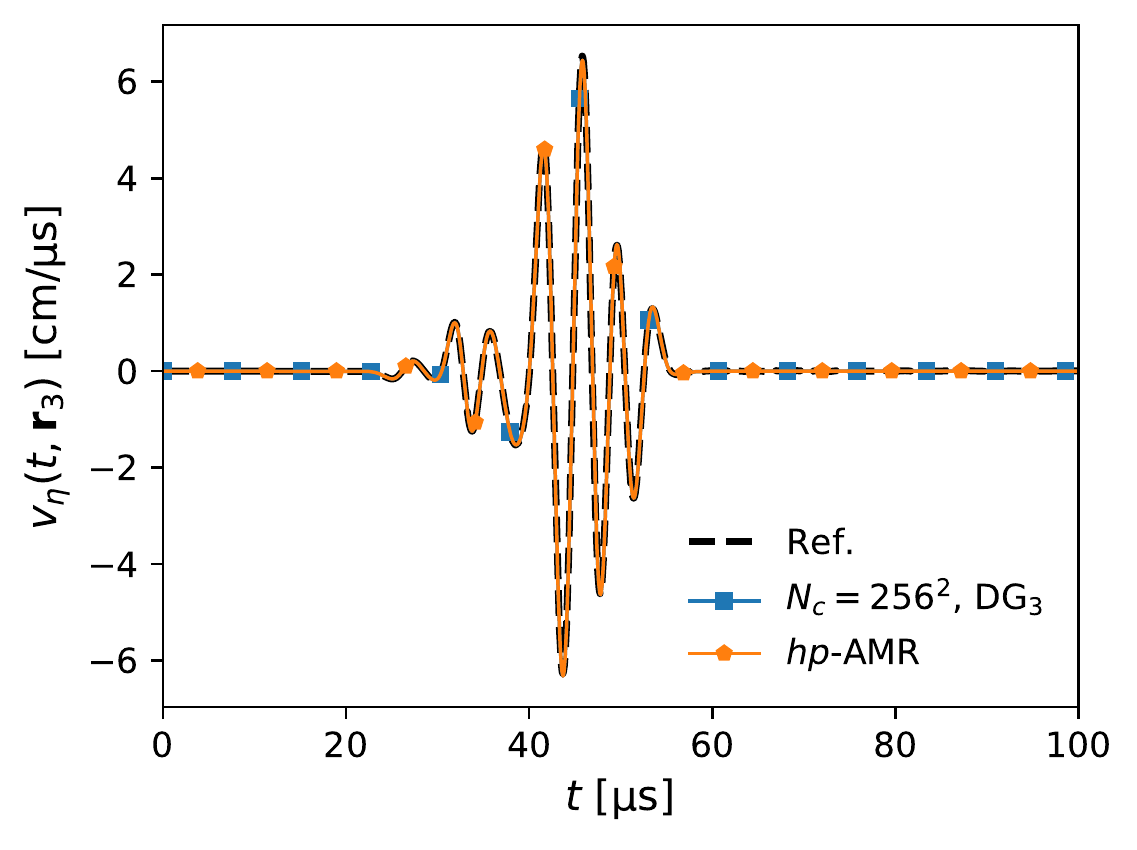}
\caption{}
\end{subfigure}
\
\begin{subfigure}{0.48\textwidth}
\includegraphics[width=\textwidth]{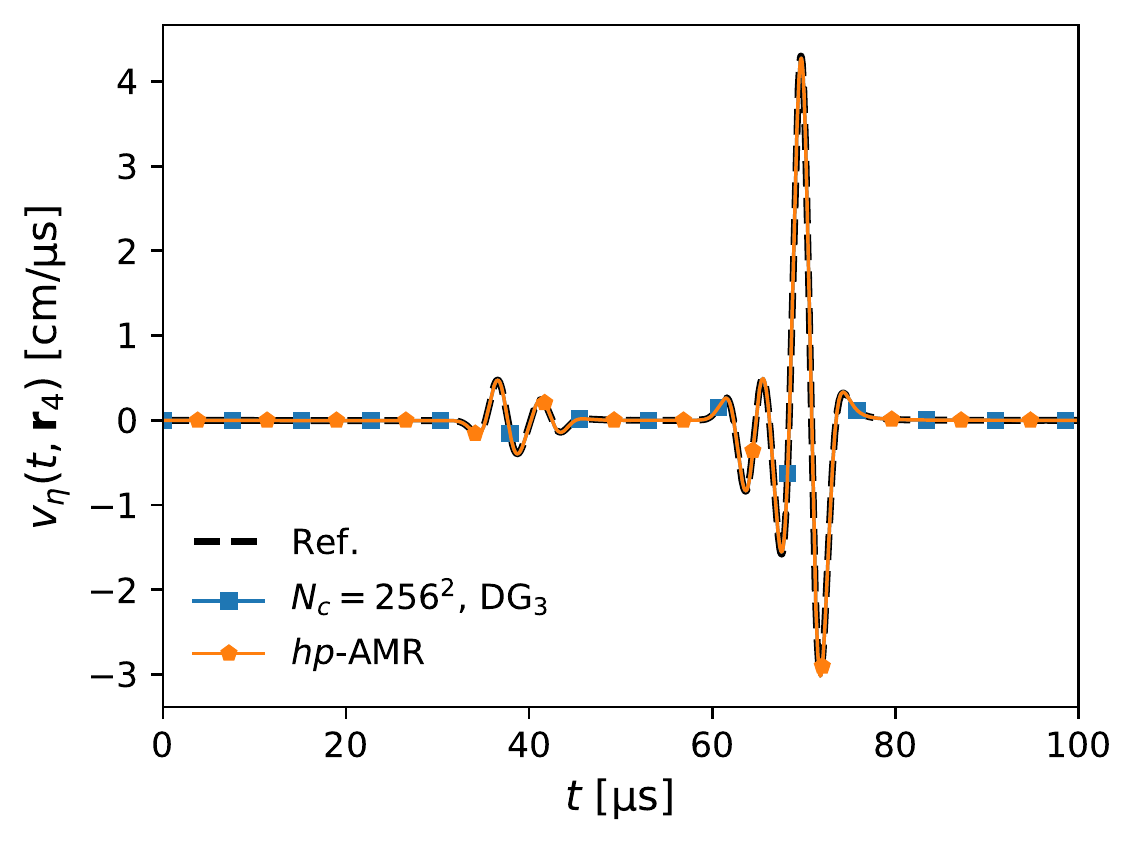}
\caption{}
\end{subfigure}
\caption
{
    Comparison between the reference solution \cite{komatitsch2000simulation} and the present formulation for the single interface problem in terms of the velocity component $v_\eta$ evaluated at the receivers (a) $\bm{r}_1$, (b) $\bm{r}_2$, (c) $\bm{r}_3$ and (d) $\bm{r}_4$ denoted by the black dots in Fig.(\ref{fig:RESULTS - ONE INTERFACE PROBLEM - GEOMETRY 2D}).
}
\label{fig:RESULTS - ONE INTERFACE PROBLEM - SOLUTION - VELOCITY}
\end{figure}

\begin{figure}
\centering
\includegraphics[width = \textwidth]{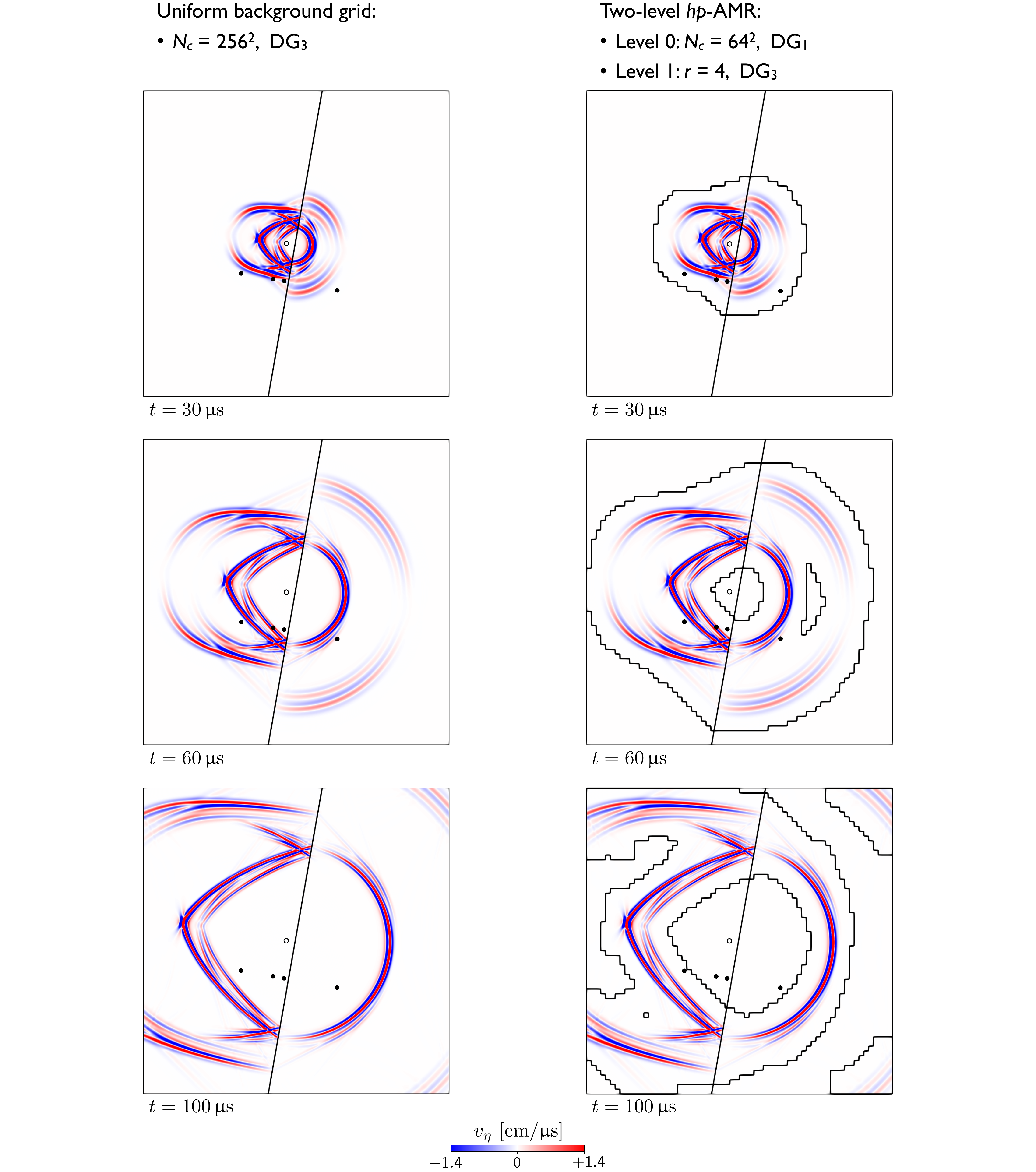}
\caption
{
    Snapshots of the velocity component $v_\eta$ at the time instants $t = 30\, \upmu$s, $60\, \upmu$s and $100\, \upmu$s for the single-interface interface problem.
    (left column) Uniform background grid.
    (right column) Two-level $hp$-AMR.
    In the right column, the stepped lines denote the boundaries between the $hp$-AMR levels.
}
\label{fig:RESULTS - ONE INTERFACE PROBLEM - SNAPSHOTS}
\end{figure}

Figures (\ref{fig:RESULTS - ONE INTERFACE PROBLEM - SOLUTION - VELOCITY}a) to (\ref{fig:RESULTS - ONE INTERFACE PROBLEM - SOLUTION - VELOCITY}d) report the velocity component ${v}_\eta$ at the receiver locations $\bm{r}_1$ to $\bm{r}_4$, respectively, and show that the results obtained with the uniform mesh and the results obtained with the $hp$-AMR are overlapping and match very well with the reference solution \cite{komatitsch2000simulation}.
Finally, the wave structure generated by the concentrated force is displayed in Fig.(\ref{fig:RESULTS - ONE INTERFACE PROBLEM - SOLUTION - VELOCITY}) at the time instants $t = 30\, \upmu$s, $60\, \upmu$s and $100\, \upmu$s in terms of the velocity component ${v}_\eta$.
From the figures, one can clearly observe the structure of the isotropic waves (characterized by semicircles) and the structure of the orthotropic waves, which propagate faster along the direction perpendicular to the interface and slower along the direction parallel to the interface.
Figure (\ref{fig:RESULTS - ONE INTERFACE PROBLEM - SOLUTION - VELOCITY}) also shows that the selected tagging criterion allows the $hp$-AMR scheme to reproduce the solution obtained with the uniform mesh.

\subsection{Structured solids}
\begin{figure}
\centering
\includegraphics[width = \textwidth]{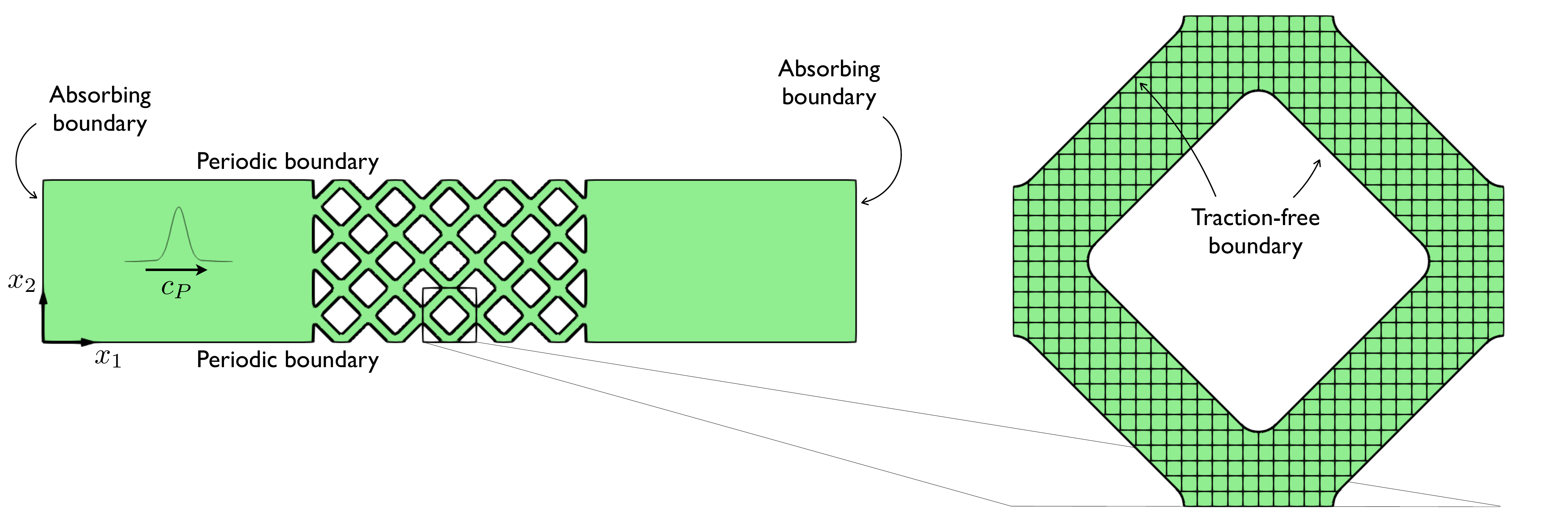}
\caption
{
    (Left) Geometry, initial conditions and boundary conditions for the 2D structured solid problem.
    (Right) Closeup on the geometry and the implicitly-defined mesh of the unit-cell.
}
\label{fig:RESULTS - PHONONIC CRYSTALS - GEOMETRY 2D}
\end{figure}

\begin{figure}
\centering
\includegraphics[width = \textwidth]{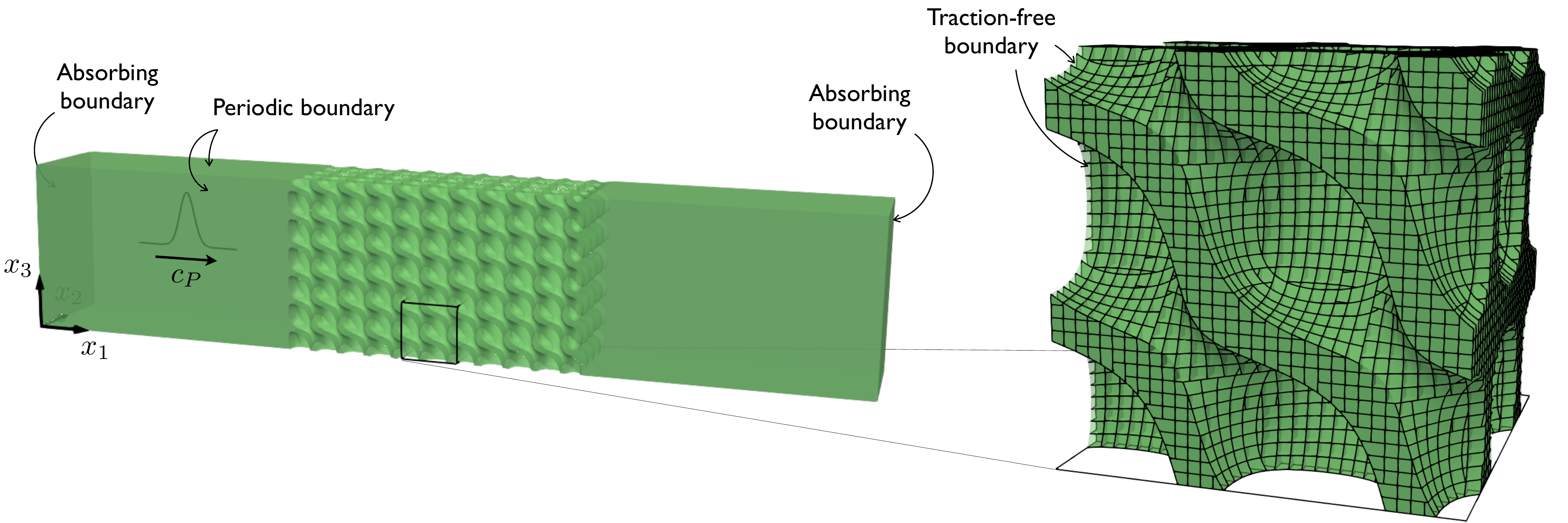}
\caption
{
    (Left) Geometry, initial conditions and boundary conditions for the 3D structured solid problem.
    (Right) Closeup on the geometry and the implicitly-defined mesh of the unit-cell.
}
\label{fig:RESULTS - PHONONIC CRYSTALS - GEOMETRY 3D}
\end{figure}

We conclude the numerical results by discussing the application of the present implicit-mesh DG framework to the analysis of an elastic wave propagating through a structured solid.
The solid is referred to as \emph{structured} because it is made of a periodically repeating structure, or unit cell, whose characteristic dimensions are smaller than the dimensions of the solid itself.
This scenario is common in the area of metamaterial design, where one is interested in obtaining non-conventional macroscopic elastodynamic properties for the solid, such as negative refractive index or negative effective density \cite{zhu2014negative,mokhtari2019emergence}, by changing the geometric features and/or by coupling dissimilar materials at the scale the unit cell.
Note that the aim of this section is not to provide an investigation of a structured solid in the context of metamaterials but to use the structured solid as an example of a (mildly) complex geometry where the elastodynamic problem can be solved with high-order accuracy in time and space by means of the present framework.
For this last set of tests, we also consider non-dimensional units.

The geometry for the considered 2D case is displayed in Fig.(\ref{fig:RESULTS - PHONONIC CRYSTALS - GEOMETRY 2D}).
The solid consists of two homogeneous ends and a central structured region, whose unit cell is displayed in the right-hand side of Fig.(\ref{fig:RESULTS - PHONONIC CRYSTALS - GEOMETRY 2D}).
To construct the geometry, we consider three level set functions:
\begin{subequations}\label{eq:structured solid - level set functions}
\begin{equation}
\varphi_1(\bm{x}) = |\sin(\pi (x_1-x_2-1/2))|-W,
\quad
\varphi_2(\bm{x}) = |\sin(\pi (x_1+x_2-1/2))|-W
\end{equation}
and
\begin{equation}
\varphi_3(\bm{x}) = (L-x_1)(x_1-2L),
\end{equation}
\end{subequations}
which are combined together to introduce a unique level set function as
\begin{equation}\label{eq:structured solid - level set function}
\varphi(\bm{x}) = 1
-\left[\max\left(0,1-\frac{\varphi_1(\bm{x})}{\delta_1}\right)\right]^{\delta_2}
-\left[\max\left(0,1-\frac{\varphi_2(\bm{x})}{\delta_1}\right)\right]^{\delta_2}
-\left[\max\left(0,1-\frac{\varphi_3(\bm{x})}{\delta_1}\right)\right]^{\delta_2}.
\end{equation}
In Eqs.(\ref{eq:structured solid - level set functions}), the functions $\varphi_1$ and $\varphi_2$ define the lattice structure of the unit cell, with $W$ controlling the width of the unit cell's struts, and the function $\varphi_3$ controls the location of the transition between the homogeneous ends and the structured region.
Meanwhile, Eq.(\ref{eq:structured solid - level set function}) provides a way to blend multiple level set functions where the positive parameters $\delta_1$ and $\delta_2$ control the sharpness of the transition among the level set functions; in particular, a low value of $\delta_1$ and a high value of $\delta_2$ make the transition sharper whereas a high value of $\delta_1$ and a low value of $\delta_2$ make the transition smoother.
Here, we use $L = 5$, $W = 0.4$, $\delta_1 = 0.2$ and $\delta_2 = 2$.

The geometry is periodic along the $x_2$ direction and, owing to its periodicity, the numerical problem is setup in the background rectangle $\mathscr{R} = [0,3L]\times[0,1]$ where the implicitly-defined mesh is generated from a background grid consisting of $480\times{32}$ cells.
The implicitly-defined mesh is partially shown in the right-hand side of Fig.(\ref{fig:RESULTS - PHONONIC CRYSTALS - GEOMETRY 2D}) for the unit cell.
Absorbing boundary conditions are prescribed at $x_1 = 0$ and $x_1 = 3L$, periodic boundary conditions are prescribed at $x_2 = 0$ and $x_2 = 1$, and zero-traction boundary conditions are prescribed on the zero-contour $\mathscr{L}_\alpha$ of the level set function $\varphi$.
Absorbing and periodic boundary conditions are chosen to minimize the reflection of the elastic waves from the background rectangle's boundaries while zero-traction boundary conditions are typical of single-phase metamaterials; different types of boundary condition could also be considered, especially if finite-size specimens are to be modelled.
The solid is assumed isotropic with density $\rho = 1$ and elastic constants defined by the velocity of the P-waves $c_P = 1$ and the velocity of the S-waves $c_S = 0.56$.
Initial conditions are prescribed as
\begin{equation}\label{eq:structured solid - ICs}
\bm{U}_{\alpha{0}} = \wtilde{\bm{U}}e^{-25(x_1-2.5)^2},
\end{equation}
where $\wtilde{\bm{U}}$ is the eigenvector solution of Eq.(\ref{eq:eigenvalue problem}) with $\bm{\kappa} = (1,0)$ and $\omega = c_P$; this initiates a wave that propagates from the homogeneous end of the solid to the structured region along the direction of the positive $x_1$ axis.

Similar to the 2D case, the considered 3D solid consists of two homogeneous ends and a central structured region as shown in Fig.(\ref{fig:RESULTS - PHONONIC CRYSTALS - GEOMETRY 3D}).
The unit cell is a Schwarz diamond \cite{maconachie2019slm} and is displayed in the right-hand side of Fig.(\ref{fig:RESULTS - PHONONIC CRYSTALS - GEOMETRY 2D}).
Consider the following functions:
\begin{subequations}
\begin{equation}
\varphi_1(\bm{x}) = \sin(\xi_1)\sin(\xi_2)\sin(\xi_3)+\sin(\xi_1)\cos(\xi_2)\cos(\xi_3)+\cos(\xi_1)\cos(\xi_2)\sin(\xi_3)+\cos(\xi_1)\cos(\xi_2)\cos(\xi_3),
\end{equation}
with $\xi_1 \equiv 2\pi x_1$, $\xi_2 \equiv 2 \pi x_2$ and $\xi_3 \equiv 2 \pi x_3$, and
\begin{equation}
\varphi_2(\bm{x}) = 
\begin{cases}
    x_1/L-1&\mathrm{if}~x_1 < 3L/2\\
    2-x_1/L&\mathrm{if}~x_1 \ge 3L/2
\end{cases},
\end{equation}
\end{subequations}
where $L = 5$.
Then, to construct the whole geometry, we define the following level set function
\begin{equation}
\varphi(\bm{x}) = w(\bm{x})\bigl(1/2-\varphi_1(\bm{x})\bigr)\bigl(1/2+\varphi_1(\bm{x})\bigr)+\bigl(1-w(\bm{x})\bigr)\varphi_2(\bm{x}),
\end{equation}
where $w(\bm{x})$ is a function controlling the transition between the homogeneous region and the structured region and is given by
\begin{equation}
w(\bm{x}) = \frac{1}{\pi}\left[\mathrm{atan}(50(x_1-L))-\mathrm{atan}(50(x_1-2L))\right].
\end{equation}

The numerical setup of the 3D problem is similar to the numerical setup of the 2D problem: the 3D geometry is periodic along the $x_2$ and the $x_3$ directions, and the implicitly-defined mesh is generated from a background grid consisting of $480\times{32}\times{32}$ cells in the background prism $\mathscr{R} = [0,3L]\times[0,1]\times[0,1]$.
The right-hand side of Fig.(\ref{fig:RESULTS - PHONONIC CRYSTALS - GEOMETRY 3D}) shows the resulting implicitly-defined mesh corresponding to the unit cell.
Absorbing boundary conditions are prescribed at $x_1 = 0$ and $x_1 = 3L$, periodic boundary conditions are prescribed at $x_2 = 0$ and $x_2 = 1$ and at $x_3 = 0$ and $x_3 = 1$, and zero-traction boundary conditions are prescribed on the zero-contour $\mathscr{L}_\alpha$ of the level set function $\varphi$.
The same elastic properties for the 2D case are employed, while initial conditions are given as in Eq.(\ref{eq:structured solid - ICs}) where $\wtilde{\bm{U}}$ is the eigenvector solution of the 3D version of Eq.(\ref{eq:eigenvalue problem}) with $\bm{\kappa} = (1,0,0)$ and $\omega = c_P$.

\begin{figure}
\centering
\includegraphics[width = \textwidth]{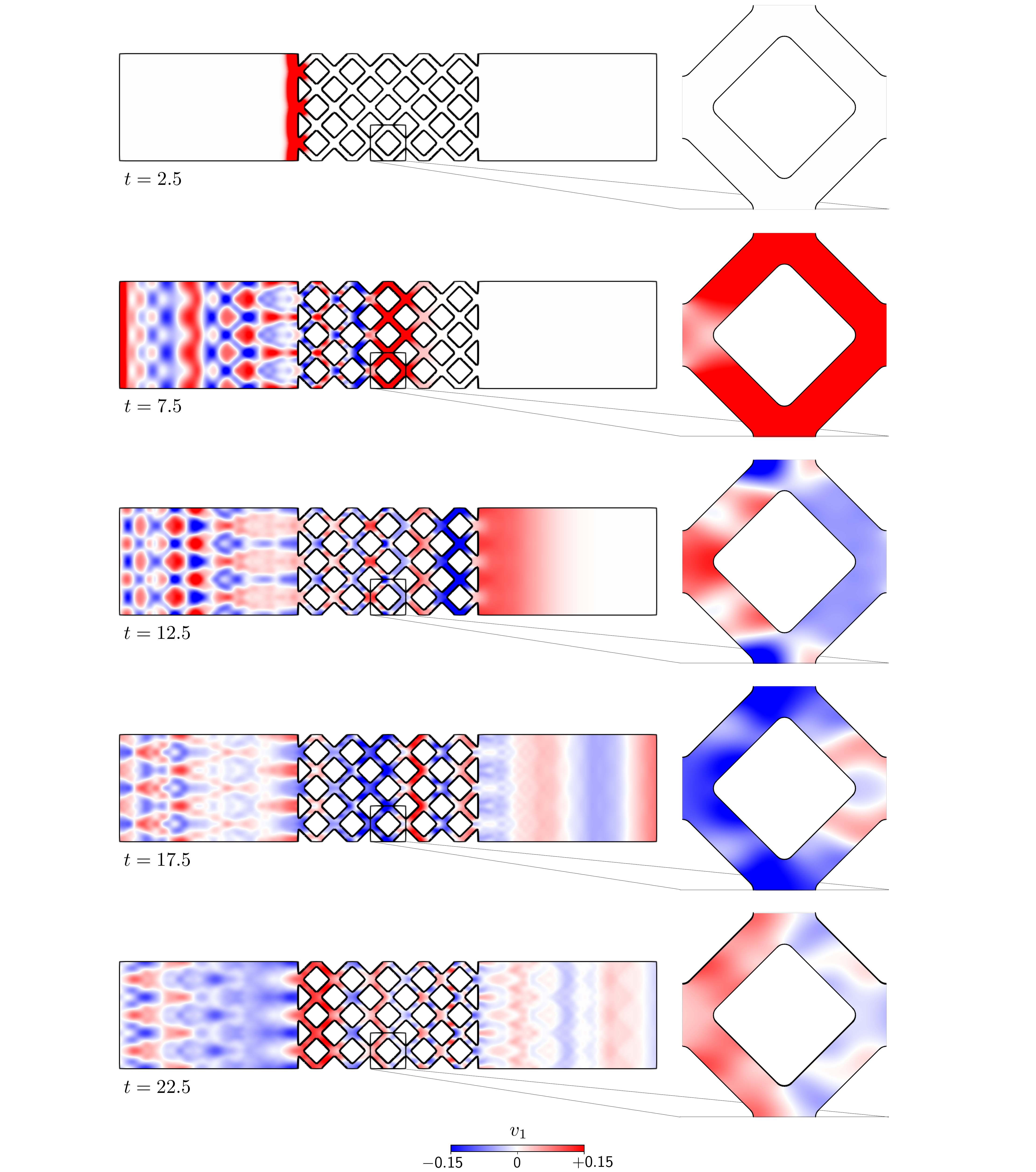}
\caption
{
    Snapshots of the velocity component $v_1$ at the time instants $t = 2.5$, $7.5$, $12.5$, $17.5$ and $22.5$ for the 2D structured solid problem.
}
\label{fig:RESULTS - PHONONIC CRYSTALS - SNAPSHOTS 2D}
\end{figure}

\begin{figure}
\centering
\includegraphics[width = \textwidth]{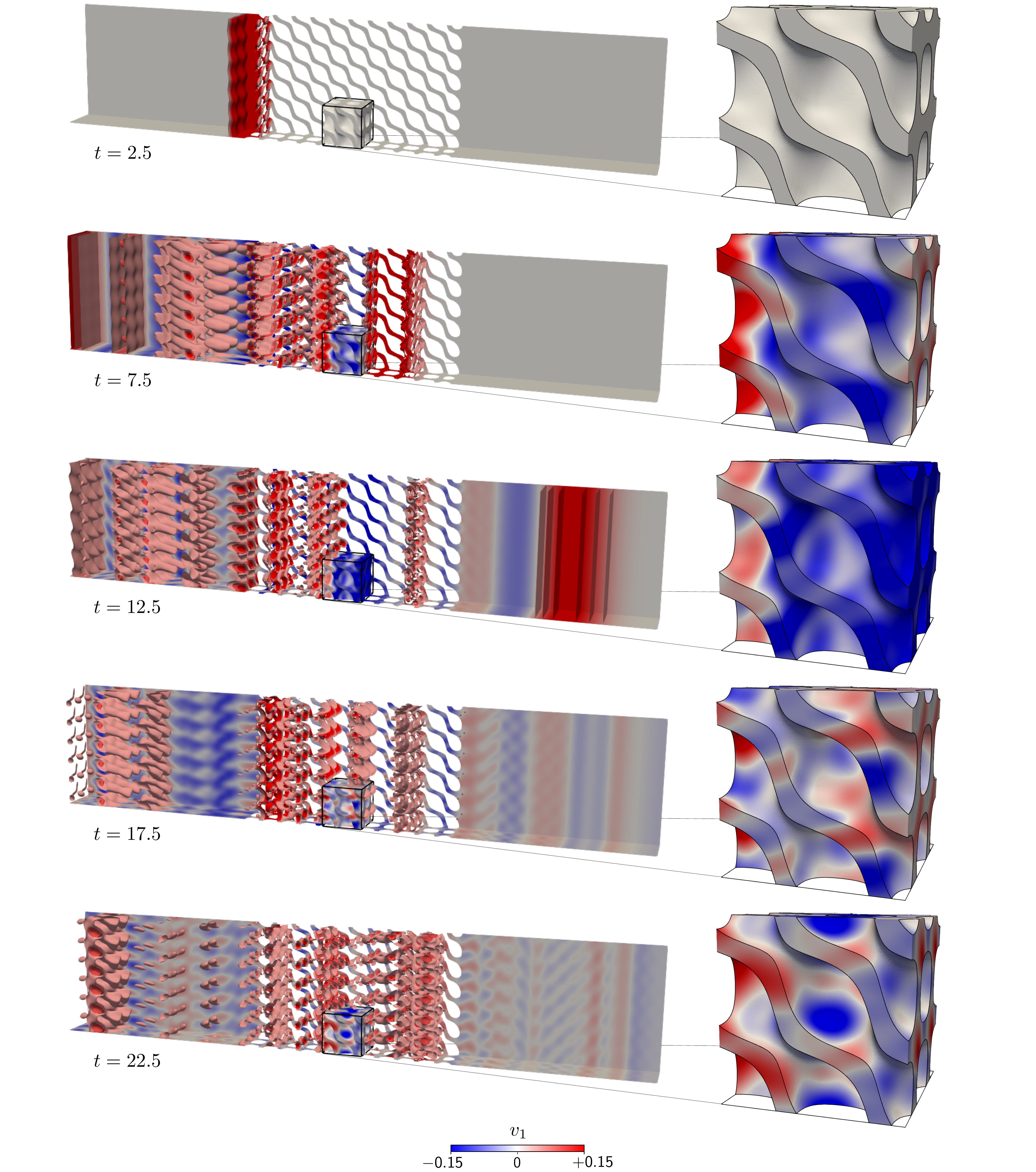}
\caption
{
    Snapshots of the velocity component $v_1$ at the time instants $t = 2.5$, $7.5$, $12.5$, $17.5$ and $22.5$ for the 3D structured solid problem.
}
\label{fig:RESULTS - PHONONIC CRYSTALS - SNAPSHOTS 3D}
\end{figure}

Figures (\ref{fig:RESULTS - PHONONIC CRYSTALS - SNAPSHOTS 2D}) and (\ref{fig:RESULTS - PHONONIC CRYSTALS - SNAPSHOTS 3D}) show a few snapshots of the velocity component $v_1$ at the time instants $t = 2.5$, $7.5$, $12.5$, $17.5$ and $22.5$ for the 2D setup and the 3D setup, respectively.
In both cases, it is possible to observe that part of the wave gets reflected by the structured region and part of it gets transmitted, whilst a complex distribution of the mechanical field is induced by the geometry of the structured solid.
Moreover, upon recalling that the length of the structured region is $L = 5$ and the wave travels at $c_P = 1$, it is interesting to notice that, between $t = 2.5$ and $t = 7.5$, the wave is not able to propagate from the beginning to the structured region to its end.
This means that, from a macroscopic viewpoint, the geometry of the structured region is responsible for slowing down the wave speed.
This is one example of several well-known features of metamaterials.
In the context of metamaterials, it would be also possible to analyse the frequency content of the reflected signal and of the transmitted signal and to investigate the stop-band properties of the structured solid, which might be considered as a filter for elastic waves.
However, these aspects are outside the scope of this paper and the application of the present framework to the analysis and design of metamaterials will be discussed elsewhere.
\section{Conclusions}\label{sec:CONCLUSIONS}
We have presented a discontinuous Galerkin framework for modeling wave propagation in single-phase and bi-phase elastic solid with complex geometries and general anisotropic constitutive behavior.
The framework belongs to embedded-boundary methods and is referred to as implicit-mesh DG method because it is based on the use of structured grids where the curved geometries are represented implicitly via level set functions and the domain discretization is generated by intersecting the level set functions with the grid cells, while a suitable cell-merging technique avoids the presence of overly small cut cells.
The novelty of the method regards the space discretization and, in particular, the use of high-order accurate quadrature rules for implicitly-defined domains and boundaries, which allow resolving the presence of the embedded geometries as well as enforcing boundary and interface conditions with high-order accuracy.

Various numerical tests have been considered and discussed, including several $hp$-convergence analyses in 2D and 3D and for single- and bi-phase solids, a few case studies involving 2D and 3D $hp$-AMR, as well as an application of the present method to the analysis of waves propagating in 2D and 3D structured solids.
The results demonstrate that the method achieves high-order accuracy in the maximum norm and is capable of dealing with implicitly-defined curved geometries, whilst taking advantage of the ease of generation and manipulation of structured grids.

The approach also offers several avenues of further research in the area of elastodynamics.
First, we recall that the present DG method has been employed to model waves in linear elastic solid with spatially constant elastic properties.
Therefore, a natural extension of the method would be to consider space-varying material properties so as to model functionally-graded materials with complex geometries; similarly, the method could be extended to account for non-linear elastic waves, see e.g.~Ref.\cite{bou2012nodal}.
Second, it is worth noting that the numerical tests feature a smooth geometry implicitly defined by a smooth level set function, including for the case of the structured solids wherein multiple level set functions were blended together to form a unique level set function; nevertheless, this does not represent a requirement (or limitation) for the present implicit-mesh DG framework, which can be used in combination with more complex geometry definitions, provided that the corresponding quadrature rules be available.
To this end, one possibility is to leverage the high-order accurate quadrature algorithms recently developed in Ref.\cite{saye2021high}, which can handle various kinds of complex geometry, such as intersecting/overlapping domains containing corners, junctions, tunnels, and multiple components, among other kinds of interfacial features.
Finally, we note that the simulations have been run using the functionalities for classic MPI parallelization implemented in \texttt{AMReX} \cite{zhang2019amrex}.
We have not discussed here the performance of the implementation as these aspects will be thoroughly investigated in future research, including a comprehensive scalability analysis involving also the use of modern accelerators, such as general-purpose graphical processing units.
\section*{Acknowledgements}
This research was supported in part by the Exascale Computing Project (17-SC-20-SC), a collaborative effort of the U.S.~Department of Energy Office of Science and the National Nuclear Security Administration, through U.S.~Department of Energy, Office of Science, Office of Advanced Scientific Computing Research, under contract DE-AC02-05CH11231.
It was also supported by the Applied Mathematics Program of the U.S.~Department of Energy Office of Advanced Scientific Computing Research under contract number DE-AC02-05CH11231.
Some computations used resources of the National Energy Research Scientific Computing Center, a DOE Office of Science User Facility supported by the Office of Science of the U.S. Department of Energy under Contract No. DE-AC02-05CH11231.


\bibliography{bibliography}

\end{document}